%\input TEXSHOP_macros_new.tex

%\input BoxedEPS.tex

  %%
 %%%%%%%%%%%%%%%%%%%%%%%%%%%%%%%%%%%%%%%%%%%%%%%%%%%%%%%%%%%%%
  %%
 %%%%%   BoxedEPS.tex FOR FIGURE INSERTS OF EPSF NORM  %%%%%
 %%%%%   (EPSF = Encapsulated PostScript File)
  %%
 %%%%%%%%%%%%%%%%%%%%%%%%%%%%%%%%%%%%%%%%%%%%%%%%%%%%%%%%%%%%%
  %%  
 %%%  AUTHOR: Laurent Siebenmann
  %%    lcs@matups.matups.fr
  %%  
 %%%  VERSIONS: Feb 1991 -- October, 1992
  %%  
 %%%  SOMMAIRE: BoxedEPS.tex d\'efinit des macro-commandes
  %%    qui permettent d'int\'egrer dans un document TeX des 
  %%    objets graphiques d\'ecrits par fichier de norme EPSF,
  %%    tout en accordant a chacun le statut d'une bo\^ite TeX ayant 
  %%    les bonnes dimensions.  La (seule!) contribution unique 
  %%    de ce fichier est de faire cela d'une fa{\c}con universelle.
  %%    C'est a dire de fa{\c}con \`a pouvoir commod\'ement 
  %%    servir avec tout pilote d'imprimante de norme 
  %%    PostScript --- malgr\'e l'absence d'une norme 
  %%    pour \special. 
  %%  
 %%%  POSTINGS: anonymous ftp 
  %%  ---  ftp 130.84.128.100 (alias rsovax.circe.fr); 
  %%  login: anonymous; password: <anything>; directory 
  %%  [anonymous.siebenmann].  This is the master copy in 1992.
  %%  
  %%  ---  ftp 129.69.1.12 (alias rusinfo.rus.uni-stuttgart.de);
  %%  login: anonymous; password: <anything>; 
  %%  directory hints .../tex/graphics/...
  %%  
 %%%% DOCUMENTATION:
  %%  --- see BoxedEPS.doc
  %%  
 %%%% ACTIVATION:
  %%  by a driver-by-driver protocol
  %%  see \SetTexturesEPSFSpecial 
  %%  and its companions below.
  %%  

 \ifx\MYUNDEFINED\BoxedEPSF
   \let\temp\relax
 \else
   \message{}
   \message{ !!! BoxedEPS %
         or BoxedArt macros already defined !!!}
   \let\temp\endinput
 \fi
  \temp
 
 \chardef\EPSFCatAt\the\catcode`\@
 \catcode`\@=11

 \chardef\C@tColon\the\catcode`\:
 \chardef\C@tSemicolon\the\catcode`\;
 \chardef\C@tQmark\the\catcode`\?
 \chardef\C@tEmark\the\catcode`\!
 \chardef\C@tDqt\the\catcode`\"

 \def\PunctOther@{\catcode`\:=12
   \catcode`\;=12 \catcode`\?=12 \catcode`\!=12 \catcode`\"=12}
 \PunctOther@

 %%temporarily suppress Plain's logging of allocations
 \let\wlog@ld\wlog 
 \def\wlog#1{\relax} 

 %% New for TOOLS
 \newif\ifIN@
 \newdimen\XShift@ \newdimen\YShift@ 
 \newtoks\Realtoks
 
 %%% New for Boxed EPSF
  %
 \newdimen\Wd@ \newdimen\Ht@
 \newdimen\Wd@@ \newdimen\Ht@@
 \newdimen\TT@
 \newdimen\LT@
 \newdimen\BT@
 \newdimen\RT@
 \newdimen\XSlide@ \newdimen\YSlide@ 
 \newdimen\TheScale  %% secretly scale in mils: 1pt= 1mil 
 \newdimen\FigScale  %% secretly scale in mils: 1pt= 1mil 
 \newdimen\ForcedDim@@

 \newtoks\EPSFDirectorytoks@
 \newtoks\EPSFNametoks@
 \newtoks\BdBoxtoks@
 \newtoks\LLXtoks@  %% useful info for Oz
 \newtoks\LLYtoks@

 \newif\ifNotIn@
 \newif\ifForcedDim@
 \newif\ifForceOn@
 \newif\ifForcedHeight@
 \newif\ifPSOrigin

 \newread\EPSFile@ 
 
 %%%% MESSAGES (separate macro needed for Europe)
  %%  
  \def\ms@g{\immediate\write16}

 %%%% WORD-PROCESSING MACROS
  %%
  %%% \IN@0#1@#2@ : Is 1st exp of #1 in 1st exp of #2 ??
   %% Answer in \ifIN@
 \newif\ifIN@\def\IN@{\expandafter\INN@\expandafter}
  \long\def\INN@0#1@#2@{\long\def\NI@##1#1##2##3\ENDNI@
    {\ifx\m@rker##2\IN@false\else\IN@true\fi}%
     \expandafter\NI@#2@@#1\m@rker\ENDNI@}
  \def\m@rker{\m@@rker}

  %%%  \SPLIT@0#1@#2@  :  Split 1st exp of #2 at 1st exp of #1
   %%  \Initialtoks@ , \Terminaltoks@ will contain pieces
  \newtoks\Initialtoks@  \newtoks\Terminaltoks@
  \def\SPLIT@{\expandafter\SPLITT@\expandafter}
  \def\SPLITT@0#1@#2@{\def\TTILPS@##1#1##2@{%
     \Initialtoks@{##1}\Terminaltoks@{##2}}\expandafter\TTILPS@#2@}

 %%%% MACROS TO TRIM  \ForeTrim@0#1@ and \Trim@0#1@  
   %% result appears in \Trimtoks@
   %% LIMITATION: assume no multiple spaces to trim

  \newtoks\Trimtoks@

  %%% \ForeTrim@0#1@ trims initial space of first erpansion of #1
   %% #1 of form \the\toks0 or \mymacro
 \def\ForeTrim@{\expandafter\ForeTrim@@\expandafter}
 \def\ForePrim@0 #1@{\Trimtoks@{#1}}
 \def\ForeTrim@@0#1@{\IN@0\m@rker. @\m@rker.#1@%
     \ifIN@\ForePrim@0#1@%
     \else\Trimtoks@\expandafter{#1}\fi}
   %%\m@rker expands here to \m@@rker since spot initial,
   %% so no confusuion with \m@rker

  %%% \Trim@0#1@ trims init and terminal spaces 
   %% Same syntax.
   %% Warns if internal spaces found.
   %% 
  \def\Trim@0#1@{%
      \ForeTrim@0#1@%
      \IN@0 @\the\Trimtoks@ @%
        \ifIN@ 
             \SPLIT@0 @\the\Trimtoks@ @\Trimtoks@\Initialtoks@
             \IN@0\the\Terminaltoks@ @ @%
                 \ifIN@
                 \else \Trimtoks@ {FigNameWithSpace}%
                 \fi
        \fi
      }

  %%%% MATH MACROS (provisional)
    %% use dimen registers for reals; unit 1pt
    %% (numerical dimension arguments OK unless contrary noted)

  %%%% One needs the point token seq (pt with cat 12) USES dimen 0
   \newtoks\pt@ks
   \def \getpt@ks 0.0#1@{\pt@ks{#1}}
   \dimen0=0pt\relax\expandafter\getpt@ks\the\dimen0@

   %%% Convert dimen to "decimal multiplier"% USES dimens 0,2
  \newtoks\Realtoks% the output!
  \def\Real#1{%
    \dimen2=#1%
      \SPLIT@0\the\pt@ks @\the\dimen2@%%  lop off the points
       \Realtoks=\Initialtoks@%\showthe\Realtoks
            }

   %%% Multiplication 
      % USES dimens 0,2,4,6; preserves args; output \Product
   \newdimen\Product
   \def\Mult#1#2{%
     \dimen4=#1\relax
     \dimen6=#2%
     \Real{\dimen4}%
     \Product=\the\Realtoks\dimen6%
        }

   %%% Inverse 
     % USES dimens 0; preserves arg; output \Inverse
 \newdimen\Inverse
 \newdimen\hmxdim@ \hmxdim@=8192pt%halfmaxdimen
 \def\Invert#1{%
  \Inverse=\hmxdim@
  \dimen0=#1%
  \divide\Inverse \dimen0%
  \multiply\Inverse 8}

 %%% \Rescale#1#2#3  % USES dimens 0,2,4,6
  %%  alters dimen register #1 by ratio #2/#3 
  %%  where #2,#3 can be raw dimensions OR dimen registers
   \def\Rescale#1#2#3{% Adequate accuracy. Can improve. 
              \divide #1 by 100\relax
              \dimen2=#3\divide\dimen2 by 100 \Invert{\dimen2}% 
              \Mult{#1}{#2}%
              \Mult\Product\Inverse 
              #1=\Product}

 %%% \Scale#1 scales dimen register #1 
   %  by dimen register real \TheScale; USES dimens 0
  \def\Scale#1{\dimen0=\TheScale %
      \divide #1 by  1280 %% 1280*5120*10=1000*2^16 
      \divide \dimen0 by 5120 % 
      \multiply#1 by \dimen0 
      \divide#1 by 10   %% max size of #1 about 32000/10 pt
     }
 
 %%% SCRUNCHING BOXES AND SHIFTING CONTENTS
  %% TeX has to do this in general
  %% since some drivers do not let 
  %% one do it readily using Postscript

 \newbox\scrunchbox

 %%% \Scrunched#1 puts #1 in an hbox
  %%    then in effect zeros the dimensions of this box
 \def\Scrunched#1{{\setbox\scrunchbox\hbox{#1}%
   \wd\scrunchbox=0pt
   \ht\scrunchbox=0pt
   \dp\scrunchbox=0pt
   \box\scrunchbox}}

  %%% \Shifted@#1 puts #1 in \hbox 
   %% then locates basepoint to bottom left corner
   %% then translates ink only by \XShift@,\YShift@
   %% with Postscript convention
   %% For simplicity use only on scrunched boxes
  %\newdimen\XShift@ 
  %\newdimen\YShift@ 
 \def\Shifted@#1{%
   \vbox {\kern-\YShift@
       \hbox {\kern\XShift@\hbox{#1}\kern-\XShift@}%
           \kern\YShift@}}

  %%% \cBoxedEPSF#1 the main macro
   %%  component macros are explained in order below

 \def\cBoxedEPSF#1{{\leavevmode 
    %% double brace for amstex \allign, \alligned, ...
   \ReadNameAndScale@{#1}%
   \SetEPSFSpec@
   \ReadEPSFile@ \ReadBdB@x  
   %% Calculations
     \TrimFigDims@ 
     \CalculateFigScale@  
     \ScaleFigDims@
     \SetInkShift@
   \hbox{$\mathsurround=0pt\relax
         \vcenter{\hbox{%
             \FrameSpider{\hskip-.4pt\vrule}%
             \vbox to \Ht@{\offinterlineskip\parindent=\z@%
                \FrameSpider{\vskip-.4pt\hrule}\vfil 
                \hbox to \Wd@{\hfil}%
                \vfil
                \InkShift@{\EPSFSpecial{\EPSFSpec@}{\FigSc@leReal}}%
             \FrameSpider{\hrule\vskip-.4pt}}%
         \FrameSpider{\vrule\hskip-.4pt}}}%
     $}%
    \CleanRegisters@ 
    \ms@g{ *** Box composed for the % 
         EPSF file \the\EPSFNametoks@}%
    }}
 
 \def\tBoxedEPSF#1{\setbox4\hbox{\cBoxedEPSF{#1}}%
     \setbox4\hbox{\raise -\ht4 \hbox{\box4}}%
     \box4
      }

 \def\bBoxedEPSF#1{\setbox4\hbox{\cBoxedEPSF{#1}}%
     \setbox4\hbox{\raise \dp4 \hbox{\box4}}%
     \box4
      }

  \let\BoxedEPSF\cBoxedEPSF% default setting

  %% Some compatibility with BoxedArt.tex
   %

  %% Some compatibility with Sweet-teX
   %
  \def\gLinefigure[#1scaled#2]_#3{%
        \BoxedEPSF{#3 scaled #2}}
    
  %% Some compatibility with Rokicki's dvips
   %

  \def\EPSFxsize{\afterassignment\ForceW@\ForcedDim@@}
      \def\ForceW@{\ForcedDim@true\ForcedHeight@false}
  
  \def\EPSFysize{\afterassignment\ForceH@\ForcedDim@@}
      \def\ForceH@{\ForcedDim@true\ForcedHeight@true}

  \def\EmulateRokicki{%
       \let\epsfbox\bBoxedEPSF \let\epsffile\bBoxedEPSF
       \let\epsfxsize\EPSFxsize \let\epsfysize\EPSFysize} 
 
 %%% \ReadNameAndScale@#1
  %
 \def\ReadNameAndScale@#1{\IN@0 scaled@#1@% DOUBLE BARRELED
   \ifIN@\ReadNameAndScale@@0#1@%
   \else \ReadNameAndScale@@0#1 scaled\DefaultMilScale @%
   \fi}
  
 \def\ReadNameAndScale@@0#1scaled#2@{% HELPER MACRO
    \let\OldBackslash@\\%
    \def\\{\OtherB@ckslash}%
    \edef\temp@{#1}%
    \Trim@0\temp@ @%
    \EPSFNametoks@\expandafter{\the\Trimtoks@ }%
    \FigScale=#2 pt%
    \let\\\OldBackslash@
    }
 
 \def\SetDefaultEPSFScale#1{%
      \global\def\DefaultMilScale{#1}}

 \SetDefaultEPSFScale{1000}

 %%% \ReadEPSFile@
  %
 \def \SetBogusBbox@{%
     \global\BdBoxtoks@{ BoundingBox:0 0 100 100 }%
     \global\def\BdBoxLine@{ BoundingBox:0 0 100 100 }%
     \ms@g{ !!! Will use placeholder !!!}%
     }

 {\catcode`\%=12\gdef\P@S@{%!}} %% %! min sign of PS file

 \def\ReadEPSFile@{%\show\EPSFSpec@%
     \openin\EPSFile@\EPSFSpec@
     \relax  %necessary to prevent precocious expansion of \ifeof
  \ifeof\EPSFile@
     \ms@g{}%
     \ms@g{ !!! EPS FILE \the\EPSFDirectorytoks@
       \the\EPSFNametoks@\space WAS NOT FOUND !!!}%
     \SetBogusBbox@
  \else%\fi
   \begingroup%%
   \catcode`\%=12\catcode`\:=12\catcode`\!=12
   \catcode`\G=14\catcode`\\=14\relax% 14 is comment
   \global\read\EPSFile@ to \BdBoxLine@%\show\BdBoxLine@
   \IN@0\P@S@ @\BdBoxLine@ @%
   \ifIN@ %% %! accepted as %!PS so do BdBox search!!
     \NotIn@true
     \loop   
       \ifeof\EPSFile@\NotIn@false 
         \ms@g{}%
         \ms@g{ !!! BoundingBox NOT FOUND IN %
            \the\EPSFDirectorytoks@\the\EPSFNametoks@\space!!! }%
         \SetBogusBbox@
       \else\global\read\EPSFile@ to \BdBoxLine@
       %\show\BdBoxLine@
       \fi
       \global\BdBoxtoks@\expandafter{\BdBoxLine@}%
       \IN@0BoundingBox:@\the\BdBoxtoks@ @%
       \ifIN@\NotIn@false\fi%
     \ifNotIn@\repeat
   \else
         \ms@g{}%
         \ms@g{ !!! \the\EPSFNametoks@\space not PS!\space !!!}%
         \SetBogusBbox@
   \fi
  \endgroup\relax
  \fi
  \closein\EPSFile@ 
   }

  %%% \ReadBdB@x
   % Rmk For simplicity 0 not used in syntax 
   %  of \ReadBdB@x@,  \ReadBdB@x@@ 
  \def\ReadBdB@x{% PART 0
   \expandafter\ReadBdB@x@\the\BdBoxtoks@ @}
  
  \def\ReadBdB@x@#1BoundingBox:#2@{% PART 1
    \ForeTrim@0#2@%
    \IN@0atend@\the\Trimtoks@ @%
       \ifIN@\Trimtoks@={0 0 100 100 }%
         \ms@g{}%
         \ms@g{ !!! BoundingBox not found in %
         \the\EPSFDirectorytoks@\the\EPSFNametoks@\space !!!}%
         \ms@g{ !!! It must not be at end of EPSF !!!}%
         \ms@g{ !!! Will use placeholder !!!}%
       \fi%% cf \SetBogusBbox@
    \expandafter\ReadBdB@x@@\the\Trimtoks@ @%
   }
    
  \def\ReadBdB@x@@#1 #2 #3 #4@{% PART 2
      \Wd@=#3bp\advance\Wd@ by -#1bp%
      \Ht@=#4bp\advance\Ht@ by-#2bp%
       \Wd@@=\Wd@ \Ht@@=\Ht@ %% useful info for Clark
       \LLXtoks@={#1}\LLYtoks@={#2}%% useful info for Oz
      \ifPSOrigin\XShift@=-#1bp\YShift@=-#2bp\fi 
     }

  %%% \SetEPSFDirectory 
   %
   \def\G@bbl@#1{}
   \bgroup
     \global\edef\OtherB@ckslash{\expandafter\G@bbl@\string\\}
   \egroup

  \def\SetEPSFDirectory{%  Part 1
           \bgroup\PunctOther@\relax
           \let\\\OtherB@ckslash
           \SetEPSFDirectory@}

 \def\SetEPSFDirectory@#1{% Part 2
    \edef\temp@{#1}%
    \Trim@0\temp@ @%  result in \Trimtoks@
    \global\toks1\expandafter{\the\Trimtoks@ }\relax
    \egroup
    \EPSFDirectorytoks@=\toks1
    }

  %%% \SetEPSFSpec@
 \def\SetEPSFSpec@{%
     \bgroup
     \let\\=\OtherB@ckslash
     \global\edef\EPSFSpec@{%
        \the\EPSFDirectorytoks@\the\EPSFNametoks@}%
     \global\edef\EPSFSpec@{\EPSFSpec@}%
     \egroup}

 %%% \TrimFigDims@ 
  % 
 \def\TrimTop#1{\advance\TT@ by #1}
 \def\TrimLeft#1{\advance\LT@ by #1}
 \def\TrimBottom#1{\advance\BT@ by #1}
 \def\TrimRight#1{\advance\RT@ by #1}

 \def\TrimBoundingBox#1{%
   \TrimTop{#1}%
   \TrimLeft{#1}%
   \TrimBottom{#1}%
   \TrimRight{#1}%
       }

 \def\TrimFigDims@{%
    \advance\Wd@ by -\LT@ 
    \advance\Wd@ by -\RT@ \RT@=\z@
    \advance\Ht@ by -\TT@ \TT@=\z@
    \advance\Ht@ by -\BT@ 
    }

 %%% \CalculateFigScale@
  %
  \def\ForceWidth#1{\ForcedDim@true
       \ForcedDim@@#1\ForcedHeight@false}
  
  \def\ForceHeight#1{\ForcedDim@true
       \ForcedDim@@=#1\ForcedHeight@true}

  \def\ForceOn{\ForceOn@true}
  \def\ForceOff{\ForceOn@false\ForcedDim@false}
  
  \def\CalculateFigScale@{%
            %Have default \FigScale or read \FigScale
     \ifForcedDim@\FigScale=1000pt% %% start afresh
           \ifForcedHeight@
                \Rescale\FigScale\ForcedDim@@\Ht@
           \else
                \Rescale\FigScale\ForcedDim@@\Wd@
           \fi
     \fi
     \Real{\FigScale}%
     \edef\FigSc@leReal{\the\Realtoks}%
     }
   
  \def\ScaleFigDims@{\TheScale=\FigScale
      \ifForcedDim@
           \ifForcedHeight@ \Ht@=\ForcedDim@@  \Scale\Wd@
           \else \Wd@=\ForcedDim@@ \Scale\Ht@
           \fi
      \else \Scale\Wd@\Scale\Ht@        
      \fi
      \ifForceOn@\relax\else\global\ForcedDim@false\fi
      \Scale\LT@\Scale\BT@  %%%\Scale\Wd@\Scale\Ht@
      \Scale\XShift@\Scale\YShift@
      }
      
  %%% \ShowReservedBoxes
   %%  shows (prints) corrected scaled and positioned
   %%  bounding boxes; for diagnostics
  %%% \HideReservedBoxes makes them invisible again
   %%
 \def\HideReservedBoxes{\global\def\FrameSpider##1{\null}}
 \def\ShowReservedBoxes{\global\def\FrameSpider##1{##1}}
 \let\HideDisplacementBoxes\HideReservedBoxes  %% some synonyms
 \let\ShowDisplacementBoxes\ShowReservedBoxes
 \let\HideFigureFrames\HideReservedBoxes
 \let\ShowFigureFrames\ShowReservedBoxes
  \ShowDisplacementBoxes
 
  %%% \hSlide#1, \vSlide#1
   %%
 \def\hSlide#1{\advance\XSlide@ by #1}
 \def\vSlide#1{\advance\YSlide@ by #1}
 
  %%% \SetInkShift@, \InkShift@#1
   %%
  \def\SetInkShift@{%
            \advance\XShift@ by -\LT@
            \advance\XShift@ by \XSlide@
            \advance\YShift@ by -\BT@
            \advance\YShift@ by -\YSlide@
             }
  \def\InkShift@#1{\Shifted@{\Scrunched{#1}}}
 
  %%% \CleanRegisters@
   %
  \def\CleanRegisters@{%
      \globaldefs=1\relax
        \XShift@=\z@\YShift@=\z@\XSlide@=\z@\YSlide@=\z@
        \TT@=\z@\LT@=\z@\BT@=\z@\RT@=\z@
      \globaldefs=0\relax}

 %%% Special syntax for several drivers. The macros 
  %% \SetTexturesEPSFSpecial  %% Textures 
  %% \SetUnixCoopEPSFSpecial %% dvi2ps early unix 
  %% \SetBechtolsheimDVI2PSEPSFSpecial and 
  %% \SetBechtolsheimDVITPSEPSFSpecial %% by S.P.Bechtolsheim
  %% \SetLisEPSFSpecial %% dvi2ps by Tony Lis
  %% \SetRokickiEPSFSpecial  %% dvips by Tom Rokicki
  %%  --- also for DVIReader, in DirectTeX by W. Ricken
  %% \SetOzTeXEPSFSpecial  %% OzTeX (>=1.42) by Andrew Trevorrow
  %% \SetPSprintEPSFSpecial %% PSprint by Andrew Trevorrow
  %%  --- also for OzTeX versions <= 1.41 !!
  %% \SetArborEPSFSpecial  %% ArborTeX DVILASER/PS
  %% \SetClarkEPSFSpecial %% dvitops by James Clark
  %% \SetDVIPSoneEPSFSpecial %% DVIPSONE of Y&Y 
  %% \SetBeebeEPSFSpecial %% DVIALW by N. Beebe
  %% \SetNorthlakeEPSFSpecial %% Northlake Software
  %% \SetStandardEPSFSpecial %% Nonexistant: Placebo below
  %% Many drivers supported roughly
  %% by (re-)defining the macro \EPSFSpecial#1#2, where
  %% #1 = EPS file pathname (use \\ for the letter backslash)
  %% #2 = scale in mils 
  %% Be wary of using strange characters in pathnames!
 
 %% Textures, Blue Sky Research, Barry Smith
 \def\SetTexturesEPSFSpecial{\PSOriginfalse%\PSOrigintrue
  \gdef\EPSFSpecial##1##2{\relax
    \edef\specialthis{##2}%
    \SPLIT@0.@\specialthis.@\relax
    \special{illustration ##1 scaled
                        \the\Initialtoks@}}}
 
  %% Unix : dvi2ps by:  Mark Senn, Stephan  Bechtolsheim,  
   % Bob  Brown, Richard, Furuta, James Schaad, Robert  Wells, 
   % Norm Hutchinson, Neal Holt, Scott Jones, Howard Trickey.
   % Introduced by B. Horn <bkph@ai.mit.edu>
  \def\SetUnixCoopEPSFSpecial{\PSOrigintrue % Please test!
   \gdef\EPSFSpecial##1##2{%
      \dimen4=##2pt% convert real to dimen
      \divide\dimen4 by 1000\relax
      \Real{\dimen4}%dimens 0,2 used here
      \edef\Aux@{\the\Realtoks}%  
      %%convert dimen to real
      \includegraphics{##1\space}}}

  %% dvi2ps and dvitps by S.P. Bechtolsheim,
   % Introduced by B. Horn <bkph@ai.mit.edu> and Carl.M.Jones, 
   % testing by R. Evans <Robert@cm.cardiff.ac.uk>
   % Note that a prolog file psfig.pro
   % specific to the driver should be available.
  \def\SetBechtolsheimEPSFSpecial@{%% tool macro only
   \PSOrigintrue
   \special{\DriverTag@ Include0 "psfig.pro"}%
   \gdef\EPSFSpecial##1##2{%
      \dimen4=##2pt %% convert real to dimen
      \divide\dimen4 by 1000\relax
      \Real{\dimen4} %% dimens 0,2 used here
      \edef\Aux@{\the\Realtoks}%% convert dimen to real
      \special{\DriverTag@ Literal "10 10 0 0 10 10 startTexFig
           \the\mag\space 1000 div 3.25 neg mul 
           \the\mag\space 1000 div .25 neg mul translate %% correction
           \the\mag\space 1000 div \Aux@\space mul 
           \the\mag\space 1000 div \Aux@\space mul scale "}%
      \special{\DriverTag@ Include1 "##1"}%
      \special{\DriverTag@ Literal "endTexFig "}%
        }}

  %% dvi2ps and dvitps by S.P. Bechtolsheim,
   % Introduced by B. Horn <bkph@ai.mit.edu> and Carl.M.Jones, 
   % testing by R. Evans <Robert@cm.cardiff.ac.uk>
   % Note that a prolog file psfig.pro
   % specific to the driver should be available.
  \def\SetBechtolsheimEPSFSpecial@{%% tool macro only
   \PSOrigintrue
   \special{\DriverTag@ Include0 "psfig.pro"}%
   \gdef\EPSFSpecial##1##2{%
      \dimen4=##2pt %% convert real to dimen
      \divide\dimen4 by 1000\relax
      \Real{\dimen4} %% dimens 0,2 used here
      \edef\Aux@{\the\Realtoks}%% convert dimen to real
      \special{\DriverTag@ Literal "10 10 0 0 10 10 startTexFig
           \the\mag\space 1000 div 
           dup 3.25 neg mul 2 index .25 neg mul translate %% correction line
           \Aux@\space mul dup scale "}%
      \special{\DriverTag@ Include1 "##1"}%
      \special{\DriverTag@ Literal "endTexFig "}%
        }}

  \def\SetBechtolsheimDVITPSEPSFSpecial{\def\DriverTag@{dvitps: }%
      \SetBechtolsheimEPSFSpecial@}

  \def\SetBechtolsheimDVI2PSEPSFSSpecial{\def\DriverTag@{DVI2PS: }%
      \SetBechtolsheimEPSFSpecial@}

  %% dvi2ps by Tony Lis,
   % implantations? ; dates?; availability?
   % Introduced by B. Horn <bkph@ai.mit.edu>
  \def\SetLisEPSFSpecial{\PSOrigintrue 
   \gdef\EPSFSpecial##1##2{%
      \dimen4=##2pt% convert real to dimen
      \divide\dimen4 by 1000\relax
      \Real{\dimen4}% dimens 0,2 used here
      \edef\Aux@{\the\Realtoks}%  
      %%convert dimen to real
      \special{pstext="10 10 0 0 10 10 startTexFig\space
           \the\mag\space 1000 div \Aux@\space mul 
           \the\mag\space 1000 div \Aux@\space mul scale"}%
      \includegraphics{##1}%
      \special{pstext=endTexFig}%
        }}

  %% dvips by Tom Rokicki; free driver in portable C 
   % Introduced by W.D. Neumann <neumann@mps.ohio-state.edu>
  \def\SetRokickiEPSFSpecial{\PSOrigintrue 
   \gdef\EPSFSpecial##1##2{%
      \dimen4=##2pt% convert real to dimen
      \divide\dimen4 by 10\relax
      \Real{\dimen4}% dimens 0,2 used here
      \edef\Aux@{\the\Realtoks}%  
      %%convert dimen to real
      \includegraphics{##1}}}

  \def\SetInlineRokickiEPSFSpecial{\PSOrigintrue 
   \gdef\EPSFSpecial##1##2{%
      \dimen4=##2pt% convert real to dimen
      \divide\dimen4 by 1000\relax
      \Real{\dimen4}% dimens 0,2 used here
      \edef\Aux@{\the\Realtoks}%  
      %%convert dimen to real
      \special{ps::[begin] 10 10 0 0 10 10 startTexFig\space
           \the\mag\space 1000 div \Aux@\space mul 
           \the\mag\space 1000 div \Aux@\space mul scale}%
      \special{ps: plotfile ##1}%
      \special{ps::[end] endTexFig}%
        }}

 %%%  OzTeX (versions 1.42 and later), by Andrew Trevorrow
 %%%  (for earlier versions see PSprint below!!)
 %%  complete public domain TeX for Macintosh
 %%  Send 10 UNFORMATTED 800K disks 
 %%  with return postage to
 %%  Peter Abbott, Computing Service, 
 %%  Aston University, Aston Triangle, Birmingham B4 7ET
 %%  Posting: ftp   midway.uchicago.edu
 %%  Nota: Version 1.42 may give
 %%  spurious "offpage" error notices on printing.
 %%  Nota: Support for MacPaint files not here yet.
 \def\SetOzTeXEPSFSpecial{\PSOrigintrue
 \gdef\EPSFSpecial##1##2{%
 \dimen4=##2pt%% convert real to dimen
 \divide\dimen4 by 1000\relax
 \Real{\dimen4}%% dimens 0,2 used here
 \edef\Aux@{\the\Realtoks}%% convert dimen to real
 \special{epsf=\string"##1\string"\space scale=\Aux@}%
 }} 

 %% PSprint, by AndrewTrevorrow for VaX VMS
 %% and OzTeX versions <= 1.41  
  % tested 2-91 by Max Calviani <ISICA@ASTRPD.infn.it>
  \def\SetPSprintEPSFSpecial{\PSOriginFALSE % artifice; see below
   \gdef\EPSFSpecial##1##2{%note order
     \special{##1\space 
       ##2 1000 div \the\mag\space 1000 div mul
       ##2 1000 div \the\mag\space 1000 div mul scale
       \the\LLXtoks@\space neg \the\LLYtoks@\space neg translate
       }}}

 %% DVILASER/PS driver originally written by David Fuchs
  % marketed and supported by ArborTeXt  535 W. William St.
  % Suite 300, Ann Arbor, MI 48103, U.S.A
  % (313) 996-3566 (313) 996-3573
  % help@arbortext.com, Andrew Dobrowolski
 \def\SetArborEPSFSpecial{\PSOriginfalse % check!
   \gdef\EPSFSpecial##1##2{%
     \edef\specialthis{##2}%
     \SPLIT@0.@\specialthis.@\relax % suppress decimals (nec!)
     \special{ps: epsfile ##1\space \the\Initialtoks@}}}

 %% dvitops, (c) James Clark <jjc@jclark.uucp>
  % public domain; distributed by UK TeX Archive
  % computers: unix, msdos, vms, primos and vm/cms,
  % introduced by S. Ratz <spqr@uk.ac.southampton.ecs>
 \def\SetClarkEPSFSpecial{\PSOriginfalse % please test!
   \gdef\EPSFSpecial##1##2{%
     \Rescale {\Wd@@}{##2pt}{1000pt}%
     \Rescale {\Ht@@}{##2pt}{1000pt}%
     \special{dvitops: import 
           ##1\space\the\Wd@@\space\the\Ht@@}}}

 %% DVIPSONE, for PC compatibles
  % Y&Y, 106 Indian Hill, Carlisle MA 01741, USA
  % (508) 371-3286
  % (introduced by B. Horn <bkph@ai.mit.edu>)
  \let\SetDVIPSONEEPSFSpecial\SetUnixCoopEPSFSpecial
  \let\SetDVIPSoneEPSFSpecial\SetUnixCoopEPSFSpecial

 %% DVIALW by N. Beebe, public domain 
  % DVI Driver Distribution, Center for Scientific Computing,
  % Department of Mathematics, 220 South Physics Building,
  % University of Utah, Salt Lake City, UT 84112, USA
  % (introduced by B. Horn <bkph@ai.mit.edu>)
  % Proposed standard; see TUGboat article 1993.
  \def\SetBeebeEPSFSpecial{%please test!
   \PSOriginfalse% 
   \gdef\EPSFSpecial##1##2{\relax
    \special{language "PS",
      literal "##2 1000 div ##2 1000 div scale",
      position = "bottom left",
      include "##1"}}}
  \let\SetDVIALWEPSFSpecial\SetBeebeEPSFSpecial

 %% Northlake software
  \def\SetNorthlakeEPSFSpecial{\PSOrigintrue
   \gdef\EPSFSpecial##1##2{%
     \edef\specialthis{##2}%
     \SPLIT@0.@\specialthis.@\relax % suppress decimals (nec!)
     \special{insert ##1,magnification=\the\Initialtoks@}}}

 \def\SetStandardEPSFSpecial{%
   \gdef\EPSFSpecial##1##2{%
     \ms@g{}
     \ms@g{%
       !!! Sorry! There is still no standard for \string%
       \special\space EPSF integration !!!}%
     \ms@g{%
      --- So you will have to identify your driver using a command}%
     \ms@g{%
      --- of the form \string\Set...EPSFSpecial, in order to get}%
     \ms@g{%
      --- your graphics to print.  See BoxedEPS.doc.}%
     \ms@g{}
     \gdef\EPSFSpecial####1####2{}
     }}

  \SetStandardEPSFSpecial %% currently gives warning
 
 \let\wlog\wlog@ld %%restore logging 

 \catcode`\:=\C@tColon
 \catcode`\;=\C@tSemicolon
 \catcode`\?=\C@tQmark
 \catcode`\!=\C@tEmark
 \catcode`\"=\C@tDqt

 \catcode`\@=\EPSFCatAt

%%%%%%%%%%%% MACROS FOR PAPERS AND REPORTS %%%%%%%%%%%%
%My little macros
\def\G{\Gamma}
\def\sg{\partial}
\def\ul{\underline}
\def\ignore#1{}
\def\remark#1{\smskip\pn{\bf Remark #1}\quad} 
%%%%%%%%%%%%%%%%%%%%%

\newcount\sectnum
\newcount\subsectnum
\newcount\eqnumber

\global\eqnumber=1\sectnum=0

% Equation labels

\def\lab{(\the\sectnum.\the\eqnumber)}

%Example of use: suppose we want to give a label \lgh to an equation
% $$ ......  \xdef\lgh{\lab} \eqnum \show{lgh}$$
% Later refer to Eq. \lgh\ ...
% Note the \ after \lgh; it seems to be needed if we want the equation number
% to be followed by a space; not needed if followed by . or ,

%The next macro is used to display labels in drafts, so that you do
%not have to remember them

\def\show#1{#1}

%The next macro is to be used for final drafts that do not display labels
%\def\show#1{}

%%%%%%%%%%%%%%%%%%%%%

\def\magnify{\magnification=1200}

\def\smskip{\vskip 5 pt}
\def\medskip{\vskip 10 pt}
\def\bigskip{\vskip 15 pt}
\def\pn{\par\noindent}
\def\br{\break}
\def\un{\underline}
\def\ov{\overline}
\def\Bl{\Bigl}
\def\Br{\Bigr}
\def\bl{\bigl} 
\def\br{\bigr} 
\def\lf{\left}
\def\ri{\right}

\def\argmin{\mathop{\arg \min}}
\def\tendsd{\downarrow}
\def\mod{\rm mod}
\def\tends{\rightarrow}
\def\thereis{\exists}
\def\implies{\Rightarrow}
\def\implied{\Leftarrow}
\def\ubar{\underline}
\def\ol#1{\overline{#1}}
\def\kth{$k^{\rm th}$ }
\def\ith{$i^{\rm th}$ }
\def\jth{$j^{\rm th}$ }
\def\lth{$\ell^{\rm th}$ }
\def\grad{\nabla}
\def\tendsd{\downarrow}
\def\tr{ ^{\prime}}
\def\half{{\scriptstyle {1\over 2}}}
\def\Ascr{{\cal A}}
\def\Fscr{{\cal F}}
\def\Pscr{{\cal P}}
\def\Zscr{{\cal Z}}
\def\Nscr{{\cal N}}
\def\Cscr{{\cal C}}

\def\a{\alpha}
\def\be{\beta}
\def\b{\beta}
\def\l{\lambda}
\def\g{\gamma}
\def\m{\mu}
\def\p{\pi}
\def\r{\rho}
\def\e{\epsilon}
\def\t{\tau}
\def\d{\delta}
\def\s{\sigma}
\def\f{\phi}
\def\o{\omega}
\def\D{\Delta}
\def\P{\Pi}
\def\G{\Gamma}

\def\re{\Re}
\def\rn{\Re^n}

\def\thereis{\exists}
\def\gr{\nabla}
\def\lt{\lim_{t\tends\infty}}
\def\lty{\lim_{t\tends\infty}}
\def\noteq{\ne}
\def\mapr{:\re\mapsto\re}
\def\mapn{:\rn\mapsto\re}

\def\bn{\hfil\break}
\def\bnt{\hfil\break\indent} %break line; horizontal space
\def\tl{\tilde}
\def\ab{\allowbreak}
\def\old#1{}% invalidates text in braces 
\def\leaderfill{\leaders\hbox to 1em{\hss.\hss}\hfill}
%Example of use: \line{1. Optimality Conditions\leaderfill p.\ 2}

% John's macros
\def\ie{i.e.}
\def\eg{e.g.}

\parindent=2pc
\baselineskip=15pt
\vsize=8.7 true in
\voffset=0.125 true in
\parskip=3pt
\def\singlespace{\baselineskip 12 pt}
\def\onehalfspace{\baselineskip 18 pt}
\def\doublespace{\baselineskip 24 pt}
\def\normalspace{\baselineskip 15 pt}

% vector/matrix macros

\def\colvect#1#2{\lf(\matrix{#1\cr#2\cr}\ri)}
\def\rowvect#1#2{\lf(\matrix{#1&#2\cr}\ri)}
\def\twomat#1#2#3#4{\lf(\matrix{#1&#2\cr #3&#4\cr}\ri)}
\def\ncolvect#1#2{\lf(\matrix{#1\cr\vdots\cr#2\cr}\ri)}
\def\nrowvect#1#2{\lf(\matrix{#1&\ldots&#2\cr}\ri)}

%eqalign macros
\def\minprob#1#2#3{$$\eqalign{&\hbox{minimize\ \ }#1\cr &\hbox{subject to\ \
}#2\cr}\ifnum 0=#3{}\else\eqno(#3)\fi$$}        
\def\minprobn#1#2#3{$$\eqalign{&\hbox{minimize\ \ }#1\cr &\hbox{subject to\ \
}#2\cr}\eqno(#3)$$}     
\def\maxprob#1#2#3{$$\eqalign{&\hbox{maximize\ \ }#1\cr &\hbox{subject to\ \
}#2\cr}\ifnum 0=#3{}\else\eqno(#3)\fi$$}        
\def\maxprobn#1#2#3{$$\eqalign{&\hbox{maximize\ \ }#1\cr &\hbox{subject to\ \
}#2\cr}\eqno(#3)$$}     
\def\aligntwo#1#2#3#4#5{$$\eqalign{#1&#2\cr #3&#4\cr}
\ifnum 0=#5{}\else\eqno(#5)\fi$$}
\def\alignthree#1#2#3#4#5#6#7{$$\eqalign{#1&#2\cr #3&#4\cr #5&#6\cr}
\ifnum 0=#7{}\else\eqno(#7)\fi$$}

% Macros to automatically advance equation and other numbers

\def\eqnum{\eqno{\hbox{(\the\sectnum.\the\eqnumber)}\global\advance\eqnumber
by1}}

\def\eqnu{\eqno{\hbox{(\the\sectnum.\the\eqnumber)}\global\advance\eqnumber
by1}}

\newcount\examplnumber
\def\examplnum{\global\advance\examplnumber by1}

\newcount\figrnumber
\def\figrnum{\global\advance\figrnumber by1}

\newcount\propnumber
\def\propnum{\global\advance\propnumber by1}

\newcount\defnumber
\def\defnum{\global\advance\defnumber by1}

\newcount\lemmanumber
\def\lemmanum{\global\advance\lemmanumber by1}

\newcount\assumptionnumber
\def\assumptionnum{\global\advance\assumptionnumber by1}

\def\exampl{\the\sectnum.\the\examplnumber}
\def\figr{\the\sectnum.\the\figrnumber}
\def\propn{\the\sectnum.\the\propnumber}
\def\defn{\the\sectnum.\the\defnumber}
\def\lemman{\the\sectnum.\the\lemmanumber}
\def\assumptionn{\the\sectnum.\the\assumptionnumber}

\def\section#1{\goodbreak\vskip 3pc plus 6pt minus 3pt\leftskip=-2pc
   \global\advance\sectnum by 1\eqnumber=1
\global\examplnumber=1\figrnumber=1\propnumber=1\defnumber=1\lemmanumber=1\assumptionnumber=1%
   \line{\hfuzz=1pc{\hbox to 3pc{\bf %\the\sectnum.\quad
   \vtop{\hfuzz=1pc\hsize=38pc\hyphenpenalty=10000\noindent\uppercase{\the\sectnum.\quad #1}}\hss}}
			\hfill}
			\leftskip=0pc\nobreak\tenf
			\vskip 1pc plus 4pt minus 2pt\noindent\ignorespaces}

% ETP Macros

%\def\section#1{\goodbreak\vskip 3pc plus 6pt minus 3pt\leftskip=-2pc
%   \global\advance\sectnum by 1\eqnumber=1
%   \line{\hfuzz=1pc{\hbox to 3pc{\bf %\the\sectnum.\quad
%   \vtop{\hfuzz=1pc\hsize=38pc\hyphenpenalty=10000\noindent\uppercase{#1}}\hss}}
%                        \hfill}
%                        \leftskip=0pc\nobreak\tenf
%                        \vskip 1pc plus 4pt minus 2pt\noindent\ignorespaces}

\def\sect#1{\noindent\leftskip=-2pc\tenf
   \goodbreak\vskip 1pc plus 4pt minus 2pt
                \global\advance\subsectnum by 1\eqnumber=1
   \line{\hfuzz=1pc{\hbox to 3pc{\bf %\the\sectnum.\quad
   \vtop{\hfuzz=1pc\hsize=38pc\hyphenpenalty=10000\noindent\uppercase{{\bf #1}}}\hss}}
                        \hfill}
   \leftskip=0pc\nobreak\tenf
                        \vskip 1pc plus 4pt minus 2pt\nobreak\noindent\ignorespaces}

\def\subsection#1{\noindent\leftskip=0pc\tenf
   \goodbreak\vskip 1pc plus 4pt minus 2pt
%               \global\advance\subsectnum by 1
   \line{\hfuzz=1pc{\hbox to 3pc{\bf %\the\sectnum.\quad
   \vtop{\hfuzz=1pc\hsize=38pc\hyphenpenalty=10000\noindent{\bf #1}}\hss}}
                        \hfill}
   \leftskip=0pc\nobreak\tenf
                        \vskip 1pc plus 4pt minus 2pt\nobreak\noindent\ignorespaces}
\def\subsubsection#1{\goodbreak\vskip 1pc plus 4pt minus 2pt
   \hfuzz=3pc\leftskip=0pc\noindent\tenit #1 \nobreak\tenf\vskip 6pt plus 1pt
                                minus 1pt\nobreak\ignorespaces\leftskip=0pc}
%
%\def\rthl{Sec. \the\chapnum.\the\sectnum}                      
%\def\rthc{#1}\nobreak\noindent\ignorespaces
%\newcount\sectnum \sectnum=0
%\newcount\subsectnum \subsectnum=0
\def\textlist#1{\par{\bf #1}\ }

\def\beginexample#1{\noindent\goodbreak\vskip 6pt plus 1pt minus 1pt
\noindent
  \hbox {\bf Example #1\hss}%\break%\noindent
  \nobreak\vskip 4pt plus 1pt minus 1pt \nobreak\noindent\ninef
  \global\advance
                \leftskip by\parindent\pn}
\def\endexample{\vskip 12pt\tenf\par
  \global\advance\leftskip by -\parindent
  }

\def\beginexercise#1{\noindent\goodbreak\vskip 6pt plus 1pt minus 1pt \noindent\global\normalbaselineskip=12pt
  \hbox {\bf Exercise #1\hss}%\break%\noindent
  \nobreak\vskip 4pt plus 1pt minus 1pt 
  \nobreak\noindent\ninef\global\advance\leftskip
                        by\parindent\pn}
\def\endexercise{\vskip 12pt\tenf\par
  \global\advance\leftskip by -\parindent
  }

\def\beginsection#1{\noindent\goodbreak\vskip 6pt plus 1pt minus 1pt \noindent\global\normalbaselineskip=12pt
  \hbox {\it #1\hss}
  \vskip 0.1pt plus 1pt minus 1pt \nobreak\noindent\ninef\global\advance
                \leftskip by\parindent\noindent\pn}
\def\endsection{\vskip 12pt\tenf\par
  \global\advance\leftskip by -\parindent
}

\def\enddisplaylist{\vskip 12pt\par}

%%%%%%%%% REVISED \section and \subsection MACROS FOR NUMBERED SUBSECTIONS %%%

\def\section#1{\goodbreak\vskip 3pc plus 6pt minus 3pt\leftskip=-2pc
   \global\advance\sectnum by 1\eqnumber=1
\global\examplnumber=1\figrnumber=1\propnumber=1\defnumber=1\lemmanumber=1\assumptionnumber=1\subsectnum=0%
   \line{\hfuzz=1pc{\hbox to 3pc{\bf %\the\sectnum.\quad
   \vtop{\hfuzz=1pc\hsize=38pc\hyphenpenalty=10000\noindent\uppercase{\the\sectnum.\quad #1}}\hss}}
			\hfill}
			\leftskip=0pc\nobreak\tenf
			\vskip 1pc plus 4pt minus 2pt\noindent\ignorespaces}

\def\subsection#1{\noindent\leftskip=0pc\tenf
   \goodbreak\vskip 1pc plus 4pt minus 2pt
               \global\advance\subsectnum by 1
   \line{\hfuzz=1pc{\hbox to 3pc{\bf \the\sectnum.\the\subsectnum\ \ \
   \vtop{\hfuzz=1pc\hsize=38pc\hyphenpenalty=10000\noindent{\bf #1}}\hss}}
                        \hfill}
   \leftskip=0pc\nobreak\tenf
                        \vskip 1pc plus 4pt minus 2pt\nobreak\noindent\ignorespaces}

\def\subsubsection#1{\goodbreak\vskip 1pc plus 4pt minus 2pt
   \hfuzz=3pc\leftskip=0pc\noindent{\bf #1} \nobreak\vskip 6pt plus 1pt
                                minus 1pt\nobreak\ignorespaces\leftskip=0pc}

%%%%%%%%% END OF REVISED \section and \subsection MACROS %%%%%%%%%

% Header/Title macros

\def\lemma#1{\smskip\pn{\bf Lemma #1}\quad}
\def\theorem#1{\smskip\pn{\bf Theorem #1}\quad}
\def\proposition#1{\smskip\pn{\bf Proposition #1}\quad}
\def\proof{\smskip\pn{\bf Proof:}\quad} 
\def\definition#1{\smskip\pn{\bf
Definition #1}\quad} \def\assumption#1{\smskip\pn{\bf Assumption #1}\quad}
\def\corollary#1{\smskip\pn{\bf Corollary #1}\quad}
\def\exercise#1{\smskip\pn{\bf Exercise #1}\quad}
\def\figure#1{\smskip\pn{\bf Figure #1}\quad}
\def\QED{\quad{\bf Q.E.D.} \par\bigskip} \def\qed{\quad{\bf
Q.E.D.} \par\bigskip}
\def\ref{\smskip\pn}

\def\chapter#1#2{{\bf \centerline{\helbigbig
{#1}}}\bigskip\bigskip{\bf \centerline{\helbigbig
{#2}}}\bigskip\bigskip} % ex. \chapter{Chapter 1}{Title of chapter}

\def\longchapter#1#2#3{{\bf \centerline{\helbigbig
{#1}}}\bigskip\bigskip{\bf \centerline{\helbigbig
{#2}}}\bigskip{\bf \centerline{\helbigbig
{#3}}}\bigskip\bigskip} % ex. \longchapter{Chapter 1}{Title of chapter}{Title of
 %chapter}

\def\papertitle#1#2{{\bf \centerline{\helbigb
{#1}}}\bigskip\bigskip{\centerline{
by}}\bigskip{\bf \centerline{
{#2}}}\bigskip\bigskip} % ex. \papertitle{Title of paper}{Names of Authors}

\def\longpapertitle#1#2#3{{\bf \centerline{\helbigb
{#1}}}\bigskip{\bf \centerline{\helbigb
{#2}}}\bigskip\bigskip{\centerline{
by}}\bigskip{\bf \centerline{
{#3}}}\bigskip\bigskip} 
% ex. \longpapertitle{First part of title of paper}
%{2nd part of title of paper}{Names of Authors}

% List macros

\def\nitem#1{\smskip\item{#1}}
\def\nitemitem#1{\smskip\itemitem{#1}}
\def\endlist{\smskip}

\newcount\alphanum
\newcount\romnum

\def\alphaenumerate{\ifcase\alphanum \or (a)\or (b)\or (c)\or (d)\or (e)\or
(f)\or (g)\or (h)\or (i)\or (j)\or (k)\fi}
\def\romenumerate{\ifcase\romnum \or (i)\or (ii)\or (iii)\or (iv)\or (v)\or
(vi)\or (vii)\or (viii)\or (ix)\or (x)\or (xi)\fi}

\def\alist{\begingroup\vskip10pt\alphanum=1% alphabetical list
\parskip=2pt\parindent=0pt \leftskip=3pc
\everypar{\llap{{\rm\alphaenumerate\hskip1em}}\advance\alphanum by1}}
\def\endalist{\vskip1pc\endgroup\parskip=0pt
  \parindent=\bigskipamount \leftskip=0pc\tenf
                        \noindent\ignorespaces}

\def\nolist{\begingroup\vskip10pt\alphanum=0% numerical list
\parskip=2pt\parindent=0pt \leftskip=3pc
\everypar{\llap{\global\advance\alphanum by1(\the\alphanum)\hskip1em}}}
\def\endnolist{\vskip1pc\endgroup\parskip=0pt\leftskip=0pc\tenf
                        \noindent\ignorespaces}

\def\romlist{\begingroup\vskip10pt\romnum=1% roman list
\parskip=2pt\parindent=0pt \leftskip=5pc
\everypar{\llap{{\rm\romenumerate\hskip1em}}\advance\romnum by1}}
\def\endromlist{\vskip1pc\endgroup\parskip=0pt\leftskip=0pc\tenf
                        \noindent\ignorespaces}
% romlist indents more than alist or nolist and can be used inside them

%Figure, table, and box macros

\long\def\fig#1#2#3{\vbox{\vskip1pc\vskip#1
\prevdepth=12pt \baselineskip=12pt
\vskip1pc
\hbox to\hsize{\hfill\vtop{\hsize=25pc\noindent{\eightbf Figure #2\ }
{\eightpoint#3}}\hfill}}}%Figure space definition. Example of use:
%\fig{16pc}{1.1}{A network with one central processor and a separate
%communication link to each device.}

\long\def\widefig#1#2#3{\vbox{\vskip1pc\vskip#1
\prevdepth=12pt \baselineskip=12pt
\vskip1pc
\hbox to\hsize{\hfill\vtop{\hsize=28pc\noindent{\eightbf Figure #2\ }
{\eightpoint#3}}\hfill}}}

\long\def\table#1#2{\vbox{\vskip0.5pc
\prevdepth=12pt \baselineskip=12pt
\hbox to\hsize{\hfill\vtop{\hsize=25pc\noindent{\eightbf Table #1\ }
{\eightpoint#2}}\hfill}}}

%Running Head Macros
\def\rightleftheadline#1#2{ifodd\pageno\rightheadline{#1}
\else\leftheadline{#2}\fi} 
\def\rightheadline#1{\headline{\tenrm\hfil #1}}
\def\leftheadline#1{\headline{\tenrm#1\hfil}}

% Concept Macros

\long\def\leftfig#1#2{\vbox{\smskip\hsize=220pt
\vtop{{\noindent {\bf #1}}}
\smskip
\noindent
\vbox{{\noindent #2}}
}}

\long\def\rightfig#1#2#3{\vbox{\smskip\vskip#1
\prevdepth=12pt \baselineskip=12pt
\hsize=210pt
\smskip
\vbox{\noindent{\eightbold #2}
\hskip1em{\eightpoint#3}}
}}

\long\def\concept#1#2#3#4#5{\bigskip\hrule
\vbox{\hbox{\leftfig{#1}{#2} \hskip3em
\rightfig{#3}{#4}{#5}} \smskip}
\hrule\bigskip}

% Example of Use: \concept{Title of Concept}{Text}
% {Figure size}{Figure number?}{Figure caption}

\long\def\bconcept#1#2#3#4#5#6#7{
\vbox{
\hbox to \hsize{\vtop{\par #1}}
\concept{#2}{#3}{#4}{#5}{#6}
\hbox to \hsize{\vtop{\par #7}}
\smskip}
}

% Example of Use: \bconcept{Preceding text}{Title of Concept}{Text}
% {Figure size}{Figure number}{Figure caption}{Following text}

\long\def\boxconcept#1#2#3#4#5{
\vbox{\hbox{\leftfig{#1}{#2} \hskip3em
\rightfig{#3}{#4}{#5}} \smskip}
}% same as concept but without the \hrule's; ready to be boxed

% Put inside a box

\def\boxit#1{\vbox{\hrule\hbox{\vrule\kern3pt
                                \vbox{\kern3pt#1\kern3pt}\kern3pt\vrule}\hrule}}
% example of use: \setbox0=\vbox{.... }; \boxit{\box0}
\def\centerboxit#1{$$\vbox{\hrule\hbox{\vrule\kern3pt
                                \vbox{\kern3pt#1\kern3pt}\kern3pt\vrule}\hrule}$$}
% example of use: \setbox0=\vbox{.... }; \centerboxit{\box0}

\long\def\boxtext#1#2{$$\boxit{\vbox{\hsize #1\noindent\strut #2\strut}}$$}
% example of use: \boxtext{462pt}{This is the boxed text.}; 462pt is max length

% Picture macros and examples
%
% figures must be pasted from mcdraw
%
% Look in the 'Windows' menu for the pictures window
% It's like the Scrapbook -- cut and paste pictures
%

\def\picture #1 by #2 (#3){
  \vbox to #2{
    \hrule width #1 height 0pt depth 0pt
    \vfill
    \special{picture #3} % this is the low-level interface
    }
  }
% The first dimension of the picture macro is the width the second is depth

\def\scaledpicture #1 by #2 (#3 scaled #4){{
  \dimen0=#1 \dimen1=#2
  \divide\dimen0 by 1000 \multiply\dimen0 by #4
  \divide\dimen1 by 1000 \multiply\dimen1 by #4
  \picture \dimen0 by \dimen1 (#3 scaled #4)}
  }

%
% Note that you can also say, e.g.,
%  \special{postscript xxx yyy zzz}
% to include PostScript graphics in your documents
%
%Examples of use
%\def\stripes{\picture 2.29in by 1.75in (AWstripes)}
% By executing \stripes 
%\def\annie{\scaledpicture 102pt by 239pt (annie scaled 2000)}
%\def\finder{\picture 260pt by 165pt (screen0 scaled 500)}
%\def\icon{\picture 7in by 7in (icon)}
%Example of use
%\annie
%Example of centered picture \line{\hfil\annie\hfil}

%Figure w/  caption macro
\long\def\captfig#1#2#3#4#5{\vbox{\vskip1pc
\hbox to\hsize{\hfill{\picture #1 by #2 (#3)}\hfill}
\prevdepth=9pt \baselineskip=9pt
\vskip1pc
\hbox to\hsize{\hfill\vtop{\hsize=24pc\noindent{\eightbold Figure #4}
\hskip1em{\eightpoint#5}}\hfill}}}

%Examples of use of Figure macros
%\captfig{8.53pc}{19.9pc}{picturename}{5}{caption.}
%\captfig{2.23in}{2in}{picturename scaled 500}{16}{Caption.}
%The macro centers the picture.
%The first two numbers should be the true width
% and height after the picture has been scaled.
% So if the picture is scaled by 50% (500), the width and height in
% the macro should onw half of what they would be if the picture
% is not scaled (1000).
%
%
%
% Postcript macros

\def\illustration #1 by #2 (#3){
  \vskip#2\hskip#1\special{illustration #3} % this is the low-level interface
    }

\def\scaledillustration #1 by #2 (#3 scaled #4){{
  \dimen0=#1 \dimen1=#2
  \divide\dimen0 by 1000 \multiply\dimen0 by #4
  \divide\dimen1 by 1000 \multiply\dimen1 by #4
  \illustration \dimen0 by \dimen1 (#3 scaled #4)}
  }

%***************************************************
%         FONTS
%***************************************************

% ROMAN
\font\hel=cmr10%
\font\helb=cmbx10%
\font\heli=cmti10%
\font\helbi=cmsl10%
\font\ninehel=cmr9%
\font\nineheli=cmti9%
\font\ninehelb=cmbx9%
\font\helbig=cmr10 scaled 1500%
\font\helbigbig=cmr10 scaled 2500%
\font\helbigb=cmbx10 scaled 1500%
\font\helbigbigb=cmbx10 scaled 2500%
\font\bigi=cmti10 scaled \magstep5%
\font\eightbold=cmbx8%
\font\ninebold=cmbx9%
\font\boldten=cmbx10%

\def\tenf{\hel}%
\def\tenit{\heli}%
\def\ninef{\ninehel}%
\def\nineb{\ninehelb}%
\def\nineit{\nineheli}%
\def\smit{\nineheli}%
\def\smbf{\ninehelb}%

%  FONT FAMILIES

\font\tenrm=cmr10%
\font\teni=cmmi10%
\font\tensy=cmsy10%
\font\tenbf=cmbx10%
\font\tentt=cmtt10%
\font\tenit=cmti10%
\font\tensl=cmsl10%

\def\tenpoint{\def\rm{\fam0\tenrm}%
\textfont0=\tenrm%
\textfont1=\teni%
\textfont2=\tensy%
\textfont\itfam=\tenit%
\textfont\slfam=\tensl%
\textfont\ttfam=\tentt%
\textfont\bffam=\tenbf%
\scriptfont0=\sevenrm%
\scriptfont1=\seveni%
\scriptfont2=\sevensy%
%\scriptfont3=\tenex%
\scriptscriptfont0=\sixrm%
\scriptscriptfont1=\sixi%
\scriptscriptfont2=\sixsy%
%\scriptscriptfont3=\tenex%
\def\it{\fam\itfam\tenit}%
\def\tt{\fam\ttfam\tentt}%
\def\sl{\fam\slfam\tensl}%
\scriptfont\bffam=\sevenbf%
\scriptscriptfont\bffam=\sixbf%
\def\bf{\fam\bffam\tenbf}%
\normalbaselineskip=18pt%
\normalbaselines\rm}%

\font\ninerm=cmr9%
\font\ninebf=cmbx9%
\font\nineit=cmti9%
\font\ninesy=cmsy9%
\font\ninei=cmmi9%
\font\ninett=cmtt9%
\font\ninesl=cmsl9%

\def\ninepoint{\def\rm{\fam0\ninerm}%
\textfont0=\ninerm%
\textfont1=\ninei%
\textfont2=\ninesy%
\textfont\itfam=\nineit%
\textfont\slfam=\ninesl%
\textfont\ttfam=\ninett%
\textfont\bffam=\ninebf%
\scriptfont0=\sixrm%
\scriptfont1=\sixi%
\scriptfont2=\sixsy%
%\scriptfont3=\tenex%
\def\it{\fam\itfam\nineit}%
\def\tt{\fam\ttfam\ninett}%
\def\sl{\fam\slfam\ninesl}%
\scriptfont\bffam=\sixbf%
\scriptscriptfont\bffam=\fivebf%
\def\bf{\fam\bffam\ninebf}%
\normalbaselineskip=16pt%
\normalbaselines\rm}%

\font\eightrm=cmr8%
\font\eighti=cmmi8%
\font\eightsy=cmsy8%
\font\eightbf=cmbx8%
\font\eighttt=cmtt8%
\font\eightit=cmti8%
\font\eightsl=cmsl8%

\def\eightpoint{\def\rm{\fam0\eightrm}%
\textfont0=\eightrm%
\textfont1=\eighti%
\textfont2=\eightsy%
\textfont\itfam=\eightit%
\textfont\slfam=\eightsl%
\textfont\ttfam=\eighttt%
\textfont\bffam=\eightbf%
\scriptfont0=\sixrm%
\scriptfont1=\sixi%
\scriptfont2=\sixsy%
%\scriptfont3=\tenex%
\scriptscriptfont0=\fiverm%
\scriptscriptfont1=\fivei%
\scriptscriptfont2=\fivesy%
%\scriptscriptfont3=\tenex%
\def\it{\fam\itfam\eightit}%
\def\tt{\fam\ttfam\eighttt}%
\def\sl{\fam\slfam\eightsl}%
%\scriptfont\bffam=\sixbf%
\scriptscriptfont\bffam=\fivebf%
\def\bf{\fam\bffam\eightbf}%
\normalbaselineskip=14pt%
\normalbaselines\rm}%

\font\sevenrm=cmr7%
\font\seveni=cmmi7%
\font\sevensy=cmsy7%
\font\sevenbf=cmbx7%
\font\seventt=cmtt8 at 7pt%
\font\sevenit=cmti8 at 7pt%
\font\sevensl=cmsl8 at 7pt%

\def\sevenpoint{%
   \def\rm{\sevenrm}\def\bf{\sevenbf}%
   \def\smc{\sevensmc}\baselineskip=12pt\rm}%

\font\sixrm=cmr6%
\font\sixi=cmmi6%
\font\sixsy=cmsy6%
\font\sixbf=cmbx6%
\font\sixtt=cmtt8 at 6pt%
\font\sixit=cmti8 at 6pt%
\font\sixsl=cmsl8 at 6pt%
\font\sixsmc=cmr8 at 6pt%

\def\sixpoint{%
   \def\rm{\sixrm}\def\bf{\sixbf}%
   \def\smc{\sixsmc}\baselineskip=12pt\rm}%

\fontdimen13\tensy=2.6pt%
\fontdimen14\tensy=2.6pt%
\fontdimen15\tensy=2.6pt%
\fontdimen16\tensy=1.2pt%
\fontdimen17\tensy=1.2pt%
\fontdimen18\tensy=1.2pt%       

\def\tenf{\tenpoint}%
\def\ninef{\ninepoint}%
\def\eightf{\eightpoint}%

%%%%%%%%%%%% END OF MACROS %%%%%%%%%%%%

%\input TEXSHOP_small_baseline.tex

%***************************************************
%         FONTS
%***************************************************

\def\tenpoint{\def\rm{\fam0\tenrm}%
\textfont0=\tenrm%
\textfont1=\teni%
\textfont2=\tensy%
\textfont\itfam=\tenit%
\textfont\slfam=\tensl%
\textfont\ttfam=\tentt%
\textfont\bffam=\tenbf%
\scriptfont0=\sevenrm%
\scriptfont1=\seveni%
\scriptfont2=\sevensy%
%\scriptfont3=\tenex%
\scriptscriptfont0=\sixrm%
\scriptscriptfont1=\sixi%
\scriptscriptfont2=\sixsy%
%\scriptscriptfont3=\tenex%
\def\it{\fam\itfam\tenit}%
\def\tt{\fam\ttfam\tentt}%
\def\sl{\fam\slfam\tensl}%
\scriptfont\bffam=\sevenbf%
\scriptscriptfont\bffam=\sixbf%
\def\bf{\fam\bffam\tenbf}%
\normalbaselineskip=14pt%
\normalbaselines\rm}%

\def\ninepoint{\def\rm{\fam0\ninerm}%
\textfont0=\ninerm%
\textfont1=\ninei%
\textfont2=\ninesy%
\textfont\itfam=\nineit%
\textfont\slfam=\ninesl%
\textfont\ttfam=\ninett%
\textfont\bffam=\ninebf%
\scriptfont0=\sixrm%
\scriptfont1=\sixi%
\scriptfont2=\sixsy%
%\scriptfont3=\tenex%
\def\it{\fam\itfam\nineit}%
\def\tt{\fam\ttfam\ninett}%
\def\sl{\fam\slfam\ninesl}%
\scriptfont\bffam=\sixbf%
\scriptscriptfont\bffam=\fivebf%
\def\bf{\fam\bffam\ninebf}%
\normalbaselineskip=13pt%
\normalbaselines\rm}%

\def\eightpoint{\def\rm{\fam0\eightrm}%
\textfont0=\eightrm%
\textfont1=\eighti%
\textfont2=\eightsy%
\textfont\itfam=\eightit%
\textfont\slfam=\eightsl%
\textfont\ttfam=\eighttt%
\textfont\bffam=\eightbf%
\scriptfont0=\sixrm%
\scriptfont1=\sixi%
\scriptfont2=\sixsy%
%\scriptfont3=\tenex%
\scriptscriptfont0=\fiverm%
\scriptscriptfont1=\fivei%
\scriptscriptfont2=\fivesy%
%\scriptscriptfont3=\tenex%
\def\it{\fam\itfam\eightit}%
\def\tt{\fam\ttfam\eighttt}%
\def\sl{\fam\slfam\eightsl}%
%\scriptfont\bffam=\sixbf%
\scriptscriptfont\bffam=\fivebf%
\def\bf{\fam\bffam\eightbf}%
\normalbaselineskip=12pt%
\normalbaselines\rm}%

\def\sevenpoint{%
   \def\rm{\sevenrm}\def\bf{\sevenbf}%
   \def\smc{\sevensmc}\baselineskip=10pt\rm}%
\def\sixpoint{%
   \def\rm{\sixrm}\def\bf{\sixbf}%
   \def\smc{\sixsmc}\baselineskip=9pt\rm}%

\def\skiptonextpage{\vfill\eject}

%%%%%%%%%% REDEFINITION OF BOX SPACING %%%%%%%%%%%%%%%%

\def\texshopbox#1{\boxtext{462pt}{\vskip-1.5pc\pshade{\vskip-1.0pc#1\vskip-2.0pc}}}
\def\texshopboxnt#1{\boxtextnt{462pt}{\vskip-1.5pc\pshade{\vskip-1.0pc#1\vskip-2.0pc}}}
\def\texshopboxnb#1{\boxtextnb{462pt}{\vskip-1.5pc\pshade{\vskip-1.0pc#1\vskip-2.0pc}}}

%%%%%%%%%%%%%%%%%%%%%%%%%%%%%%%%%%%%%%%%%%%%

\input miniltx

\ifx\pdfoutput\undefined
  \def\Gin@driver{dvips.def} % we are not running PDFTeX
\else
  \def\Gin@driver{pdftex.def} % we are running PDFTeX
\fi

\input graphicx.sty
\resetatcatcode

\long\def\fig#1#2#3{\vbox{\vskip1pc\vskip#1
\prevdepth=12pt \baselineskip=12pt
\vskip1pc
\hbox to\hsize{\hskip3pc\hfill\vtop{\hsize=35pc\noindent{\eightbf Figure #2\ }
{\eightpoint#3}}\hfill}}}

\def\show#1{}
\def\dpbshow#1{#1}
\def\frac#1#2{{#1\over #2}}

\rightheadline{\botmark}

\pageno=1

\rightheadline{\botmark}

\pn {\bf October 2016 (Revised Jan.\ 2017 and Nov.\ 2017)}\hfill{\bf  Report LIDS-P-3205}
\bigskip 
\bigskip\bigskip

\bigskip\bigskip

\def\longpapertitle#1#2#3{{\bf \centerline{\helbigb
{#1}}}\medskip{\bf \centerline{\helbigb
{#2}}}\medskip{\centerline{
by}}\medskip{\bf \centerline{
{#3}}}\bigskip}

\def\longpapertitle#1#2#3{{\bf \centerline{\helbigb
{#1}}}\medskip{\bf \centerline{\helbigb
{#2}}}\medskip{\centerline{
}}\medskip{\bf \centerline{
{#3}}}\bigskip}

\vskip-4.3pc
\longpapertitle{Proximal Algorithms and Temporal Differences for Large }{Linear Systems: Extrapolation, Approximation, and Simulation\footnote{\dag}
{\ninepoint  This report is an extended version of a paper titled ``Proximal Algorithms and Temporal Difference Methods for Solving Fixed Point Problems," which will appear in Computational Optimization and Applications Journal.}}
{Dimitri P.\ Bertsekas
\footnote{\ddag}
{\ninepoint  The author is with the Dept.\ of Electr.\ Engineering and
Comp.\ Science, and the Laboratory for Information and Decision Systems, M.I.T., Cambridge, Mass., 02139.}
}

\centerline{\bf Abstract}
We consider large linear and nonlinear fixed point problems, and solution with proximal algorithms. We show that there is a close connection between two seemingly different types of methods from distinct fields: 1) Proximal iterations for linear systems of equations, which are prominent in numerical analysis and convex optimization, and 2) Temporal difference (TD) type methods, such as TD($\l$), LSTD($\l$), and LSPE($\l$), which are central in simulation-based approximate dynamic programming/reinforcement learning (DP/RL), and its recent prominent successes in large-scale game contexts, among others. 

One benefit of this connection is a new and simple way to accelerate the standard proximal algorithm by extrapolation towards the TD iteration, which generically has a faster convergence rate. Another benefit is the potential integration into the proximal algorithmic context of several new ideas that have emerged in the DP/RL context.  We discuss some of the possibilities, and in particular, algorithms that  project each proximal iterate onto the subspace spanned by a small number of  basis functions, using low-dimensional calculations and simulation. A third benefit is that insights and analysis from proximal algorithms can be brought to bear on the enhancement of TD methods.

The linear fixed point methodology can be extended to nonlinear fixed point problems involving a contraction, thus providing guaranteed and potentially substantial acceleration of the proximal and forward backward splitting algorithms at no extra cost. Moreover, the connection of proximal and TD methods can be extended to nonlinear (nondifferentiable) fixed point problems through new proximal-like algorithms that involve successive linearization, similar to policy iteration in DP.

\def\old#1{}% invalidates text in braces 
\xdef\z{\zeta}

\def\L{\Lambda}

\xdef\td{T^{(\l)}}
\xdef\ptd{\Pi T^{(\l)}}

\xdef\pb{P^{(c)}}
\xdef\ppb{\Pi P^{(c)}}
\xdef\ala{A^{(\l)}}
\xdef\bla{b^{(\l)}}

\vskip  -4mm
\section{Introduction}
\vskip  -2mm

\pn In this paper we focus primarily on systems of linear equations of the form
$$x=Ax+b,\xdef\linsys{\lab}\eqnum\show{oneo}$$
where $A$ is an $n\times n$ matrix and $b$ is a column vector in the $n$-dimensional space $\rn$. We denote by $\s(M)$ the  spectral radius of a square matrix $M$ (maximum over the moduli of the eigenvalues of $M$), and we assume the following.

\xdef\assumptiona{\assumptionn}\assumptionnum\show{myproposition}

%\texshopbox
{\assumption{\assumptiona}The matrix $I-A$ is invertible and $\s(A)\le 1$.} 

We consider the proximal algorithm, originally proposed for the solution of monotone variational inequalities by Martinet [Mar70] (see also the textbook treatments by Facchinei  and Pang [FaP03],  Bauschke and Combettes [BaC11], and the author's [Ber15]). This algorithm has the form 
$$x_{k+1}=\pb x_k,$$
where $c$ is a positive scalar, and for a given $x\in\rn$, $\pb x$ denotes the solution of the following equation in the vector $y$:
$$y=Ay+b+{1\over c}(x-y).$$   
Under Assumption \assumptiona, this equation has the unique solution 
$$\pb x=\Bigg({c+1\over c}I-A\Bigg)^{-1}\lf(b+{1\over c} x\ri),\xdef\pbformula{\lab}\eqnum\show{oneo}$$
because the matrix ${c+1\over c}I-A$ is invertible, since its eigenvalues lie within the unit circle that is centered at ${c+1\over c}$, so they do not include 0. 

When $A$ is symmetric, the system \linsys\ is the optimality condition for the minimization
$$\min_{x\in\rn}\left\{{1\over 2} x'Qx-b'x\right\},\xdef\quadmin{\lab}\eqnum\show{oneo}$$
where $Q=I-A$ and a prime denotes transposition. The proximal algorithm $x_{k+1}=\pb x_k$ can then be implemented through the minimization
$$x_{k+1}\in \argmin_{x\in\rn}\left\{{1\over 2} x'Qx-b'x+{1\over 2c}\|x-x_k\|^2\right\},$$
or 
$$x_{k+1}=\lf({1\over c} I+Q\ri)^{-1}\lf(b+{1\over c} x_k\ri).$$
In this case, Assumption \assumptiona\ is equivalent to $Q$ being positive definite, with all eigenvalues in the interval $(0,2]$. Note, however, that for the minimization problem \quadmin, the proximal algorithm is convergent for any positive semidefinite symmetric $Q$, as is well known. Thus Assumption \assumptiona\ is not the most general assumption under which the proximal algorithm can be applied.\footnote{\dag}
{\ninepoint  It is possible to scale the eigenvalues of $Q$ to lie in the interval $(0,2]$ without changing the problem, by multiplying $Q$ and $b$ with a suitable positive scalar. This, however, requires some prior knowledge about the location of the eigenvalues of $Q$.}
 Still, however,  the assumption covers important types of problems, including the case where $A$ is a contraction with respect to some norm, as well as applications in dynamic programming (DP for short), to be discussed shortly.

Let us denote by $T$ the mapping whose fixed point we wish to find,
$$Tx=Ax+b.$$
We  denote by $T^\ell$ the $\ell$-fold composition of $T$, where $\ell$ is a positive integer, and we define addition of a finite number and an infinite number of linear operators in the standard way. We introduce the multistep mapping  $\td$ given by
$$\td =(1-\l)\sum_{\ell=0}^\infty \l^{\ell}T^{\ell+1},\xdef\tddef{\lab}\eqnum\show{oneo}$$
where $\l$ is a scalar with $0<\l<1$. The series defining $\td$ is convergent under Assumption \assumptiona, as we will discuss later. A principal aim of this paper is to establish the relation between the mappings $\td$ and $\pb$, and the ways in which this relation can be exploited algorithmically to compute its fixed point.

The mapping $\td$ has been central in the field that we will refer to as ``approximate DP" (the name ``reinforcement learning" is also often used in artificial intelligence, and the names ``neuro-dynamic programming" and ``adaptive dynamic programming" are often used in automatic control, with essentially the same meaning). In particular, $\td$ is involved in methods for finding a fixed point of the mapping $\ptd$, where $\P$ is either the identity or some form of projection onto a low-dimensional subspace $S$.\footnote{\dag}{\ninepoint In approximate DP it is common to replace a fixed point equation of the form $x=F(x)$ with the equation $x=\Pi \bl(F(x)\br)$. This approach comes under the general framework of Galerkin approximation, which is widely used in a variety of numerical computation contexts (see e.g., the books by Krasnoselskii [Kra72] and Fletcher [Fle84], and the DP-oriented discussion in the paper by Yu and Bertsekas [YuB10]). A distinguishing feature of approximate DP applications is that $F$ is a linear mapping and the equation $x=\Pi \bl(F(x)\br)$ is typically solved by simulation-based methods.} In the DP context, $A$ is a substochastic matrix related to the Markov chain of a policy and the equation $x=Ax+b$ is the Bellman equation for the cost function $x$ of the policy. Equations of the form $x=T x$ are solved repeatedly within the exact policy iteration method, which generates a sequence of improved cost functions and associated policies. Equations of the form  $x=\ptd x$ are solved within a corresponding approximate policy iteration method. Detailed accounts of the approximate DP context are given in several books, including the ones by Bertsekas and Tsitsiklis [BeT96], Sutton and Barto [SuB98], Si et al.\ [SBP04], Powell [Pow07], Busoniu et al.\ [BBD10], Szepesvari [Sze10], Bertsekas [Ber12a],  Lewis and Liu [LeL13], and Vrabie, Vamvoudakis, and Lewis [VVL13]. Substantial computational experience has been accumulated with this methodology, and considerable success has been obtained (including prominent achievements with programs that play games, such as Backgammon, Go, and others, at impressive and sometimes above human level; see Tesauro [Tes94], Tesauro and Galperin [TeG96], Scherrer et al.\ [GMS13], [SMG15], and Silver et al.\ [MKS15], [SHM16]). In challenging approximate DP applications, the dimension of $A$ is very high, the dimension of the approximation subspace $S$ is low by comparison, and the large-scale computations involved in calculating the fixed point of $\ptd$ are handled by Monte-Carlo simulation schemes.

A variety of simulation-based methods that involve the mapping $\td$, such as TD($\l$), LSTD($\l$), and LSPE($\l$), have been proposed in approximate DP. In particular,  the fixed point iteration $x_{k+1}=\ptd x_k$ (where $\P$ is orthogonal projection with respect to a weighted Euclidean norm) has been called PVI($\l$)  in the author's DP textbook [Ber12a] (PVI stands for Projected Value Iteration).
 Its simulation-based implementation is the LSPE($\l$) method  (LSPE stands for Least Squares Policy Evaluation) given in joint works of the author with his collaborators Ioffe, Nedi\'c, Borkar, and Yu [BeI96], [NeB03], [BBN04], [YuB06], [Ber11b]. 
The simulation-based matrix inversion method that solves the fixed point equation $x=\ptd x$ is the LSTD($\l$) method, given by 
Bradtke and Barto  [BrB96], 
and further discussed, extended, and analyzed by Boyan [Boy02], Lagoudakis  and Parr [LaP03], Nedi\'c and Bertsekas [NeB03], Bertsekas and Yu [BeY09], [YuB12], and Yu [Yu10], [Yu12] (LSTD stands for Least Squares Temporal Differences). TD($\l$), proposed by  Sutton [Sut88] in the approximate DP setting, is a stochastic approximation method for solving the equation $x=\ptd x$. It has the form
$$x_{k+1}=x_k+\g_k\Big({sample}\big({\ptd x_k}\big)-x_k\Big),\xdef\tdlambda{\lab}\eqnum\show{oneo}$$
where ${sample}\big({\ptd x_k}\big)$ is a stochastic simulation-generated sample of $\ptd x_k$, and $\g_k$ is a diminishing stepsize satisfying standard conditions for stochastic iterative algorithms,  such as $\g_k=1/(k+1)$.\footnote{\dag}{\ninepoint   The precise nature of the problem that TD($\l$) is aiming to solve was unclear for a long time. The paper by Tsitsiklis and VanRoy [TsV97] showed that it aims to find a fixed point of $\td$ or $\P\td$, and gave the first convergence analysis (also replicated in the book [BeT96]). The paper by Bertsekas and Yu [BeY09] (Section 5.3) generalized TD($\l$), LSTD($\l$), and LSPE($\l$) to the linear system context of this paper.}  The computation of the samples in Eq.\ \tdlambda\ involves simulation using Markov chains and the notion of temporal differences, which originated in reinforcement
learning with the works of Samuel [Sam59], [Sam67] on a checkers-playing program. Mathematically, temporal differences are residual-type terms of the form $A^\ell (Ax+b-x)$, $\ell\ge0$, which can be used to streamline various computations within the aforementioned methods. We refer to the sources given above for methods to generate samples of temporal differences in the DP/policy evaluation context, and to [BeY09] for corresponding methods and analysis within the more general linear fixed point context of the present paper.

 \xdef \figoneone{\figr}\figrnum\show{myfigure}

A central observation of this paper, shown  in Section 2, is that the proximal mapping $\pb$ is closely related to the multistep mapping $\td$, 
where
$$\l={c\over c+1}.$$
In particular $\pb$ is an interpolated mapping between $\td$ and the identity, or reversely, $\td$ is an extrapolated form of $\pb$; see  Fig.\ \figoneone. Moreover, we show in Section 2 that under Assumption \assumptiona, $\td$ has a smaller spectral radius than $\pb$, and as a result extrapolation of the proximal iterates by a factor ${c+1\over c}$ results in convergence acceleration at negligible computational cost. We also  characterize the region of extrapolation factors that lead to acceleration of convergence, and show that it is an interval that contains $(1,1+1/c]$, but may potentially be substantially larger. These facts are new to the author's knowledge, and they are somewhat unexpected as they do not seem to readily admit an intuitive explanation.

\topinsert
\centerline{\hskip2pc\includegraphics[width=5.5in]{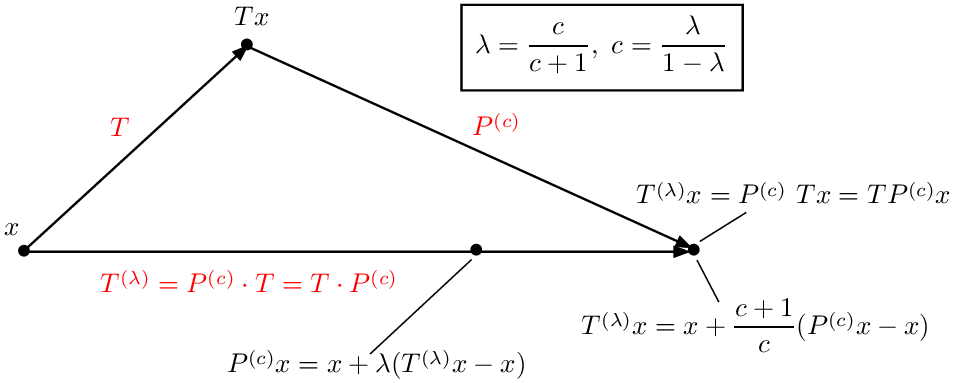}}
\vskip-1.0pc
\hskip-4pc\fig{0pc}{\figoneone.} {Relation of the mappings $\pb$ and $\td$. The mapping $\pb$ is obtained by interpolation  between $\td$ and the identity. Reversely, $\td$ is an extrapolated form of $\pb$.}\endinsert

Aside from its conceptual value and its acceleration potential, the relation between $\pb$ and $\td$ suggests the possibility of new algorithmic approaches for large scale applications where the proximal algorithm can be used conveniently.
In particular, one may consider the projected proximal algorithm,
 $$x_{k+1}=\ppb x_k,$$
which aims to converge to a fixed point of $\ppb$. 
The algorithm may be based on simulation-based computations of $\ptd x$, and such computations have been discussed in the approximate DP context as part of the LSPE($\l$) method (noted earlier), and the $\l$-policy iteration method (proposed in [BeI96], and further developed in the book [BeT96], and the papers [Ber12b] and [Sch13]). The simulation-based methods for computing $\ptd x$ have been adapted to the more general linear equation context in [BeY07], [BeY09]; see also [Ber12a], Section 7.3. Another possibility is to use simulation-based matrix inversion to solve the fixed point equation $x=\ptd x$. In the approximate DP context this is the LSTD($\l$) method, which has also been extended to the general linear equations context in [BeY09].

In Section 3 of this paper, we selectively summarize without much analysis how to adapt and transfer algorithms between the TD/approximate DP and the proximal contexts. Our aim is to highlight the algorithmic possibilities that may allow us to benefit from the accumulated implementation experience within these contexts. To this end, we will draw principally on the individual and joint works of the author, M.\ Wang, and H.\ Yu; see [BeY07], [BeY09], [YuB10], [Yu10], [Yu12], [Ber11a], [Ber11b], [YuB12], [WaB13], [WaB14], and the textbook account of [Ber12a], Section 7.3, where extensions and analysis of TD($\l$), LSTD($\l$), and LSPE($\l$) for solution of the general linear system $x=\ptd x$ were given.
This includes criteria for $\td$ and $\ptd$ to be a contraction, error bounds, simulation-based implementations, algorithmic variations, dealing with singularity or near singularity of $\ppb$, etc. 

In Section 4 we extend our analysis and algorithms of Section 2 to nonlinear fixed point problems. In particular, we show that an extrapolated form of the proximal algorithm provides increased reduction of the distance to the fixed point over the standard proximal algorithm, provided the fixed point problem has a unique solution and involves a nonexpansive mapping (cf.\ Assumption \assumptiona). In Section 4, we also consider forward-backward splitting algorithms and provide a natural  generalization of the extrapolation ideas. To our knowledge, these are the first simple extensions of the proximal and forward-backward algorithms for major classes of nonlinear problems, which guarantee acceleration. Other extrapolation methods, such as the ones of [Ber75] and [Ber82], Section 2.3.1 (for convex optimization), or [EcB92] (for monotone operator problems), guarantee convergence but not acceleration, in the absence of additional prior knowledge. 

The convergence theory of temporal difference methods is restricted to linear systems that satisfy Assumption \assumptiona. Thus, for nonlinear fixed point problems, the connection of temporal difference and proximal algorithms seems considerably weaker. To address this situation, we introduce in Section 5 algorithmic ideas based on linearization whereby $T$ is linearized at each iterate $x_k$, and  the next iterate $x_{k+1}$ is obtained with a temporal differences-based (exact, approximate, or extrapolated) proximal iteration using the linearized mapping. This approach is similar to Newton's method for solving nonlinear fixed point problems, where the linearized system is solved exactly, rather than approximately (using a single proximal iteration). It is also related to earlier work by Bertsekas and Yu [BeY10], [BeY12], [YuB13] on distributed asynchronous policy iteration algorithms in exact and approximate DP. Let us also mention that minimization algorithms with a composite objective function have been addressed with a prox-linear approach, whereby at each iteration the (differentiable) inner composite function is first linearized and a proximal iteration is applied to the outer composite function; see Lewis and Wright [LeW08], and subsequent works, including the recent papers by Drusvyatskiy and Lewis [DrL16], Duchi and  Ryan [DuR17], and the references quoted there. However, these methods deal with minimization rather than fixed point problems and do not seem to be related to the linearized proximal methods of this paper.

\vskip-1.5pc

\section{Interpolation and Extrapolation Formulas}
 
\xdef\proptwoone{\propn}\propnum\show{myproposition}

 \pn We first review a known result from [BeY09] (Prop.\ 3) regarding the multistep mapping $\td$. By repeatedly applying the formula $x=Ax+b$, we can verify that
$$\td x=A^{(\l)}x+b^{(\l)},\xdef\alambdaform{\lab}\eqnum\show{oneo}$$ 
where 
$$A^{(\l)}= (1-\l)\sum_{\ell=0}^\infty \l^\ell A^{\ell+1},\ \qquad  b^{(\l)}= \sum_{\ell=0}^\infty\l^\ell A^\ell b,\xdef\seriesformula{\lab}\eqnum\show{oneo}
$$
assuming that the series above are convergent.
The following proposition shows that under Assumption \assumptiona, $\td$ is well defined by the power series $(1-\l)\sum_{\ell=0}^\infty \l^{\ell}T^{\ell+1}$ [cf.\ Eq.\ \tddef], and that it is a contraction with respect to some norm.

%\texshopbox
{\proposition{\proptwoone:} Let Assumption \assumptiona\ hold, and let $\l\in (0,1)$.
\nitem{(a)} The matrix $A^{(\l)}$ and the vector $b^{(\l)}$ are well-defined in the sense that the series in Eq.\ \seriesformula\ are convergent.%}\texshopboxnt{
\nitem{(b)} The eigenvalues of $A^{(\l)}$ have the form 
$$\theta_i=(1-\l)\sum_{\ell=0}^\infty \l^\ell \zeta_i^{\ell+1}={\zeta_i(1-\l)\over 1-\z_i\l},\qquad i=1,\ldots,n,\xdef\eigenform{\lab}\eqnum\show{oneo}$$%}\texshopboxnt{
where $\z_i$, $i=1,\ldots,n,$ are the eigenvalues of $A$. Furthermore, we have
$\sigma(A^{(\l)})<1$ and  $\lim_{\l\to1}\sigma\bl(A^{(\l)}\br)=0.$
}

\xdef\proptwotwo{\propn}\propnum\show{myproposition}
 
The property $\sigma(A^{(\l)})<1$ asserted in the preceding proposition, is critical for the subsequent development and depends on the eigenvalues of $A$ being different than 1 (cf.\ Assumption \assumptiona). For an intuitive explanation, note that the eigenvalues of $A^{(\l)}$ can be viewed as convex combinations of complex numbers from the unit circle at least two of which are different from each other since $\zeta_i\ne1$ [the nonzero corresponding eigenvalues of $A$ and $A^2$ are different from each other, cf.\ Eqs.\ \alambdaform, \eigenform]. As a result the eigenvalues of $A^{(\l)}$ lie within the interior of the unit circle under Assumption \assumptiona.

 The relation between the proximal mapping $\pb$ and the multistep mapping $\td$ is established in the following proposition, which is illustrated in Fig.\ \figoneone.

%\texshopbox
{\proposition{\proptwotwo:} Let Assumption \assumptiona\ hold, and let $c>0$ and  $\l={c\over c+1}.$
Then:
\nitem{(a)} $\pb$ is given by 
$$\pb=(1-\l)\sum_{\ell=0}^\infty \l^{\ell}T^\ell,\xdef\pbpower{\lab}\eqnum\show{oneo}$$
 and can be written as%}\texshopboxnt{\item{}
$$\pb x=\ol A^{(\l)} x+ \ol b^{(\l)},\qquad x\in\rn,\xdef\alambdaformbar{\lab}\eqnum\show{oneo}$$
where 
$$\ol A^{(\l)}= (1-\l)\sum_{\ell=0}^\infty \l^\ell A^{\ell},\ \qquad  \ol b^{(\l)}= \sum_{\ell=0}^\infty\l^{\ell+1} A^\ell b.\xdef\seriesformulabar{\lab}\eqnum\show{oneo}
$$
\nitem{(b)} We have
$$\td= T\pb=\pb T,\xdef\tdpb{\lab}\eqnum\show{oneo}$$
and for all $x\in\rn$,
$$\pb x=(1-\l)x+\l \td x,\qquad \td x=-{1\over c} x+{c+1\over c}\pb x,\xdef\interp{\lab}\eqnum\show{oneo}$$%}\texshopboxnt{\item{}
or equivalently
$$\pb x=x+\l \big(\td x-x\big),\qquad \td x=x+{c+1\over c}\big(\pb x-x\big).\xdef\interpt{\lab}\eqnum\show{oneo}$$
}

 \xdef \figtwoone{\figr}\figrnum\show{myfigure}

\xdef \figtwotwo{\figr}\figrnum\show{myfigure}

\proof (a) The inverse in the definition of $\pb$ [cf.\ Eq.\ \pbformula] is written as
$$\lf({c+1\over c}I-A\ri)^{-1}=\left({1\over \l}I -A\right)^{-1}=\l(I-\l A)^{-1}=\l \sum_{\ell=0}^\infty (\l A)^{\ell},$$
where the power series above is convergent by Prop.\ \proptwoone(a). Thus, from Eq.\ \pbformula\ and the equation ${1\over c}={1-\l\over \l}$,
$$\pb x=\lf({c+1\over c}I-A\ri)^{-1}\lf(b+{1\over c}x\ri)=\l\sum_{\ell=0}^\infty (\l A)^{\ell}\lf (b+{1-\l\over \l} x\ri)=(1-\l)\sum_{\ell=0}^\infty (\l A)^{\ell}x+\l\sum_{\ell=0}^\infty (\l A)^{\ell}b,$$
which from Eq.\ \seriesformulabar, is equal to $\ol A^{(\l)}x+\ol b^{(\l)}$, thus proving Eq.\ \alambdaformbar.
\smskip
\pn (b) We have, using Eqs.\ \alambdaform, \seriesformula, \alambdaformbar\ and \seriesformulabar,
$$T\pb x=A\big(\ol A^{(\l)} x+ \ol b^{(\l)}\big)+b=(1-\l)\sum_{\ell=0}^\infty \l^\ell A^{\ell+1}x+\sum_{\ell=0}^\infty\l^{\ell+1} A^{\ell+1} b+b=A^{(\l)}x+b^{(\l)}=\td x,$$
thus proving the left side of Eq.\ \tdpb. The right side is proved similarly.
The interpolation/extrapolation formulas \interp\   and \interpt\ follow by a straightforward calculation from Eq.\ \pbpower\ and the definition $\td =(1-\l)\sum_{\ell=0}^\infty \l^{\ell}T^{\ell+1}$ [cf.\ Eq.\ \tddef]. As an example, the following calculation shows the left side of Eq.\ \interpt:
$$\eqalignno{x+\l \big(\td x-x\big)&=(1-\l)x+\l \td x\cr
&=(1-\l)x+\l \left((1-\l)\sum_{\ell=0}^\infty \l^\ell A^{\ell+1}x+\sum_{\ell=0}^\infty \l^\ell A^{\ell}b\right)            \cr
&=(1-\l)\left(x+\sum_{\ell=1}^\infty \l^\ell A^{\ell}x\right)+\sum_{\ell=0}^\infty \l^{\ell+1} A^{\ell}b            \cr
&=\ol A^{(\l)}x+\ol b^{(\l)}\cr
&=\pb x.\cr}$$
{\bf Q.E.D}
\smskip

We will now use the extrapolation formulas of Prop.\ \proptwotwo(b) to construct variants of the proximal algorithm with interesting properties. The next proposition establishes the convergence and convergence rate properties of the proximal and multistep iterations, and shows how the proximal iteration can be accelerated by extrapolation or interpolation.

\xdef\proptwothree{\propn}\propnum\show{myproposition}
 
%\texshopbox
{\proposition{\proptwothree:} Let Assumption \assumptiona\ hold, and let $c>0$ and  $\l={c\over c+1}.$
Then the eigenvalues of $\ol A^{(\l)}$ are
$$\ol \theta_i={1-\l\over 1-\z_i\l},\qquad i=1,\ldots,n,\xdef\eigenformbar{\lab}\eqnum\show{oneo}$$%}\texshopboxnt{\item{} 
where $\z_i$, $i=1,\ldots,n,$ are the eigenvalues of $A$. Moreover, $A^{(\l)}$  and $\ol A^{(\l)}$  have the same eigenvectors. Furthermore, we have
$${\sigma(A^{(\l)})\over \sigma(A)}\le \sigma(\ol A^{(\l)})<1,\xdef\specrad{\lab}\eqnum\show{oneo}$$
so $\sigma(A^{(\l)})<\sigma(\ol A^{(\l)})$ if $\sigma(A)<1$.
}

\proof  Let $e_i$ be an eigenvector of $A^{(\l)}$ corresponding to the eigenvalue $\theta_i$. By using the interpolation formula \interp\ and the eigenvalue formula \eigenform\ for $\theta_i$, we have
$$\ol A^{(\l)}e_i=(1-\l)e_i+\l A^{(\l)}e_i=\big((1-\l)+\l\theta_i\big)e_i=\lf((1-\l)+\l{\zeta_i(1-\l)\over 1-\z_i\l}\ri)e_i={1-\l\over 1-\z_i\l}e_i.$$
Hence, $\ol \theta_i={1-\l\over 1-\z_i\l}$ and $e_i$ are the corresponding eigenvalue and eigenvector of $\ol A^{(\l)}$, respectively.

\topinsert
\centerline{\includegraphics[width=4.4in]{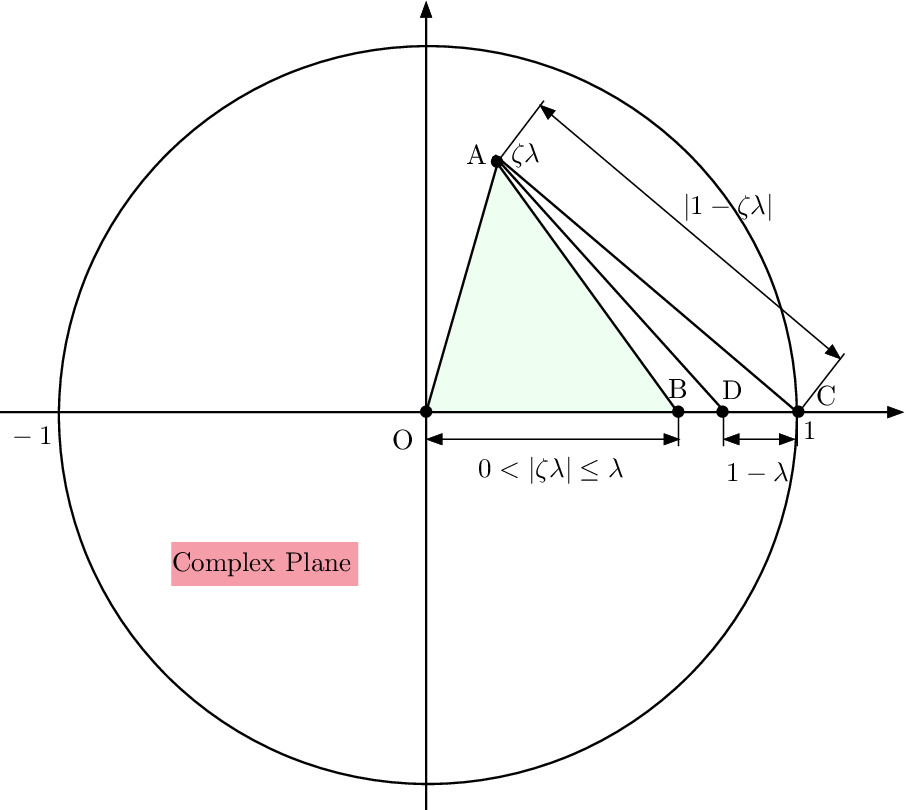}}
\vskip-1pc
\fig{0pc}{\figtwoone.} {Proof of the inequality $|\ol \theta_i|<1$, or equivalently that $1-\l<|1-\zeta \l|$ for all complex numbers $\zeta\ne 1$ with $|\zeta|\le 1$, and $\l\in(0,1)$. We consider the unit circle of the complex plane and the complex number $\zeta \l$, and we note that $0<|\zeta\l |\le \l<1$. If $\zeta$ is in the left-hand side of the plane or on the vertical axis, we clearly have 
$$1-\l<1\le |1-\zeta \l|,$$
so it is sufficient to consider the case where $\zeta\ne0$ and the real part of $\zeta$ is positive, which is depicted in the figure. If $\z$ is real, we have $\z>0$ as well as $\l>0$, so 
$$|1-\zeta \l|=1-\zeta \l >1-\l,$$
 and we are done. If $\z$ is not real, we consider the isosceles triangle OAB (shaded in the figure), and note that the angles of the triangle bordering the side AB are less than 90 degrees. It follows that the angle ABC and hence also the angle ADC shown in the figure is greater than 90 degrees. Thus the side AC of the triangle ADC is strictly larger than the side DC. This is equivalent to the desired result $1-\l<|1-\zeta \l|$.%\vskip1pc
}\endinsert

The proof that $\sigma(\ol A^{(\l)})<1$, or equivalently that $|\ol \theta_i|<1$ for all $i$, follows from a graphical argument on the complex plane, which is given in the caption of Fig.\ \figtwoone. We also have from Eqs.\ \eigenform\ and \eigenformbar\
$$|\ol \theta_i|={|\theta_i|\over |\zeta_i|},\qquad i=1,\ldots,n,$$
which implies that
$$|\ol \theta_i|\ge {|\theta_i|\over \sigma(A)},\qquad i=1,\ldots,n.$$
By taking the maximum of both sides over $i$, we obtain the left side of Eq.\ \specrad.
 \qed
 
 \xdef\proptwofour{\propn}\propnum\show{myproposition}
 
An interesting conclusion can  be drawn from Prop.\ \proptwothree\ about the convergence and the rate of convergence of the proximal iteration $x_{k+1}=\pb x_k$ and  the multistep iteration $x_{k+1}=\td x_k$.
Under Assumption \assumptiona, {\it both iterations are convergent, but the multistep iteration is faster when $A$ is itself a contraction} (with respect to some norm) and is not slower otherwise; cf.\ Prop.\ \proptwothree(a). In the case where $A$ is not a contraction [$\sigma(A)=1$] it is possible that $\sigma(A^{(\l)})=\sigma(\ol A^{(\l)})$ (as an example consider a case where all the eigenvalues $\zeta_i$ have modulus 1). 

Even in the case where $\sigma(A^{(\l)})=\sigma(\ol A^{(\l)})$, however, it is possible to accelerate the proximal iteration by interpolating strictly between $\pb x_k$ and $\td x_k$. This is shown in the next proposition, which establishes the convergence rate properties of the extrapolated proximal iteration, and quantifies the range of extrapolation factors that lead to acceleration.

%\texshopbox
{\proposition{\proptwofour} Let Assumption \assumptiona\ hold, and let $c>0$ and  $\l={c\over c+1}.$ Consider any iteration that extrapolates  from $\pb$ in the direction of $\td$, i.e.,
$$x_{k+1}=(1-\g)\pb x_k+\g \td x_k,\qquad \g>0,\xdef\interpolatediter{\lab}\eqnum\show{oneo}$$%}\texshopboxnt{\pn 
and write it in matrix-vector form as
$$x_{k+1}=A(\l,\g) x_k+b(\l,\g),$$%}\texshopboxnt{\pn
where $A(\l,\g)$ is an $n\times n$ matrix and $b(\l,\g)\in\rn$. The eigenvalues of $A(\l,\g)$ are given by
$$\theta_i(\g)=(1-\g) \ol \theta_i+\g\theta_i,\qquad i=1,\ldots,n,\xdef\eigenconvex{\lab}\eqnum\show{oneo}$$
and we have
$$\sigma\big(A(\l,\g)\big)<\sigma(\ol A^{(\l)}),\xdef\spectralineq{\lab}\eqnum\show{oneo}$$%}\texshopboxnt{\pn
for all $\g$ in the interval $(0,\g_{\eightpoint{max}})$, where
$$\g_{\eightpoint{max}}=\max\Big\{\g>0\ \big|\ \big|\theta_i(\g)\big|\le\ol\theta_i,\ \forall\ i=1,\ldots,n\Big\}.$$
Moreover, we have $\g_{\eightpoint{max}}\ge 1$, with equality holding if and only if $\s(A)=1$. 
}

\proof  The eigenvalue formula \eigenconvex\ follows from the interpolation formula
$$A(\l,\g)=(1-\g)\ol A^{(\l)}+ \g A^{(\l)},$$
and the fact that $A^{(\l)}$  and $\ol A^{(\l)}$  have the same eigenvectors (cf.\ Prop.\ \proptwothree). For each $i$, the scalar 
$$\max\Big\{\g>0\ \big|\ \big|\theta_i(\g)\big|\le |\ol \theta_i|\Big\}$$
 is the maximum extrapolation factor $\g$ for which $\theta_i(\g)$ has at most as large modulus as $\ol \theta_i$ (cf.\ Fig.\ \figtwotwo), and the inequality \spectralineq\ follows. The inequality $\g_{\eightpoint{max}}\ge 1$ follows from the construction of Fig.\ \figtwotwo, since $|\theta_i|\le |\ol \theta_i|$, $\theta_i\ne \ol \theta_i$,  and $\g=1$ corresponds to the iteration $x_{k+1}=\td x_k$. Finally, we have $\g_{\eightpoint{max}}=1$ if and only if $|\theta_i|= |\ol \theta_i|$ for some $i$, which happens if and only if $|\zeta_i|=1$ for some $i$, i.e.,  $\s(A)=1$.
\qed

\topinsert
\centerline{\includegraphics[width=4.3in]{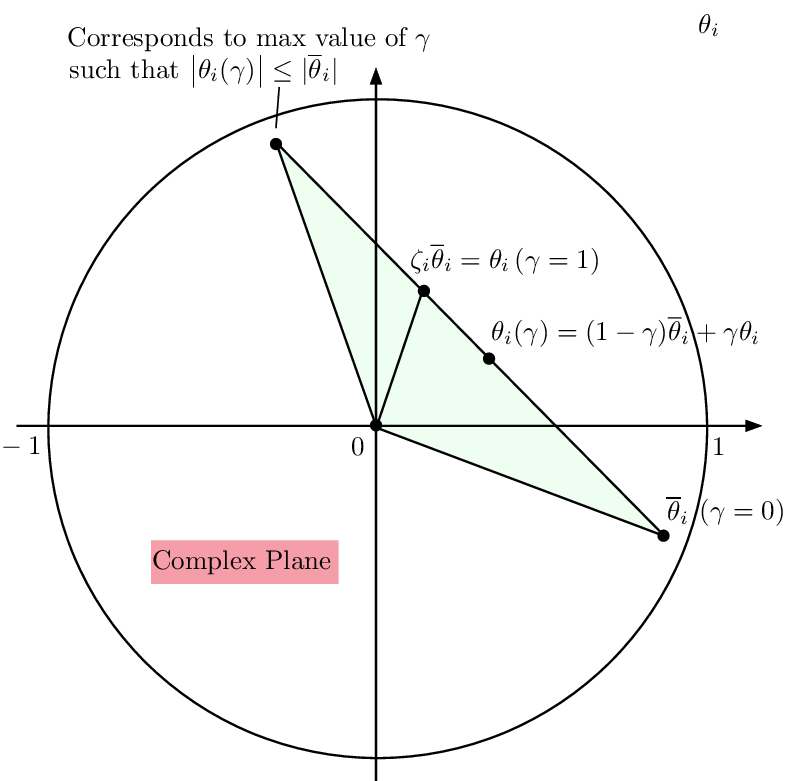}}
%\vskip-16pc
\fig{0pc}{\figtwotwo.} {Illustration of the proof of Prop.\ \proptwofour. The eigenvalues $\theta_i(\g)$ of $A(\l,\g)$ are linear combinations (with $\g>0$) of the eigenvalues $\ol \theta_i$ and $\theta_i=\zeta_i\ol\theta_i$ of $\ol A^{(\l)}$ and $A^{(\l)}$, respectively, and we have $\theta_i\le \ol \theta_i$.}\endinsert

We may implement the extrapolation/interpolation iteration \interpolatediter\ by first implementing the proximal iteration $x_{k+1}=\pb x_k$ and then the multistep iteration according to
$$x_{k+1}=\td x_k=x_k+{c+1\over c} \big(\pb x_k-x_k\big).$$
In this way, unless $\sigma(A)=1$ [cf.\ Eq.\ \specrad], we achieve acceleration over the proximal iteration. We may then aim for greater acceleration by extrapolating or interpolating  between $\pb x_k$ and $\td x_k$ with some  factor, possibly determined by experimentation [strict acceleration can always be achieved with $\g\in(0,1)$]. This provides a simple and reliable method to accelerate the convergence of the proximal algorithm without knowledge of the eigenvalue structure of $A$ beyond Assumption \assumptiona.\footnote{\dag}{\ninepoint  It is well known that the proximal iteration can be extrapolated by a factor of as much as two while maintaining convergence. This was first shown for the special case of a convex optimization problem by the author in [Ber75], and  for the general case of finding a zero of a monotone operator in Eckstein and Bertsekas [EcB92]; see also a more refined analysis, which quantifies the effects of extrapolation, for the case of a quadratic programming problem, given in the book [Ber82], Section 2.3.1. However, we are not aware of any earlier proposal of a simple and general scheme to choose an extrapolation factor that maintains convergence {\it and} also guarantees acceleration. Moreover, this extrapolation factor, $(c+1)/c$, may be much larger than two.} Conversely, we may implement the proximal iteration by interpolating the multistep iteration. Some of the possibilities along this direction will be reviewed in the next section.

Finally, let us show that the multistep and proximal iterates $\pb x_k$ and $\td x_k$ can be computed by solving fixed point problems involving a contraction of modulus $\l \s(A)$.

 \xdef\proptwofive{\propn}\propnum\show{myproposition}
 
%\texshopbox
{\proposition{\proptwofive} Let Assumption \assumptiona\ hold, and let $c>0$ and  $\l={c\over c+1}.$ The multistep and proximal iterates $\td x_k$ and $\pb x_k$ are the unique fixed points of the contraction mappings $W_{x_k}$ and $\ol W_{x_k}$  given by%}\texshopboxnt{\pn
$$W_{x_k}x=(1-\l)T x_k+\l T x,\qquad x\in \rn,$$%}\texshopboxnt{\pn
and
$$\ol W_{x_k}x=(1-\l)x_k+\l T x,\qquad x\in \rn,$$
respectively.
}

\proof Clearly $W_{x_k}$ and $\ol W_{x_k}$ are contraction mappings, since they are linear with spectral radius $\l \s(A)\le \l<1$. To show that $\td x_k$ is the fixed point of  $W_{x_k}$, we must verify that $\td x_k=W_{x_k}\bl(\td x_k\br)$, or equivalently that
$$\td x_k = (1-\l)Tx_k+\l T\bl(\td x_k)=(1-\l)Tx_k+\l \td (Tx_k)\xdef\firsteq{\lab}\eqnum\show{oneo}$$
[here we are applying the formula $T\bl(\td x)=\td (Tx)$, which is easily verified using Eqs.\ \alambdaform\ and \seriesformula]. In view of the interpolation formula
$$(1-\l)x+\l \td x=\pb x,\qquad \forall\ x\in\rn,\xdef\pcinterp{\lab}\eqnum\show{oneo}$$
[cf.\ Eq.\ \interp],
the right-hand side of Eq.\ \firsteq\ is equal to $P^{(c)}(Tx_k)$,
which from the formula
$\td= \pb T$
[cf.\ Eq.\ \tdpb], is equal to $\td x_k$, the left-hand side of Eq.\ \firsteq.

Similarly, to show that $\pb x_k$ is the fixed point of  $\ol W_{x_k}$, we must verify that $\pb x_k=\ol W_{x_k}(\pb x_k)$, or equivalently that
$$\pb x_k = (1-\l)x_k+\l T\bl(\pb x_k\br).$$
This is proved by combining the formula
$\td= T\pb$
[cf.\ Eq.\ \tdpb], and the interpolation formula \pcinterp. \qed

The fixed point property of the preceding proposition states that $\td x$ is the unique solution of the following equation in $y$:
$$y=(1-\l)Tx+\l Ty,$$
or
$$y=(1-\l)(Ax+b)+\l (Ay+b),$$
and thus provides an explicit matrix inversion formula for the multistep mapping $\td$:  
$$\td x=(1-\l A)^{-1}\big(b+(1-\l)Ax\big).\xdef\tdmapformula{\lab}\eqnum\show{oneo}$$
This formula should be compared with the formula
\pbformula\ for the proximal mapping, which can be written in terms of $\l$ as
$$P^{(c)}x=(1-\l A)^{-1}\big(\l b+(1-\l)x\big).$$
The fact that the multistep iterate $x_{k+1}=\td x_k$ is the fixed point of $W_{x_k}$ is known in exact and approximate DP, and forms the basis for the $\l$-policy iteration method; see [BeI96], [BeT96], Section 2.3.1. This is a variant of policy iteration where policy evaluation is done by performing a single multistep iteration using the mapping $\td$, where $A$ corresponds to the policy being evaluated. In fact the formula \tdmapformula\ is given in [Ber18], Section 4.3.3. In view of our analysis in this paper, it follows that $\l$-policy iteration is the approximate version of exact policy iteration, where the exact policy evaluation phase of the latter (which is to find the fixed point of $T$), is approximated with a single (extrapolated) iteration of the proximal algorithm. The $\l$-policy iteration method admits some interesting simulation-based implementations, which have been discussed in the approximate DP literature ([Ber12b], [Sch13]), but will not be discussed further here.  Based on Prop.\ \proptwofive, the proximal iteration $x_{k+1}=\pb x_k$ admits similar implementations. 

Proposition \proptwofive\ also suggests the iteration
$$x_{k+1}=V_m x_k,\xdef\altapprox{\lab}\eqnum\show{oneo}$$
where $V_m x_k$ is obtained by $m>1$ iterations of the mapping $W_{x_k}$ starting with $x_k$, i.e.,
$$V_mx_k=(W_{x_k})^m x_k,$$
so $V_mx_k$ is an approximate evaluation of $\td x_k$,  the fixed point of $W_{x_k}$.
It can be verified by induction that
$$V_mx_k=(1-\l)(Tx_k+\l T^2 x_k+\cdots+\l^{m-1}T^m x_k)+\l^m T^m x_k,$$
and that $V_m$ is a contraction mapping [the preceding formula is given in [BeT96], Prop.\ 2.7(b), while the contraction property of $V_m$ is proved similar to Prop.\ 3(a) of [BeY09]]. There is also the similar iteration $x_{k+1}=\ol V_m x_k,$
where $\ol V_m x_k$ is obtained by $m>1$ iterations of the mapping $\ol W_{x_k}$ starting with $x_k$. This iteration may be viewed as an iterative approximate implementation of the proximal algorithm that does not require matrix inversion.

\vskip-1pc

\section{Projected Proximal, Proximal Projected, and Temporal \hfill\break Difference Algorithms}

 \xdef \figonetwo{\figr}\figrnum\show{myfigure}

\pn In this section we aim to highlight situations where analysis and experience from approximate DP can be fruitfully transferred to the solution of linear equations by proximal-like algorithms. In particular, we discuss the  simulation-based approximate solution of the equation $x=Tx$ within a subspace $S$ spanned by a relatively small number of basis functions $\phi_1,\ldots,\phi_s\in \rn$. Note that while in DP the matrix $A$ is assumed substochastic, the methodology described in this section assumes only Assumption \assumptiona. The extension beyond the DP framework was first given in [BeY09] and described in further detail in Section 7.3 of the textbook [Ber12a] under the name ``Monte-Carlo Linear Algebra."  

We denote by $\Phi$ the $n\times s$ matrix whose columns are $\phi_1,\ldots,\phi_s$, so $S$ can be represented as
$$S=\{\Phi r\mid r\in \re^s\}.$$
Instead of solving $x=Tx$ or the multistep system $x=\td x$ (which has the same solution as $x=Tx$) we consider the projected form
$$x=\ptd x,$$
where the mapping $\P:\rn\mapsto\rn$ is some form of projection onto $S$, in the sense that $\P$
is linear, $\P x\in S$ for all $x\in\rn$, and $\P x=x$ for all $x\in S$.

A general form of such $\P$ is the oblique projection
$$\P=\Phi(\Psi'\Xi\Phi)^{-1}\Psi'\Xi,\xdef\projoblique{\lab}\eqnum\show{oneo}$$
where $\Xi$ is a diagonal $n\times n$ positive semidefinite matrix with components $\xi_1,\ldots,\xi_n$ along the diagonal, and $\Psi$ is an $n\times s$ matrix such that $\Psi'\Xi\Phi$ is invertible. Most often in  approximate DP applications, the projection is orthogonal with respect to a Euclidean norm [cf.\ case (1)] below. However, there are situations where a more general type of projection is interesting [see cases (2) and (3) below]. The use of oblique (as opposed to orthogonal) projections in approximate DP was suggested by Scherrer [Sch10].
We note some special cases that have received attention in the approximate DP literature:

\nitem{(1)} {\it Orthogonal projection\/}: Here $\Psi=\Phi$ and $\Xi$ is positive definite. Then $\P$ is the orthogonal projection onto  $S$ with respect to the norm corresponding to $\Xi$, i.e., $\|x\|^2=\sum_{i=1}^n\xi_i x_i^2$. Solving the projected system $x=\ptd x$ corresponds to a form of Galerkin approximation, as noted earlier.\footnote{\dag}{\ninepoint An important fact here is that if $\s(A)<1$ there exists a (weighted) orthogonal projection $\P$ such that $\ptd$ is a contraction for all $\l\in(0,1)$ (see [BeY09] for a more detailed discussion of this issue). While finding such $\P$ may not be simple in general, in approximate DP, a suitable projection $\P$ is implicitly known in terms of the stationary distribution of a corresponding Markov chain (see e.g., [BeT96], Lemma 6.4, or [Ber12a], Lemma 6.3.1).  Another important fact is that if $\s(A)\le 1$, the mapping $\ptd$ can be made to be a contraction with respect to any norm by taking $\l$ sufficiently close to 1 [cf.\ Prop.\ \proptwoone(b)].}

\nitem{(2)} {\it Seminorm projection\/}: Here $\Psi=\Phi$, the matrix $\Phi'\Xi\Phi$ is invertible, but $\Xi$ is only positive semidefinite (so some of the components $\xi_i$ may be zero). Then $\P$ is a seminorm projection with respect to the seminorm defined by $\Xi$. The seminorm projection $\P x$ can be computed as the unique solution of the least squares problem
$$\min_{r\in \re^s}\sum_{i=1}^n\xi_i(x_i-\phi_i'r)^2,\xdef\semileastsq{\lab}\eqnum\show{oneo}$$
where $\phi_1',\ldots,\phi_n'$ are the rows of $\Phi$ (the solution is unique in view of the assumed invertibility of $\Phi'\Xi\Phi$). This type of least squares problem arises when we are trying to approximate a high dimensional vector $x$ onto $S$ but we know only some of the components of $x$ (but enough to ensure that the preceding minimization has a unique solution). Seminorm projections were first considered in the approximate DP context in the paper by Yu and Bertsekas [YuB12], and we refer to that reference for more discussion on their applications.

\nitem{(3)} {\it Aggregation\/}: Here $\P=\Phi D$, where $D$ is an $s\times n$ matrix, and we assume that the rows of $\Phi$ and $D$ are probability distributions, and that $\Phi$ and $D$ have a full rank $s$. We replace the solution of the system $x=Ax+b$ with the projected system $x=\P(Ax+b)$, or $\Phi r=\Phi D(A\Phi r +b)$, which equivalently (since $\phi$ has full rank) can be written as
$$r=DA\Phi r+Db.$$
This system is obtained by forming convex combinations of rows and columns of $A$ to construct the ``aggregate" matrix $DA\Phi$. Aggregation has a long history in numerical computation, and it is discussed in detail in the context of approximate DP in [Ber11b] and [Ber12a] (Sections 6.5 and 7.3.7), and the references quoted there. It turns out that for a very broad class of aggregation methods, called ``aggregation with representative features" ([Ber12a], Example 6.5.4, and Exercise 7.11), the matrix $\Phi D$ is a seminorm projection, as first shown in [YuB12]. The matrix $\Phi D$ can also be viewed as an oblique projection in this case (see [Ber12a], Section 7.3.6). 
\smskip

While the solution of $x=\td x$ is the fixed point $x^*$ of $T$ for all $\l$, the solution of the projected equation $x=\ptd x$ depends on $\l$. Also, because of the linearity of $\P$ and the extrapolation property \interp\ shown in the preceding section, the projected proximal equation
$x=\ppb x$
has the same solution as the system $x=\ptd x$ where $\l={c\over c+1}$. Let us denote by $x_\l$ this solution. Generally, for any norm $\|\cdot\|$ we have  the error bound
$$\|x^*-x_\l\|\le \bl\|(I-\Pi \ala)^{-1}\br\|\,\|x^*-\Pi x^*\|,\xdef\basicbound{\lab}\eqnum\show{oneo}$$
which is derived from the following  calculation:
$$x^*-x_\l= x^*-\Pi x^*+\Pi x^*-x_\l=x^*-\Pi x^*+\Pi Tx^*-\Pi \td x_\l
=x^*-\Pi x^*+\Pi \ala(x^*-x_\l).$$
Thus the approximation error $\|x^*-x_\l\|$ is proportional to the ``distance" $\|x^*-\Pi x^*\br\|$ of the solution $x^*$ from the approximation subspace. It is well known that for values of $\l$ close to 1, $x_\l$  tends to approximate better the projection $\P x^*$. However, values of $\l$ close to 1 correspond to large values of $c$, resulting in a less-well conditioned projected proximal equation $x=\ppb x$.
There is also a related tradeoff that is well-documented in the DP literature: as $\l$ is increased towards 1,  solving the projected equation $x=\ptd x$ by simulation is more time-consuming because of increased simulation noise and an associated need for more simulation samples.
For further discussion of the choice of $\l$, we refer to the references cited earlier, and for error bounds that are related but sharper than Eq.\ \basicbound, we refer to Yu and Bertsekas [YuB10], and Scherrer [Sch10], [Sch13]. 

For the case of the projection formula \projoblique\ and assuming that $\Phi$ has full rank, the solution $x_\l$ of the system
$$x=\ptd x=\P\bl(\ala x+b\br)=\Phi(\Psi'\Xi\Phi)^{-1}\Psi'\Xi\big(\ala x+\bla\big),$$
can be written as
$$x_\l=\Phi r_\l,$$
where $r_\l$ is the unique solution of the low-dimensional system of equations
$$r=(\Psi'\Xi\Phi)^{-1}\Psi'\Xi\big(\ala \Phi r+\bla\big).$$
Equivalently, this system is written as  
$$r=Q^{(\l)}r+d^{(\l)},\xdef\lowdimsystlaz{\lab}\eqnum\show{oneo}$$
where
$$Q^{(\l)}=(\Psi'\Xi\Phi)^{-1}\Psi'\Xi \ala \Phi,\qquad d^{(\l)}=(\Psi'\Xi\Phi)^{-1}\Psi'\Xi b^{(\l)}.\xdef\lowdimsystlao{\lab}\eqnum\show{oneo}$$
By defining 
$$C^{(\l)}=I-Q^{(\l)},\xdef\lowdimsystlac{\lab}\eqnum\show{oneo}$$
 this system can also be written as
$$C^{(\l)}r=d^{(\l)},\xdef\lowdimsystla{\lab}\eqnum\show{oneo}$$
where 
$$C^{(\l)}=(\Psi'\Xi\Phi)^{-1}\Psi'\Xi\big(I-\ala\big)\Phi,\qquad d^{(\l)}=(\Psi'\Xi\Phi)^{-1}\Psi'\Xi b^{(\l)}.\xdef\lowdimsystlat{\lab}\eqnum\show{oneo}$$

\old{where
$$C^{(\l)}=\Psi'\Xi\big(I-A^{(\l)}\big)\Phi,\qquad d^{(\l)}=(\Psi'\Xi\Phi)^{-1}\Psi'\Xi b^{(\l)}.$$
We also have
$$C^{(\l)}= (1-\l)\sum_{\ell=0}^\infty \l^\ell C^{\ell+1},\ \qquad  d^{(\l)}= \sum_{\ell=0}^\infty\l^\ell C^\ell d.\xdef\seriesformlow{\lab}\eqnum\show{oneo}
$$}

\subsection{Projected Proximal and Proximal Projected Algorithms}

\pn Let us now consider the solution of the projected equation $x=\ptd x$, and the equivalent systems $r=Q^{(\l)}r+d^{(\l)}$ [cf.\ Eq.\ \lowdimsystlaz] and $C^{(\l)}r=d^{(\l)}$ [cf.\ Eq.\ \lowdimsystla]  with proximal-type algorithms, assuming that $\Phi$ has full rank. There are two different approaches here (which coincide when there is no projection, i.e., $S=\rn$ and $\P=I$): 

\nitem{(a)} Use the algorithm
 $$x_{k+1}=\ptd x_k,\xdef\projtdla{\lab}\eqnum\show{oneo}$$
or equivalently [cf.\ Eqs.\ \lowdimsystlaz\ and \lowdimsystlac],
 $$r_{k+1}=Q^{(\l)}r_k+d^{(\l)}=r_k-\big(C^{(\l)}r_k-d^{(\l)}\big).\xdef\lspeapprox{\lab}\eqnum\show{oneo}$$ 
Another possibility is to use the interpolated version, which is the projected proximal algorithm
 $$x_{k+1}=\ppb x_k,\xdef\projprox{\lab}\eqnum\show{oneo}$$
where $c={\l\over 1-\l}$ (cf.\ Fig.\ \figonetwo), or equivalently, based on the interpolation formula \interpt,
 $$r_{k+1}=r_k+\l\big(Q^{(\l)} r_k+d^{(\l)}-r_k\big)=r_k-\l \big(C^{(\l)}r_k-d^{(\l)}\big).\xdef\lspeapproxprox{\lab}\eqnum\show{oneo}$$ 
 
\midinsert
\centerline{\hskip2pc\includegraphics[width=4.2in]{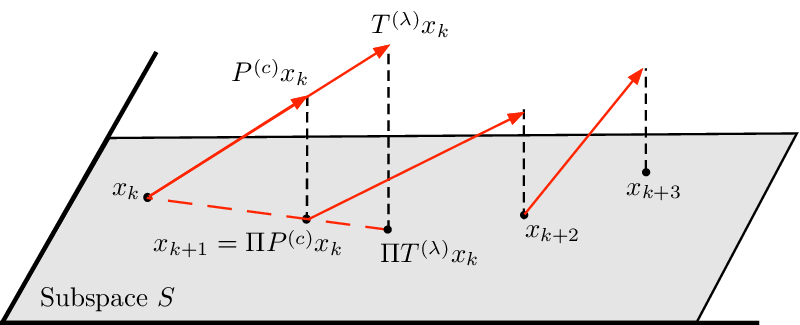}}
\vskip-1pc
\hskip-4pc\fig{0pc}{\figonetwo.} {Illustration of the projected proximal algorithm \projprox\ in relation to the projected multistep iteration \projtdla. All iterates $x_k$, except possibly $x_0$, lie on $S$.}\endinsert

\nitem{(b)} Apply the proximal algorithm \pbformula\ to the system $r=Q^{(\l)}r+d^{(\l)}$ or the system $C^{(\l)}r=d^{(\l)}$ [cf.\ Eqs.\ \lowdimsystlaz\ and \lowdimsystla]:
$$\eqalign{r_{k+1}&=\lf({\hat c+1\over \hat c} I-Q^{(\l)}\ri)^{-1}\lf(d^{(\l)}+{1\over \hat c} r_k\ri)\cr
&=\lf({1\over \hat c} I+C^{(\l)}\ri)^{-1}\lf(d^{(\l)}+{1\over \hat c} r_k\ri)\cr
&=r_k-\lf({1\over \hat c} I+C^{(\l)}\ri)^{-1}\lf(C^{(\l)}r_k-d^{(\l)}\ri),\cr}\xdef\proxpvi{\lab}\eqnum\show{oneo}$$
where ${\hat c}$ is a positive parameter that need not be related to $\l$ [cf.\ Eq.\ \pbformula]. We may also use the extrapolated version that was discussed in the preceding section:
$$r_{k+1}=r_k-{\hat c+1\over \hat c}\lf({1\over \hat c} I+C^{(\l)}\ri)^{-1}\lf(C^{(\l)}r_k-d^{(\l)}\ri),\xdef\proxpviextra{\lab}\eqnum\show{oneo}$$
(cf.\ Fig.\ \figoneone). Extrapolation factors that are intermediate between 1 and ${\hat c+1\over \hat c}$ or larger than ${\hat c+1\over \hat c}$ may be considered. Note that this is a {\it two-parameter algorithm\/}: the  parameter ${\hat c}$ is used for regularization, and may be different from ${1-\l\over \l}$. This allows some flexibility of implementation: the choice of $\l$ should aim to strike a balance between small approximation error $\|x_\l-\P x^*\|$ and implementation difficulties due to ill-conditioning and/or simulation overhead, while the choice of $\hat c$ should aim at guarding against near singularity of $C^{(\l)}$. For the extrapolated algorithm \proxpviextra\ to be valid, not only should $C^{(\l)}$ be invertible, but we must also have $\s\big(I-C^{(\l)}\big)\le 1$. This is true under mild conditions; see [BeY09], Prop.\ 5.

\smskip

The algorithms in (a) and (b) above are related but different. The algorithms in (a) are {\it projected proximal} algorithms (possibly with extrapolation), like the ones in Fig.\ \figonetwo. The algorithms in (b), use the projection and proximal operations in reverse order: they are {\it proximal projected} algorithms, i.e., proximal algorithms applied to the projected equation. Both types of algorithms have the generic form
$$r_{k+1}=r_k-\g G\big(C^{(\l)}r_k-d^{(\l)}\big),$$ 
where $\g$ is a nonnegative stepsize and $G$ is a matrix such that $GC^{(\l)}$ has eigenvalues with positive real parts. This algorithm and its simulation-based implementations, for both cases where $C^{(\l)}$ is nonsingular and singular, has been studied in detail in the papers by Wang and Bertsekas [WaB13] and [WaB14]. Its convergence properties have been analyzed in these references under the assumption that the stepsize $\g$ is sufficiently small. The algorithms \lspeapprox\ and \proxpvi\ have been explicitly noted in these references. The accelerated proximal projected algorithm \proxpviextra\ is new.

Note  that the algorithms \proxpvi\ and \proxpviextra\ make sense also when $\l=0$, in which case 
$$C^{(0)}=(\Psi'\Xi\Phi)^{-1}\Psi'\Xi(I-A)\Phi,\qquad d^{(0)}=(\Psi'\Xi\Phi)^{-1}\Psi'\Xi b,$$
[cf.\ Eq.\ \lowdimsystlat]. 
Then, the algorithm \proxpvi\ (which is known in approximate DP) can be viewed as the proximal algorithm applied to the projected system $x=\P Tx$, while the algorithm \proxpviextra\ can be viewed as the corresponding faster multistep algorithm 
$$r_{k+1}=r_k-{\hat c+1\over \hat c}\lf({1\over \hat c} I+C^{(0)}\ri)^{-1}\lf(C^{(0)}r_k-d^{(0)}\ri),$$
which has not been considered earlier. 

\subsection{Simulation-Based Methods}

\pn The difficulty with the preceding algorithms \lspeapprox-\proxpviextra\ is that when $n$ is very large, the computation of $C^{(\l)}$ and $d^{(\l)}$ involves high-dimensional inner products, whose exact computation is impossible. 
This motivates the replacement of $C^{(\l)}$ and $d^{(\l)}$ with Monte-Carlo simulation-generated estimates, which can be computed with low-dimensional calculations. This simulation-based approach has been used and documented extensively in approximate DP since the late 80s, although the algorithms \lspeapproxprox\ and \proxpviextra\ are new, to our knowledge. Moreover, sampling and simulation for solution of linear systems have a long history, starting with a suggestion
by von Neumann and Ulam (recounted by Forsythe and Leibler [FoL50]); see also Curtiss [Cur54], [Cur57], and the survey by Halton [Hal70].   More recently,
work on simulation methods has focused on using low-order calculations for solving large least squares and other problems. In this connection we note the papers by Strohmer and Vershynin [StV9], Censor, Herman, and Jiang [CeH09], and Leventhal and Lewis [LeL10] on randomized versions of coordinate descent and iterated projection methods for overdetermined least squares problems, and the series of papers by Drineas, Kannan, Mahoney, Muthukrishnan, Boutsidis, and Magdon-Ismail, who consider the use of simulation methods for linear least squares problems and low-rank matrix approximation problems; see [DKM06a], [DKM06b], [DMM06], [DMM08], [DMMS11], and [BDM14].

Let us denote by $C_k^{(\l)}$ and $d_k^{(\l)}$, respectively, the simulation-based estimates to $C^{(\l)}$ and $d^{(\l)}$, which are available at iteration $k$. Then the multistep iteration $x_{k+1}=\ptd x_k$ [cf.\ Eq.\ \projtdla] can be implemented in approximate form as
$$r_{k+1}=r_k-\big(C_k^{(\l)}r_k-d_k^{(\l)}\big),\xdef\lspeapproxsim{\lab}\eqnum\show{oneo}$$ 
[cf.\ Eq.\ \lspeapprox]. We implicitly assume here that at each iteration $k$, one or more Monte-Carlo samples are collected and added to samples from preceding iterations, in order to form $C_k^{(\l)}$ and $d_k^{(\l)}$, and that the sample collection method is such that 
$$\lim_{k\to\infty}C_k^{(\l)}=C^{(\l)},\qquad \lim_{k\to\infty}d_k^{(\l)}=d^{(\l)},$$
with probability 1. The iteration \lspeapproxsim\ is known as LSPE($\l$) in approximate DP. Thus, based on the analysis of this paper, LSPE($\l$) can be viewed as an extrapolated form of the projected proximal algorithm of Fig.\ \figonetwo, implemented by simulation.

A popular alternative to LSPE($\l$) is the LSTD($\l$) algorithm, which approximates the solution of the projected equation $C^{(\l)}r=d^{(\l)}$ with the solution of the equation
$$C_k^{(\l)}r=d_k^{(\l)}\xdef\lstdapprox{\lab}\eqnum\show{oneo}$$ 
that iteration \lspeapproxsim\ aims to solve. In  LSTD($\l$) this is done by simple matrix inversion, but the main computational burden  is the calculation of $C_k^{(\l)}$ and $d_k^{(\l)}$. Analysis and simulation suggests that overall the  LSPE($\l$) and LSTD($\l$) algorithms are computationally competitive with each other; see the discussions in [BBN04], [YuB06], [Ber11b], [Ber12a]. Another prominent algorithm in approximate DP is TD($\l$) [cf.\ Eq.\ \tdlambda], which may be viewed as a stochastic approximation method for solving $C^{(\l)}r=d^{(\l)}$. This algorithm has a long history in approximate DP, as noted earlier, and has been extended to the general linear system context in [BeY09], Section 5.3.

A major problem for the preceding simulation-based algorithms is that when $C^{(\l)}$ is nearly singular, $C_k^{(\l)}$ should be a very accurate estimate of $C^{(\l)}$, so that $\s\big(I-C_k^{(\l)}\big)<1$, which is a requirement for the methods \lstdapprox\ and \lspeapprox\ to make sense. The papers [WaB13] and [WaB14] address this issue and provide stabilization schemes to improve the performance of the LSPE($\l$) and  LSTD($\l$) methods, as well as other iterative algorithms for solving singular and near singular systems of equations by simulation. On the other hand the proximal implementations \proxpvi\ and \proxpviextra\ are more tolerant of near-singularity of $C^{(\l)}$ and simulation noise. This is a generic property of the proximal algorithm and the main motivation for its use. The simulation-based version of the  algorithm \proxpviextra, is
$$r_{k+1}=r_k-{\hat c+1\over \hat c}\lf({1\over \hat c} I+C_k^{(\l)}\ri)^{-1}\big(C_k^{(\l)}r_k-d_k^{(\l)}\big).$$
Its use of the extrapolation factor ${\hat c+1\over \hat c}$ may provide significant acceleration, particularly for small values of ${\hat c}$, while simultaneously guarding against near singularity of $C^{(\l)}$. 

Let us also mention another  approach of the proximal projected type, which can also be implemented by simulation. This is to convert the projected equation $C^{(\l)}r=d^{(\l)}$ to the equivalent equation
$${C^{(\l)}}'\Sigma^{-1}C^{(\l)}r={C^{(\l)}}'\Sigma^{-1}d^{(\l)},$$
where $\Sigma$ is a symmetric positive definite matrix, and then apply the proximal algorithm to its solution. The analog of the simulation-based proximal algorithm \proxpvi\ is
$$r_{k+1}=r_k-\lf({1\over \hat c} I+{C_k^{(\l)}}'\Sigma^{-1}C_k^{(\l)}\ri)^{-1}{C_k^{(\l)}}'\Sigma^{-1}\big(C_k^{(\l)}r_k-d_k^{(\l)}\big),\xdef\simprojprox{\lab}\eqnum\show{oneo}$$
and its extrapolated version [cf.\ Eq.\ \proxpviextra] is 
$$r_{k+1}=r_k-{\hat c+1\over \hat c}\lf({1\over \hat c} I+{C_k^{(\l)}}'\Sigma^{-1}C_k^{(\l)}\ri)^{-1}{C_k^{(\l)}}'\Sigma^{-1}\big(C_k^{(\l)}r_k-d_k^{(\l)}\big).$$
These algorithms are valid assuming that $C^{(\l)}$ is invertible,  $\s\big(I-{C^{(\l)}}'\Sigma^{-1}C^{(\l)}\big)\le 1$,  and $\lim_{k\to
\infty}C_k^{(\l)}=C^{(\l)}$. The algorithm \simprojprox\ has been considered both in the approximate DP and the more general linear system context in [WaB13], [WaB14]; see also the textbook presentation of [Ber12a], Sections 7.3.8 and 7.3.9. However, remarkably,  if $C^{(\l)}$ is singular, its iterate sequence $\{r_k\}$ may diverge as shown by Example 9 of the paper [WaB14] (although the residual sequence $\big\{C_k^{(\l)}r_k-d_k^{(\l)}\big\}$ can be shown  to converge to 0 generically; cf.\ Prop.\ 9 of [WaB14]).

\subsection{The Use of Temporal Differences}

\pn The preceding discussion has outlined the general ideas of simulation-based methods, but did not address specifics. Some of the most popular implementations are based on the notion of temporal differences, which are residual-like terms of the form
$$d(x,\ell)=A^\ell(Ax+b-x),\qquad x\in\rn,\ \ell=0,1,\ldots.$$
In particular, it can be shown using the definition \tddef\ that
$$\td x=x+\sum_{\ell=0}^\infty \l^\ell d(x,\ell).\xdef\tdexpression{\lab}\eqnum\show{oneo}$$
This can be verified with the following calculation 
$$\eqalign{T^{(\l)}x&=\sum_{\ell=0}^\infty(1-\l)\l^\ell(A^{\ell+1}x+A^\ell b+A^{\ell-1}b+\cdots+b)\cr
&=x+(1-\l)\sum_{\ell=0}^\infty \l^\ell\sum_{m=0}^\ell(A^mb+A^{m+1}x-A^mx)\cr
&=x+(1-\l)\sum_{m=0}^\infty \lf(\sum_{\ell=m}^\infty \l^\ell\ri)(A^mb+A^{m+1}x-A^mx)\cr
&=x+\sum_{m=0}^\infty \l^m(A^mb+A^{m+1}x-A^mx)\cr}$$
from [BeY09], Section 5.2. 

Based on the least squares implementation \semileastsq\ and the temporal differences expression \tdexpression, and assuming that $\Phi$ has full rank $s$, the projected algorithm
$$\Phi r_{k+1}=\Pi \td \Phi r_k,$$
[cf.\ the LSPE($\l$) algorithm] is given by
$$r_{k+1}=\arg\min_{r\in\re^s}\sum_{i=1}^n\xi_i\lf(\phi_i'r-\phi_i'r_k-\sum_{\ell=0}^\infty \l^\ell d_i(\Phi r_k,\ell)\ri),$$
where $(\xi_1,\ldots,\xi_n)$ is a probability distribution over the indices $1,\ldots,n$, $\phi_1',\ldots,\phi_n'$ are the rows of $\Phi$, and $d_i(\Phi r_k,\ell)$ is the $i$th component of the $n$-dimensional vector $d(\Phi r_k,\ell)$. Equivalently, this algorithm is written as
$$r_{k+1}=r_k+\lf(\sum_{i=1}^n\xi_i \phi_i\phi_i'\ri)^{-1}\sum_{i=1}^n\xi_i \phi_i\lf(\sum_{\ell=0}^\infty \l^\ell d_i(\Phi r_k,\ell)\ri).\xdef\semileastsqt{\lab}\eqnum\show{oneo}$$
In the simulation-based implementation of this iteration, the terms
$$\sum_{i=1}^n\xi_i \phi_i\phi_i'$$
and 
$$\sum_{i=1}^n\xi_i \phi_i\lf(\sum_{\ell=0}^\infty \l^\ell d_i(\Phi r_k,\ell)\ri)$$
are viewed as expected values with respect to the probability distribution $(\xi_1,\ldots,\xi_n)$, and are approximated by sample averages.

The samples may be collected in a variety of ways. Typically, they are obtained by simulating a suitable $n$-state Markov chain to produce one infinite sequence of indexes $(i_0,i_1,\ldots)$ or multiple finite sequences  of indexes 
$(i_0,i_1,\ldots,i_m)$, where $m$ is an integer ($m$ may be different for different sequences). Corresponding samples of the temporal differences are also collected during this process. The transition probabilities of the Markov chain are related to the elements of the matrix $A$, and are chosen in a way to preserve convergence of the iteration \semileastsqt. The details of this, as well as the convergence analysis are complicated, and further review is beyond the scope of this paper. While the formalism of temporal differences is commonly used in practice, simulation-based algorithms may be implemented in other ways (see, e.g., [BeY09], Section 5.1, [Ber12b], [YuB12]). Moreover, similar ideas can be used in simulation-based versions of other related multistep algorithms, such as the one of Eq.\ \altapprox\ as well as analogs of the LSTD($\l$) and TD($\l$) methods. We refer to the literature cited earlier for more detailed discussions.

\vskip  -4mm
\section{Extensions to Nonlinear Fixed Point Problems}
\vskip  -2mm

\pn In this section we consider the solution of the fixed point problem
$$x=T(x),\xdef\nonlinsys{\lab}\eqnum\show{oneo}$$
where $T:\rn\mapsto\rn$ is a possibly nonlinear mapping.\footnote{\dag}
{\ninepoint For mappings $T$ that may be nonlinear, we use the notation $T(x)$ rather than $Tx$, which is used only when $T$ is known to be linear.} The proximal algorithm for this problem has the form
$$x_{k+1}=P^{(c)}(x_k),\xdef\nonlinproxalg{\lab}\eqnum\show{oneo}$$
where $c$ is a positive scalar, and for a given $x\in\rn$, $P^{(c)}(x)$ solves the following equation in the vector $y$:
$$y=T(y)+{1\over c}(x-y).\xdef\proximalit{\lab}\eqnum\show{oneo}$$
We will operate under assumptions guaranteeing that this equation has a unique solution, so that $P^{(c)}$ will be well defined as a point-to-point mapping.

\xdef \figtdproximal_nonlin{\figr}\figrnum\show{myfigure}

We focus on the following extrapolated version of the proximal algorithm \nonlinproxalg:
$$x_{k+1}=E^{(c)}(x_k),\xdef\extrapalg{\lab}\eqnum\show{oneo}$$
where
$$E^{(c)}(x)=x+{c+1\over c}\bl(P^{(c)}(x)-x\br).\xdef\extrap{\lab}\eqnum\show{oneo}$$
When $T$ is linear as in Section 1, this algorithm coincides with the multistep method $x_{k+1}=T^{(\l)}x_k$ (cf.\ Fig.\ \figoneone). 

The key fact for our purposes is that
$$E^{(c)}(x)=T\bl(P^{(c)}(x)\br),\qquad \forall\ x\in\rn;\xdef\nonlintdpb{\lab}\eqnum\show{oneo}$$
see Fig.\ \figtdproximal_nonlin.
Indeed from Eq.\ \proximalit\ we have
$$P^{(c)}(x)+{1\over c}\bl(P^{(c)}(x)-x\br)=T\bl(P^{(c)}(x)\br).\old{\xdef\criticaleq{\lab}\eqnum\show{oneo}}$$
Using the form \extrap\ of $E^{(c)}(x)$ and the preceding equation, we obtain
$$E^{(c)}(x)=x+{c+1\over c}\bl(P^{(c)}(x)-x\br)=P^{(c)}(x)+{1\over c}\bl(P^{(c)}(x)-x\br)=T\bl(P^{(c)}(x)\br),$$
showing Eq.\ \nonlintdpb.

\topinsert
\centerline{\hskip2pc\includegraphics[width=4in]{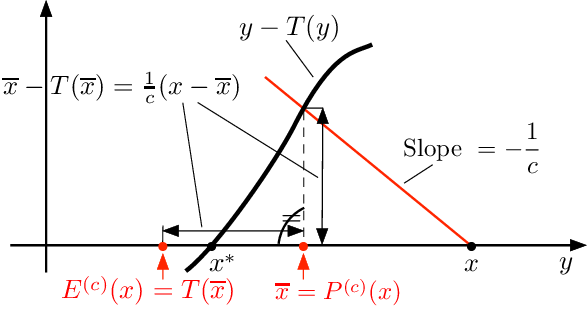}}
\vskip-1pc
\hskip-4pc\fig{0pc}{\figtdproximal_nonlin.} {Illustration of the extrapolated algorithm \extrapalg. The proximal iterate $P^{(c)}(x)$, denoted $\ol x$ in the figure, is extrapolated by ${1\over c}(\ol x-x)$. From the definition of $P^{(c)}(x)$, the extrapolated iterate is equal to $T(\ol x)$ [cf.\ Eq.\ \nonlintdpb], and its distance to $x^*$ is strictly smaller than the distance of $\ol x$ when $T$ is a contraction.}\endinsert

The form of Eq.\ \nonlintdpb\ suggests that the extrapolated iteration \extrapalg\ has faster convergence than the proximal iteration \nonlinproxalg, within contexts where $T$ is contractive with respect to a suitable norm. In particular, if the solution $P^{(c)}(x)$ of Eq.\ \proximalit\ exists and is unique for all $x\in\rn$, and $P^{(c)}$ and $T$ are contractions with respect to the same norm, then both iterations \nonlinproxalg\ and \extrapalg\ converge to the unique fixed point of $T$, and the extrapolated iteration converges faster. The following proposition provides specific conditions guaranteeing that this is so.

\xdef\propfourtwo{\propn}\propnum\show{myproposition}
 
%\texshopbox
{\proposition{\propfourtwo:} Assume that $T$ is a contraction mapping with modulus $\g\in(0,1)$ with respect to a Euclidean norm $\|\cdot\|$, i.e., 
$$\big\|T(x_1)-T(x_2)\big\|\le \g \|x_1-x_2\|,\qquad \forall\ x_1,x_2\in\rn.\xdef\contractprop{\lab}\eqnum\show{oneo}$$
Then the solution $P^{(c)}(x)$ of Eq.\ \proximalit\ exists and is unique for all $x\in\rn$, and the mappings $P^{(c)}$
and $E^{(c)}$ are contraction mappings with respect to $\|\cdot\|$. In particular, we have%}\texshopboxnt{\pn
$$\big\|P^{(c)}(x_1)-P^{(c)}(x_2)\big\|\le {1\over 1+c(1-\g)}\|x_1-x_2\|,\qquad\forall\ x_1,x_2\in\rn,$$
$$\big\|E^{(c)}(x_1)-E^{(c)}(x_2)\big\|\le {\g\over 1+c(1-\g)}\|x_1-x_2\|,\qquad\forall\ x_1,x_2\in\rn.$$%}\texshopboxnt{\pn
Moreover, every sequence
$\{x_k\}$ generated by either the proximal algorithm \nonlinproxalg\ or its extrapolated version \extrapalg\ converges geometrically to the unique fixed point $x^*$ of $T$, and the convergence of the extrapolated version is faster in the sense that
$$\bl\|E^{(c)}(x)-x^*\br\|\le \g\bl\|P^{(c)}(x)-x^*\br\|,\qquad \forall\ x\in\rn.\xdef\extraprate{\lab}\eqnum\show{oneo}$$
}

\proof We first verify that the mapping $x\mapsto x-T(x)$ satisfies the standard strong monotonicity assumption under which the proximal mapping is a contraction. In particular, denoting by $\langle\cdot,\cdot\rangle$ the inner product that defines the Euclidean norm $\|\cdot\|$, and using the Cauchy-Schwarz inequality and Eq.\ \contractprop, we have
$$\eqalign{\langle x_1-x_2,x_1-T(x_1)-x_2+T(x_2)\rangle&= \|x_1-x_2\|^2-\langle x_1-x_2,T(x_1)-T(x_2)\rangle\cr
&\ge \|x_1-x_2\|^2-\|x_1-x_2\|\cdot\big\|T(x_1)-T(x_2)\big\|\cr
&\ge \|x_1-x_2\|^2-\g\|x_1-x_2\|^2\cr 
&=(1-\g) \|x_1-x_2\|^2,\cr}\qquad \forall\ x_1,x_2\in\rn.$$
This relation shows that the mapping $x\mapsto x-T(x)$ is strongly monotone and from standard results, $P^{(c)}$ is well-defined as a point-to-point mapping and we have
$$\big\|P^{(c)}(x_1)-P^{(c)}(x_2)\big\|\le {1\over 1+c(1-\g)}\|x_1-x_2\|,\qquad\forall\ x_1,x_2\in\rn,$$
(see [Roc76] or [Ber15], Exercise 5.2). In view of Eq.\ \nonlintdpb\ and the contraction property of $T$, the corresponding contraction property of $E^{(c)}$ and Eq.\ \extraprate\ follow.
\qed

\subsubsection{Extrapolation of the Forward-Backward and Proximal Gradient Algorithms}

 \xdef\figforwback{\figr}\figrnum\show{myfigure}
\xdef\figforwbackext{\figr}\figrnum\show{myfigure}

\pn The forward-backward splitting algorithm applies to the fixed point problem $x=T(x)-H(x)$, where $T$ is a maximally monotone point-to-set mapping in a Euclidean space with inner product $\langle\cdot,\cdot\rangle$ defining a Euclidean norm $\|\cdot\|$, and $H$ is single-valued and strongly monotone, in the sense that for some scalar $\b>0$, we have
$$\langle x_1-x_2,H(x_1)-H(x_2)\rangle\ge  \b\|x_1-x_2\|^2,\qquad \forall\ x_1,x_2\in\rn.$$
The algorithm has the form
$$x_{k+1}=P^{(\a)}\big(x_k-\a H(x_k)\big),\qquad \a>0,$$
where $P^{(\a)}$ is the proximal mapping corresponding to $T$; see Fig.\ \figforwback. This algorithm was analyzed at various levels of generality, by Lions and Mercier [LiM79], Gabay [Gab83], and Tseng [Tse91]. It has been shown to converge to $x^*$ if $\a$ is sufficiently small.
For a minimization problem where $H$ is the gradient of a strongly convex function, it becomes the popular proximal gradient algorithm; for recent surveys, see Beck and Teboulle [BeT10], Parikh and Boyd [PaB13], and the author's textbook [Ber15] (Ch.\ 6), among others.

\midinsert
\centerline{\hskip2pc\includegraphics[width=4in]{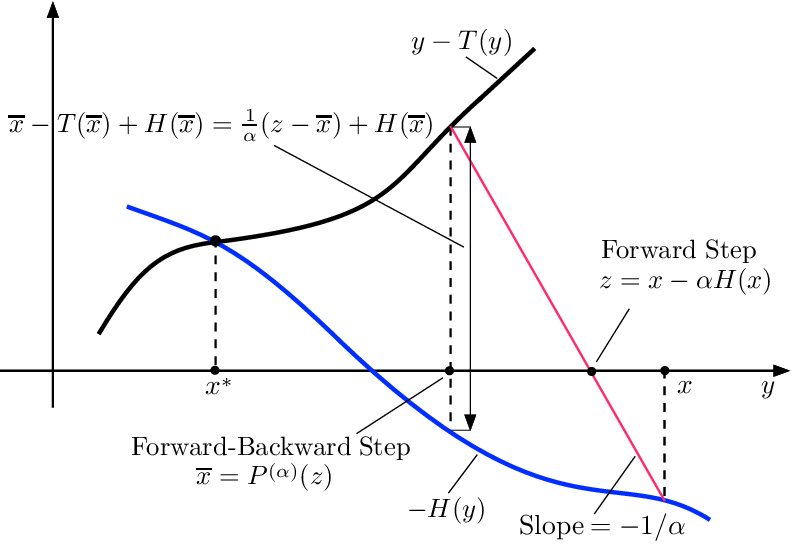}}
\vskip-1pc
\hskip-4pc\fig{0pc}{\figforwback.} {Illustration of an iteration of the forward backward algorithm.}\endinsert
 
The extrapolated forward-backward algorithm has the form
$$z_k=x_k-\a H(x_k),\qquad \ol x_{k}=P^{(\a)}(z_k),$$
$$x_{k+1}=\ol x_{k}+{1\over \a}(\ol x_{k}-z_k)-H(\ol x_{k}),$$
and is illustrated  in Fig.\ \figforwbackext. It can be seen that 
$$x_{k+1}=T(\ol x_{k})-H(\ol x_{k})$$
so there is acceleration if the mapping $T-H$ is contractive. In the linear case where $T(x)=Ax+b$, $H(x)=Bx$, where $A$ and $B$ are $n\times n$ matrices, the algorithm can be related to temporal difference methods, and may be implemented using simulation-based techniques.

\vskip-1pc

\section{Linearized Proximal and Temporal Difference Methods for\hfill\break Nonlinear Problems}

\pn The proximal algorithm \nonlinproxalg\ and its extrapolated version \extrap\ cannot be related to multistep temporal difference algorithms when $T$ is nonlinear, because then the mapping $P^{(c)}$ does not admit a power series expansion; cf.\ Eq.\ \pbpower. It is possible, however, to consider algorithmic ideas based on linearization whereby $T$ is linearized at each iterate $x_k$, and  the next iterate $x_{k+1}$ is obtained with a temporal differences-based (exact, approximate, or extrapolated) proximal iteration using the linearized mapping. This type of algorithm bears similarity to Newton's method for solving nonlinear fixed point problems, the difference being that the linearized system is solved approximately, using a {\it single proximal iteration\/}, rather than exactly (as in Newton's method). The algorithm does not seem to have been considered earlier, to the author's knowledge, although related ideas underlie the $\l$-policy iteration and optimistic policy iteration methods in DP (see the end of Section 2 and the subsequent discussion).

\topinsert
\centerline{\hskip2pc\includegraphics[width=4.5in]{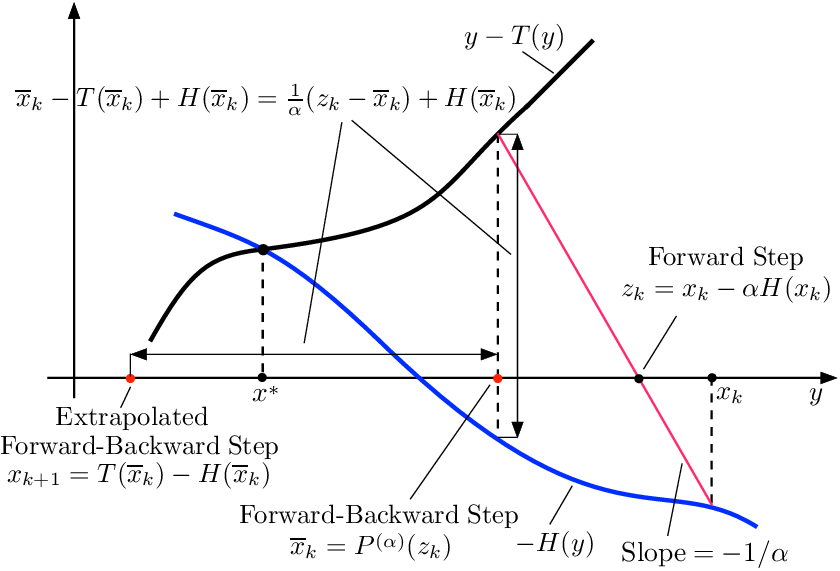}}
\vskip-1pc
\hskip-4pc\fig{0pc}{\figforwbackext.} {Illustration of the forward backward algorithm with extrapolation.}\endinsert

We focus on the fixed point problem $x=T(x)$ with the $i$th component of $T(x)$ having the form 
$$T(i,x)=\min_{\m(i)\in  M(i)}\big\{a\big(i,\m(i)\big)'x+b\big(i,\m(i)\big)\big\},\qquad x\in \rn,\ i=1,\ldots,n,\xdef\minform{\lab}\eqnum\show{oneo}$$
where for each $i$, $M(i)$ is some set, and  $a\big(i,\m(i)\big)$ and $b\big(i,\m(i)\big)$ are a (column) vector in $\rn$ and scalar, respectively, for each $\m(i)\in  M (i)$. We assume that the minimum in Eq.\ \minform\ is attained throughout this section. For a given $i$, this form of  $T(i,x)$ includes a very broad class of concave functions of $x$. Intuition suggests that an algorithmic analysis for more general forms of $T$ may be possible, but this is beyond the scope of the present paper (see the discussion of Section 5.2).

Let ${\cal M}$ be the Cartesian product $ M (1)\times\cdots\times  M (n)$. Given a vector $\m=\big((\m(1),\ldots,\m(n)\big)\in  {\cal M} $, the matrix whose $i$th row is the vector $a\big(i,\m(i)\big)'$ is denoted by $A_\m$ and the vector whose $i$th component is $b\big(i,\m(i)\big)$ is denoted by $b_\m$. We denote by $T_\m$ the linear  mapping given by
$$T_\m x=A_\m x+b_\m,$$
and we write in shorthand
$$T(x)=\min_{\m\in {\cal M}}T_\m x,$$
where the minimum over $\m$ above is meant  separately, for each of the $n$ components of $T_\m x$. 
Our notation here and later is inspired by notation widely used in DP and policy iteration contexts, where $i$ corresponds to state, the components $x(i)$ of $x$ correspond to cost at state $i$, $\m(i)$ corresponds to control at state $i$, $ M (i)$ corresponds to the control constraint set at state $i$, $\m$ corresponds to policy, $A_\m$ and $b_\m$ correspond to the transition probability matrix and cost per stage vector for policy $\m$, $T_\m$ is the mapping that defines Bellman's equation for the policy $\m$, and $T$ is the mapping that defines Bellman's equation for the corresponding Markovian decision problem. 

\xdef\figfproxlin{\figr}\figrnum\show{myfigure}

\xdef\figvalueit{\figr}\figrnum\show{myfigure}

\subsection{A Linearized Proximal Algorithm Under a Monotonicity Assumption}

\pn In this subsection, we will introduce a linearized algorithm, which we will analyze under some assumptions. Chief among these assumptions is that for all $\m\in {\cal M}$, the mapping $T_\m$ is monotone, in the sense that the matrix $A_\m$ has nonnegative components, and that the initial condition $x_0$ must satisfy $x_0\ge T(x_0)$. In the next subsection, we will consider a randomized variant of this algorithm under an alternative set of assumptions.

At the typical iteration, given the current iterate $x_k$, we find $\m_{k}(i)$ that attains the minimum over $\m(i)\in  M (i)$ of  $a\big(i,\m(i)\big)'x_k+b\big(i,\m(i)\big)$, $i=1,\ldots,n$, and we let $\m_k=\big(\m_k(1),\ldots,\m_k(n)\big)$. We denote this by writing
$$T_{\m_k}x_k=\min_{\m\in {\cal M}}T_\m x_k=T(x_k),$$
or
$$\m_k\in\arg\min_{\m\in {\cal M}}T_\m x_k,\xdef\tdpolimproveo{\lab}\eqnum\show{oneo}$$
where the minimum is meant separately, for each of the $n$ components of $T_\m x_k$. 
 We obtain $x_{k+1}$ via the multistep (extrapolated proximal) iteration
$$x_{k+1}=T_{\m_k}^{(\l)}x_k,\xdef\tdproxlin{\lab}\eqnum\show{oneo}$$
where for a given $\l\in(0,1)$, $T_{\m_k}^{(\l)}$ is the multistep mapping corresponding to the linear mapping $T_{\m_k}$, 
$$T_{\m_k} x=A_{\m_k} x+b_{\m_k},\qquad x\in\rn;\xdef\linsys{\lab}\eqnum\show{oneo}$$
cf.\ Eq.\ \tddef. The algorithm is illustrated in Fig.\ \figfproxlin, together with its proximal version. Note that $a\big(i,\m_k(i)\big)$ is the gradient of $T(i,\cdot)$ at $x_k$ if $T(i,\cdot)$ is differentiable, and otherwise it is a subgradient  of $T(i,\cdot)$ at $x_k$. This justifies the terms ``linearization" and ``linearized mapping."

\topinsert
\centerline{\hskip2pc\includegraphics[width=4.9in]{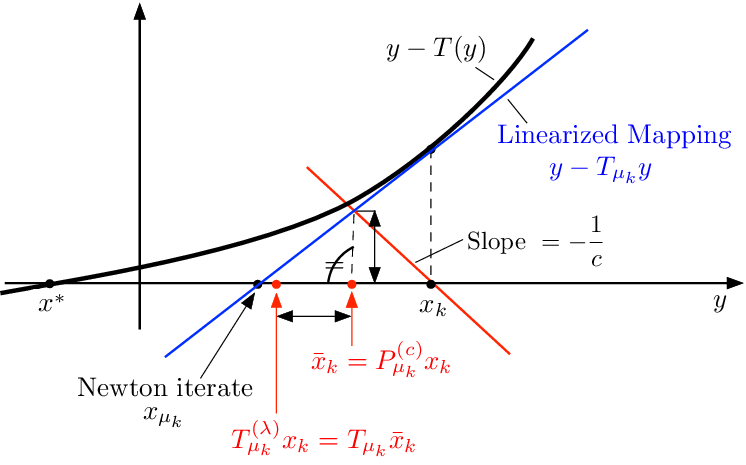}}
\vskip-1pc
\hskip-4pc\fig{0pc}{\figfproxlin.} {Illustration of the linearized multistep algorithm \tdpolimproveo-\tdproxlin, and its proximal version. At the current iterate $x_k$, we linearize $T$ and find the proximal iterate $\bar x_k=P_{\m_k}^{(c)}x_k$ that aims to find the fixed point $x_{\m_k}$ of the linearized mapping $T_{\m_k}$. We can then find the multistep iterate by extrapolation 
$$T_{\m_k}^{(\l)}x_k=T_{\m_k}\bar x_k=\bar x_k+{1\over c}(\bar x_k-x_k);$$
[cf.\ Eq.\ \tdproxlin]. Alternatively, $T_{\m_k}^{(\l)}x_k$ can be found by a temporal differences-based calculation. Note the similarity with the form of Newton's method that finds $x_{\m_k}$, the unique fixed point of $T_{\m_k}$, i.e., the iteration $x_{k+1}=x_{\m_k}$. Newton's method is generally faster but may require much more overhead than the linearized proximal or multistep iteration.}\endinsert

The algorithm \tdproxlin-\linsys\ is related to the $\l$-policy iteration method for DP problems where $\{\m_k\}$ is the sequence of generated policies, and the fixed point equation $x=T_{\m_k}x$ corresponds to Bellman's equation for the policy $\m_k$ (see the discussion and references given at the end of Section 2). The algorithm admits several variations where $T_{\m_k}^{(\l)}$ is replaced in Eq.\ \tdproxlin\ by some other related mapping; for example the iteration 
$$x_{k+1}=P_{\m_k}^{(c)}x_k,\xdef\proxlin{\lab}\eqnum\show{oneo}$$
where $P_{\m_k}^{(c)}$ is the proximal mapping corresponding to $T_{\m_k}$, or the iteration  \altapprox. Another related possibility is the iteration 
$$x_{k+1}=T_{\m_k}^m x_k,\xdef\optpolit{\lab}\eqnum\show{oneo}$$
 where $T_{\m_k}^m$ is the composition of $T_{\m_k}$ with itself $m$ times ($m\ge1$). This is related to the optimistic policy iteration method of DP; see [BeT96], [Ber12a], or [Ber18], Section 2.5. 
 
 Figure \figvalueit\ illustrates the fixed point iterate $T_{\m_k}x_k$ [Eq.\ \optpolit\ with $m=1$], which is equal to $T(x_k)$. A comparison with Fig.\ \figfproxlin\ shows that the multistep iterate $T_{\m_k}^{(\l)}x_k$ is closer to the Newton iterate $x_{\m_k}$ than the  fixed point iterate $T_{\m_k}x_k$ for large values of $\l$, and in fact converges to $x_{\m_k}$ as $\l\to1$. The multistep iterate $T_{\m_k}^{(\l)}x_k$ approaches the fixed point iterate $T(x_k)=T_{\m_k}x_k$ as $\l\to0$. The proximal iterate $P_{\m_k}^{(c)}x_k$ [cf.\ Eq.\ \proxlin] also approaches $x_{\m_k}$ as $c\to\infty$ (i.e., $\l\to1$), but approaches $x_k$ as $c\to0$ (i.e., $\l\to0$).
 
 \topinsert
\centerline{\hskip2pc\includegraphics[width=4.9in]{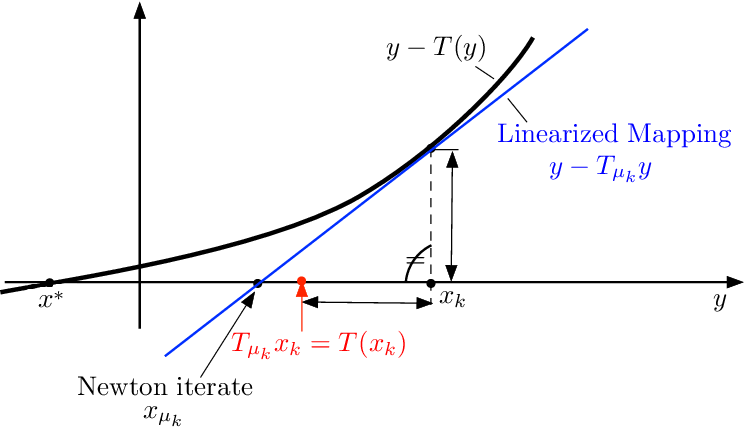}}
\vskip-1pc
\hskip-4pc\fig{0pc}{\figvalueit.} {Illustration of the iterate $T_{\m_k}x_k$ [Eq.\ \optpolit\ with $m=1$]. It is equal to the fixed point iterate $T(x_k)$. }\endinsert

We now introduce a framework for the convergence analysis of the algorithm \tdpolimproveo-\tdproxlin. We introduce a set $X\subset \rn$, within which we will require that the iterates $x_k$  lie. It is possible that $X=\rn$ but some interesting special cases are obtained if $X$ is a strict subset of $\rn$, such as for example the nonnegative orthant. Let us say that a vector $\m\in {\cal M} $ is {\it proper} if  $A_\m$ has eigenvalues strictly within the unit circle and the unique fixed point of $T_\m$, denoted by $x_\m$, lies within $X$. If $\m$ is not proper it is called {\it improper\/}. The names ``proper" and ``improper" relate to notions of proper and improper policies in DP, and stochastic shortest path problems in particular; see [BeT91], [BeT96], [Ber12a], [Ber18]. Note that for proper $\m$ the algorithmic results of Section 2 come into play, in view of the linearity of $T_\m$. In particular, the multistep and proximal iterations $x_{k+1}=T_\m^{(\l)} x_k$ and $x_{k+1}=P_\m^{(c)} x_k$  converge to $x_\m$ starting from any $x_0\in\rn$.

We will assume the following.

\xdef\assumptionlin{\assumptionn}\assumptionnum\show{myproposition}

%\texshopbox
{\assumption{\assumptionlin:}\nitem{(a)} For  all $x\in \rn$ and $i=1,\ldots,n$, the minimum over $ M (i)$ in Eq.\ \minform\ is attained.%}\texshopboxnt{
\nitem{(b)} For all $\m\in  {\cal M} $, the matrix $A_\m$ has nonnegative components.%}\texshopboxnt{
\nitem{(c)} There exists at least one proper vector, and for each improper vector $\m$ and $x\in X$, at least one component of the sequence $\{T_\m^k x\}$ diverges to $+\infty$.
}
\smskip

Assumption \assumptionlin(a) is needed to ensure that the algorithm is well defined. Assumption \assumptionlin(b) implies a monotonicity property typically encountered in  DP problems, whereby we have for all $\m\in {\cal M} $, 
$$T_\m x\le T_\m y\qquad \forall\ x,y\in\rn\hbox{ such that }x\le y,\xdef\montone{\lab}\eqnum\show{oneo}$$ 
as well as
$$T (x)\le T (y)\qquad \forall\ x,y\in\rn\hbox{ such that }x\le y.\xdef\montonet{\lab}\eqnum\show{oneo}$$
[The relation  \montonet\ follows from the relation \montone\ by first taking the infimum of the left side to obtain 
$T(x)\le T_\m y$ and then by taking the infimum of the right side over $\m\in {\cal M} $.]  The monotonicity property can be replaced by a sup-norm contraction assumption on $A_\m$, but this requires substantial algorithmic modifications that will be addressed in part in the next subsection. Note that parts (b) and (c) of Assumption \assumptionlin\ are satisfied if $X=\rn$, and $A_\m$ has nonnegative components and eigenvalues strictly within the unit circle for all $\m\in {\cal M}$, in which case all $\m\in {\cal M}$ are proper.
 
We have the following proposition.

\xdef\propproper{\propn}\propnum\show{myproposition}
 
%\texshopbox
{\proposition{\propproper:} Let Assumption \assumptionlin\ hold. 
For a given $\m\in {\cal M} $, if for some $x\in X$ we have $T_\m x\le x$, then $\m$ is proper.}

\proof By the monotonicity of $T_\m$ [cf.\ Eq.\ \montone], we have $T_\m^k x\le x$ for all $k$, so if $\m$ were improper, 
Assumption \assumptionlin(c) would be violated. \qed

We now consider a restricted optimization problem over the proper policies, which is inspired from the theory of semicontractive problems in abstract DP (see [Ber18]). We introduce the componentwise minimum vector $x^*$, which has components $x^*(i)$ given by
$$x^*(i)=\inf_{\m:\hbox{\eightpoint proper}}x_\m(i),\qquad i=1,\ldots,n,\xdef\hatxdef{\lab}\eqnum\show{oneo}$$
where $x_\m(i)$ is the $i$th component of the vector $x_\m$. The following proposition gives the central results of this subsection.

\xdef\propproxlin{\propn}\propnum\show{myproposition}
 
%\texshopbox
{\proposition{\propproxlin:} Let Assumption \assumptionlin\ hold.
\nitem{(a)} If $T$ has a fixed point within $X$, then this fixed point is equal to $x^*$. Moreover, there exists a proper $\m$ that attains the infimum in Eq.\ \hatxdef.
\nitem{(b)} If $x^*\in X$, then $x^*$ is the unique fixed point of $T$ within $X$. 
\nitem{(c)} A sequence $\{x_k\}$ generated by the algorithm \tdpolimproveo-\tdproxlin\ starting from an initial condition $x_0$ such that $x_0\ge T(x_0)$ is monotonically nonincreasing and converges to $x^*$. Moreover, in the corresponding sequence $\{\m_k\}$ all  $\m_k$ are proper.
}

\proof (a) Let $\hat x\in X$ be a fixed point of $T$ within $X$. We will show that $\hat x=x^*$. Indeed, using also the monotonicity of $T_\m$ and $T$ [cf.\ Eqs.\ \montone\ and \montonet], we have for every $m\ge 1$ and proper $\m$ 
$$\hat x=T(\hat x)\le T_\m\hat x\le T_\m^m \hat x\le \lim_{m\to\infty}T_\m^m \hat x=x_\m,$$
where the second and third inequalities follow from the first inequality and the monotonicity of $T_\m$. By taking the infimum of the right side over all proper $\m$, we obtain $\hat x\le x^*$. For the reverse inequality, let $\hat \m$ be such that $\hat x=T(\hat x)=T_{\hat \m}\hat x$ 
[there exists such $\hat\m$ by Assumption \assumptionlin(a)]. Using Prop.\ \propproper, it follows that $\hat\m$ is proper, so that $\hat x$ is the unique fixed point of $T_{\hat\m}$, i.e., $\hat x=x_{\hat \m}\ge x^*$. Thus $\hat x=x_{\hat \m}=x^*$ and the proof is complete.

\smskip

\pn (b) For every proper $\m$ we have $x_\m\ge x^*$, so by the monotonicity of $T_\m$,
$$x_\m=T_\m x_\m\ge T_\m x^*\ge T x^*.$$
Taking the infimum over all proper  $\m$, we obtain
$x^*\ge T x^*.$
Let $\m$ be such that 
$T x^*=T_\m x^*$. The preceding relations yield
$x^*\ge T_\m x^*$, so by Prop.\ \propproper, $\m$ is proper. Therefore, we have
$$x^*\ge T (x^*)= T_\m x^*\ge \lim_{k\to\infty}T_\m^k x^*=x_\m\ge x^*,$$
where the second equality holds since $\m$ is proper, and $x^*\in X$ by assumption. Hence equality holds throughout in the above relation, which proves that $x^*$ is a fixed point of $T$, which is unique by part (a).

\smskip

\pn (c) For all $x$ and $\m\in {\cal M} $ such that $x\ge T_\m x=T(x)$, we claim that
$$x\ge T_\m x\ge T_\m^{(\l)} x\ge T_\m\cdot T_\m^{(\l)} x\ge T(T_\m^{(\l)} x),\xdef\lambdaineq{\lab}\eqnum\show{oneo}$$
To see this, note that the second inequality follows from the power series expansion
$$T_\m^{(\l)} x=(1-\l)(T_\m x+\l T_\m^2x+\l^2 T_\m^3x+\cdots),\xdef\powerseries{\lab}\eqnum\show{oneo}$$
and the fact that $x\ge T_\m x$ implies $T_\m^m x\ge T_\m^{m+1} x$ for all $m\ge 0$. The latter fact also implies  the third inequality in Eq.\ \lambdaineq, by applying $T_\m$ to both sides of Eq.\ \powerseries.
We  also  have 
$$T_\m^{(\l)} x\ge x_\m,\xdef\lambdaineqt{\lab}\eqnum\show{oneo}$$
 which follows by taking the limit as $m\to\infty$ in the relation
$$\sum_{\t=0}^{m} \l^tT_{\m}^{t+1}x\ge {1-\l^{m+1}\over 1-\l}T_\m^{m+1} x.$$

Thus from the relations \lambdaineq\ and \lambdaineqt\ we have that if $x\in X$ and $x\ge T(x)=T_\m x$, then $\m$ is proper and the vector $\bar x=T_\m^{(\l)} x$ satisfies $\bar x\in X$, $x\ge \bar x\ge T(\bar x)$, and $x\ge x_\m\ge x^*$.
Applying repeatedly this argument and using the definition of the algorithm we have
$$x_k\ge T(x_k)= T_{\m_k}x_k\ge T_{\m_k}^{(\l)}x_k=x_{k+1},\qquad x_k\ge  x_{\m_k}\ge x^*,\qquad k=0,1,\ldots.\xdef\ineqstring{\lab}\eqnum\show{oneo}$$
It follows that the sequence $\{\m_k\}$ consists of proper vectors and $x_k\downarrow x_{\infty}$, where $x_{\infty}$ is some vector with $x_{\infty}\ge x^*$. 

Next we show that $x_\infty$ is a fixed point of $T$.\footnote{\dag}
{\ninepoint Note here that $x_\infty$ may have some components that are equal to $-\infty$. Still, however, because the components of $A_\m$ are assumed nonnegative, the vectors $T_\m x_\infty$ and $T(x_\infty)$ are well-defined as $n$-dimensional vectors with components that are either real or are equal to $-\infty$.} Indeed, from the relation $x_k\ge T(x_k)$ [cf.\ Eq.\ \ineqstring], we have $x_k\ge T(x_\infty)$ for all $k$, so that $x_\infty\ge T(x_\infty)$. Also we note that from the relation $T(x_k)\ge x_{k+1}$ we have 
$$\lim_{k\to\infty}T(x_k)\ge x_\infty.$$
Moreover for all $\m\in {\cal M} $, in view of the linearity of $T_\m$, we have
$$T_\m x_\infty=\lim_{k\to\infty}T_\m x_k\ge \lim_{k\to\infty}T(x_k).$$
By combining the preceding two relations, we obtain $T_\m x_\infty\ge x_\infty,$ so by taking the componentwise minimum over  $\m\in {\cal M} $, we have $T(x_\infty) \ge x_\infty$. This relation, combined with the relation $x_\infty\ge T(x_\infty)$ shown earlier, proves that $x_\infty$ is a fixed point of $T$. 

Finally we show that $x_\infty=x^*$. Indeed, since $x_\infty$ is a fixed point of $T$ and $x^*\le x_\infty$, as shown earlier, we have
$$x^*\le x_\infty=T^k (x_\infty)\le T_{\m}^k x_\infty\le T_{\m}^k x_{\m^0},\qquad \forall\ \m:\hbox{proper},\ k=0,1,\ldots.$$
By taking the limit as $k\to\infty$, and using the fact that $x_{\m_0}\in X$, it follows that
$x^*\le x_\infty\le x_\m$ for all $\m$ proper.
By taking the infimum over proper $\m$, it follows that $x_\infty=x^*$, so $x^*$ is a fixed point of $T$.\qed

The initial condition requirement $x_0\ge T(x_0)$ is somewhat restrictive and will be removed in the next subsection after the algorithm \tdpolimproveo-\tdproxlin\ is modified. We can also prove results that are similar to the preceding proposition, but where the nonnegativity Assumption \assumptionlin(b) and the condition $x_0\ge T(x_0)$ are replaced by alternative conditions. For example, linearization algorithms for finding a fixed point of $T$ are given in the author's monograph [Ber18], Section 2.6.3, under just the assumption that all the matrices $A_\m$, $\m\in M$, are contractions with respect to a common weighted sup-norm (see also the papers by Bertsekas and Yu [BeY10], [BeY12], [YuB13], which relate to the discounted DP and stochastic shortest path contexts). These algorithms, however, do not use proximal iterations.  In any case, an initial vector $x_0$ with $x_0\ge T(x_0)$ [cf.\ the assumption of part (b)] may be obtained in some important cases by adding sufficiently large scalars to the components of some given vector $x$. For example, suppose that for some $\m$, the mapping $T_\m$ is a sup-norm contraction. Then given any $x\in\rn$, it can be shown that the vector $x_0$ with components equal to $x(i)+r$ satisfies $x_0\ge T_\m x_0\ge T(x_0)$, provided that the scalar $r$ is sufficiently large. 

The monotone nonincreasing property of the sequence $\{x_k\}$, shown in Prop.\ \propproxlin(c), suggests that the algorithm is convergent even when implemented in distributed asynchronous fashion. Indeed this can be shown with an analysis similar to the one given in [BeY10]. Some well-known counterexamples from DP by Williams and Baird [WiB93] suggest that the assumption $x_0\ge T(x_0)$ is essential for asynchronous convergence (unless the alternative algorithmic framework of [Ber18], Section 2.6.3, is used).  Proposition \propproxlin(c) also shows that even without a guarantee of existence of a fixed point of $T$ within $X$, we have $x_k\downarrow x_\infty$, where $x_\infty$ is some vector that may have some components that are equal to $-\infty$. For an example of this type, consider the one-dimensional problem of finding a fixed point of the mapping
$$T(x)=\min_{\m\in (0,1]}\big\{(1-\m^2)x-\m\big\}.$$
Then Assumption \assumptionlin\ is satisfied with $X=\re$, we have $x_\m=-1/\m$, $x^*=-\infty$, $x_k\downarrow x^*$ starting from any $x_0\in\re$, while $T$ has no real-valued fixed point.

To see what may happen when Assumption \assumptionlin(c) is not satisfied, let $X=\re$ and consider the mapping $T$ given by
$$T(x)=\min\{1,x\},$$
which is of the form \minform\ but has multiple fixed points, thus violating the conclusion of Prop.\ \propproxlin(a). 
Here ${\cal M}$ consists of two vectors: one is $\hat \m$ with $T_{\hat \m}x=1$, which is proper, and the other is $\bar \m$ with $T_{\bar \m}x=x$, which is improper but does not satisfy Assumption \assumptionlin(c).

\subsection{A Linearized Proximal Algorithm Under a Contraction Assumption}

\pn We will now discuss a randomized version of the linearized proximal algorithm \tdpolimproveo-\tdproxlin, which we will analyze under an alternative set of assumptions. While we remove the monotonicity Assumption \assumptionlin(b) and the restriction $x_0\ge T(x_0)$, we introduce a finiteness assumption on ${\cal M}$, and a contraction assumption on $T_\m$ and $T$. In particular, the components of $T$ are required to be concave polyhedral functions of $x$, with ``gradients" whose norm is strictly less than one. We will assume the following.

\xdef\assumptionlino{\assumptionn}\assumptionnum\show{myproposition}

%\texshopbox
{\assumption{\assumptionlino:}Let $\|\cdot\|$ be some norm in $\rn$, and assume the following: \nitem{(a)}  The set ${\cal M}$ is finite.%}\texshopboxnt{\pn
\nitem{(b)} The mappings $T_\m$, $\m\in{\cal M}$, and $T$ are contractions with respect to $\|\cdot\|$, with modulus $\r$, and fixed points $x_\m$ and $x^*$, respectively. }
\smskip

In the special case where  the mappings $T_\m$, $\m\in{\cal M}$, are contractions, with respect to a weighted sup-norm 
$$\|x\|=\max_{i=1,\ldots,n}{|x^i|\over v^i},\qquad x=(x^1,\ldots,x^n)\in\rn,$$
where $(v^1,\ldots,v^n)$ is a vector of positive scalar weights, we have that $T$ is also a contraction with respect to the same weighted sup-norm. This is well known in the theory of abstract DP (see [Den67], [BeS78], [Ber18]). For a verification, we write for all $x=(x^1,\ldots,x^n)$, $\tl  x=(\tl  x^1,\ldots,\tl  x^n)$, and  $\m\in {\cal M}$,
$$(T_\m x)^i\le (T_\m \tl  x)^i+ \r \|x-\tl  x\|\,v^i,\qquad \forall\ i=1,\ldots,n.$$
By taking infimum of both sides over $\m\in {\cal M}$, we have
$${\big(T(x)\big)^i-(T\big(\tl  x)\big)^i}\le \r \|x-\tl  x\|\,v^i,\qquad \forall\ i=1,\ldots,n.$$
Reversing the roles of $x$ and $\tl  x$, we also have
$${\big(T(\tl  x)\big)^i-\big(T(x)\big)^i\over v^i}\le \r \|x-\tl  x\|,\qquad \forall\ i=1,\ldots,n.$$
By combining the preceding two relations, and taking the maximum of the left side over $i=1,\ldots,n$, we obtain $\big\|T(x)-T(\tl  x)\big\|\le \r \|x-\tl  x\|$. For other norms, however, $T$ may not be a contraction even if all the mappings $T_\m$ are, so its contraction property must be verified separately.

 To construct a convergent linearized proximal algorithm for finding the fixed point of $T$, we address a fundamental difficulty of the algorithm \tdpolimproveo-\tdproxlin, which is that the iteration $x_{k+1}=T_{\m_k}^{(\l)} x_k$ [cf.\ Eq.\ \optpolit] approaches the fixed point $x_{\m_k}$ of $T_{\m_k}$,  thus constantly aiming at a ``moving target" that depends on $k$. This can cause an oscillatory behavior, whereby the algorithm can get locked into a repeating cycle. We hence introduce modifications that involve a randomization, which allows the algorithm to recover from such a cycle. 

In particular, we introduce a fixed probability $p\in(0,1)$, and for a given $\l\in(0,1)$, we consider the following algorithm 
$$x_{k+1}=\cases{
T_{\m_k}^{(\l)} x_k& with probability $p$,\cr
T (x_k)& with probability $1-p$,\cr
}\xdef\randomoptpolito{\lab}\eqnum\show{oneo}$$
 where $\m_k$ is detremined according to 
 $$\m_{k}\in\arg\min_{\m\in {\cal M}}T_\m x_k,\xdef\mupdate{\lab}\eqnum\show{oneo}$$
each time the iteration $x_{k+1}=T(x_k)$ is performed. A similarly randomized version of the algorithm \tdpolimproveo, \optpolit\ was introduced in Section 2.6.2 of the abstract DP monograph [Ber18]. It has the form
 $$x_{k+1}=\cases{
T_{\m_k}^{m_k} x_k& with probability $p$,\cr
 T (x_k)& with $1-p$,\cr
}\xdef\randomoptpolitz{\lab}\eqnum\show{oneo}$$
 where $m_k$ is a positive integer from a bounded range, and $\m_k$ is updated by Eq.\ \mupdate\
 each time the iteration $x_{k+1}=T(x_k)$ is performed. A variant of this algorithm is to select to each iteration $k$ a positive integer $m_k$ by randomization using a probability distribution $\big\{p(1),p(2),\ldots\big\}$ over the positive integers, with $p(1)>0$, and update  $\m_k$ by Eq.\ \mupdate\ each time $m_k=1$.
 
In what follows, we denote by ${\cal M}^*$ the set of all $\m\in {\cal M}$ such that $T_{\m}x^*=T(x^*)$:
$${\cal M}^* =\big\{\m\in{\cal M}\mid T_{\m}x^*=T(x^*)\big\}.$$
An important fact (which relies strongly on the finiteness of ${\cal M}$) is given in the following proposition.

\xdef\propcontactiont{\propn}\propnum\show{myproposition}
 
%\texshopbox
{\proposition{\propcontactiont:} Let Assumption \assumptionlino\ hold. Then  for all $\m\in{\cal M^*}$, we have $x_\m=x^*$.
Moreover, there exists an open sphere $S_{x^*}$ that is centered at $x^*$ and is such that for all $x\in S_{x^*}$ we have $T_{\m}x=T(x)$ only if $\m\in{\cal M^*}$.
}

\proof The relation $T_{\m}x^*=T(x^*)=x^*$ implies that $x^*$ is the unique fixed point $x_\m$ of $T_\m$. The remaining statement follows from the form of the components of $T$ as the minimum of a finite number of linear functions; cf.\ Eq.\ \minform\ and Assumption \assumptionlino(a). \qed

The preceding proposition illustrates the key idea of the algorithm \randomoptpolito-\mupdate, which is that for $\m\in{\cal M}^*$, the mappings $T_\m^{(\l)}$ are all  contractions with a common fixed point equal to $x^*$, the fixed point of $T$, and modulus
$${\r(1-\l)\over 1-\r\l}<\r$$
[cf.\ the definition \tddef\ of $T_\m^{(\l)}$]. Thus within the sphere $S_{x^*}$, the iterates \randomoptpolito\ aim consistently at $x^*$. Moreover, because of the randomization, the algorithm is guaranteed (with probability one) to eventually enter the sphere $S_{x^*}$, as shown in the following proposition. 

\xdef\propcontactionth{\propn}\propnum\show{myproposition}
 
%\texshopbox
{\proposition{\propcontactionth:} Under Assumption \assumptionlino, for any starting point $x_0$, a sequence $\{x_k\}$ generated by the algorithm \randomoptpolito-\mupdate\ converges to $x^*$ with probability one.
}

\proof We will show that $\{x_k\}$ is bounded by showing that for all $k$, we have
$$\max_{\m\in{\cal M}}\|x_k-x_\m\|\le \r^k\max_{\m\in{\cal M}}\|x_0-x_\m\|+{2\over 1-\r}\max_{\m,\m'\in{\cal M}}\|x_{\m}-x_{\m'}\|.\xdef\boundineq{\lab}\eqnum\show{oneo}$$
To this end we will consider separately the two cases where  $x_k=T_{\m^{k-1}}^{(\l)} x_{k-1}$ (probability $p$) and $x_k=T(x_{k-1})$ (probability $1-p$). In the former case, we have for all $\m\in{\cal M}$, 
$$\eqalign{\|x_k-x_{\m}\|&\le \|x_k-x_{\m^{k-1}}\|+\|x_{\m^{k-1}}-x_{\m}\|\cr
&=\|T_{\m^{k-1}}^{(\l)} x_{k-1}-x_{\m^{k-1}}\|+\|x_{\m^{k-1}}-x_\m\|\cr
&\le \r\|x_{k-1}-x_{\m^{k-1}}\|+\|x_{\m^{k-1}}-x_{\m}\|\cr
&\le \r\big(\|x_{k-1}-x_\m\|+\|x_\m-x_{\m^{k-1}}\|\big)+\|x_{\m^{k-1}}-x_{\m}\|\cr
&\le \r\max_{\m\in{\cal M}}\|x_{k-1}-x_\m\|+2\max_{\m,\m'\in{\cal M}}\|x_\m-x_{\m'}\|,\cr}$$
and finally, for all $k$,
$$\max_{\m\in{\cal M}}\|x_{k}-x_\m\|\le \r\max_{\m\in{\cal M}}\|x_{k-1}-x_\m\|+2\max_{\m,\m'\in{\cal M}}\|x_\m-x_{\m'}\|.\xdef\boundineqt{\lab}\eqnum\show{oneo}$$ 
A similar calculation shows that this inequality holds also in the complementary case where $x_k=T(x_{k-1})$.
By using the relation \boundineqt\ repeatedly, we obtain Eq.\ \boundineq.
Thus in conclusion, we have $\{x_k\}\subset D$, where $BD$ is the bounded set
$$D=\left\{x\ \Big|\ \max_{\m\in{\cal M}}\|x-x_\m\|\le \max_{\m\in{\cal M}}\|x_0-x_\m\|+{2\over 1-
\r}\max_{\m,\m'\in{\cal M}}\|x_{\m}-x_{\m'}\|\right\}.$$

\old{Consider for any $b>0$, the set
$${X}_b=\big\{x\mid \|x-x_\m\|\le b\ \hbox{for some } \m\in{\cal M} \big\}.$$
We note that ${X}_b$ is bounded in view of the finiteness of ${\cal M}$. Moreover, we have $T(x)\in X_b$ as well as $T_\m^{(\l)} x\in X_b$ for all $\m\in {\cal M}$ and $x\in X_b$. The reason is that $x^*$ is equal to all $x_\m$ with $\m\in{\cal M^*}$ and 
$$\big\|T(x)-x^*\big\|=\big\|T(x)-T(x^*)\big\|\le \r\|x-x^*\|\le \r b<b,$$
while
$$\big\|T_\m^{(\l)} x-x_\m\big\|=\big\|T_\m^{(\l)} x-T_\m^{(\l)} x_\m\big\|\le {\r(1-\l)\over 1-\r}\|x-x_\m\|\le {\r(1-\l)b\over 1-\r}<b,$$
[cf.\ the definition \tddef\ of $T_\m^{(\l)}$]. 
Thus for $b$ sufficiently large so that $x_0\in {X}_b$, the entire sequence $\{x_k\}$ remains within the bounded set ${X}_b$. 
}

Each time the iteration $x_{k+1}=T (x_k)$ is performed, the distance of the iterate $x_k$ to $x^*$ is reduced by a factor $\r$, i.e., $\|x_{k+1}-x^*\|\le \r \|x_k-x^*\|$. In view of our randomization scheme, the algorithm is guaranteed (with probability one) to eventually execute a sufficient number of contiguous iterations of the form $x_{k+1}=T (x_k)$ to enter the sphere $S_{x^*}$. Once this happens, all the subsequent iterations generate $\m_k$ within ${\cal M}^*$, and aim towards $x^*$ in the sense that they reduce the distance $\|x_k-x^*\|$ by a factor of at least $\r$. Thus the iterates $x_k$ will eventually enter and remain within  $S_{x^*}$, while $x_k\to x^*$. \qed

Note that the proof just given also applies to the algorithm \randomoptpolitz. 
The proof suggests that a potentially important issue is the choice of the randomization probability $p$ in the linearized algorithms \randomoptpolitz\ and \randomoptpolito. If $p$ is too large, convergence may be slow because oscillatory behavior may go unchecked for a long time. On the other hand if $p$ is small, a correspondingly large number of fixed point iterations $x_{k+1}=T(x_k)$ may be performed, and the hoped for benefits of the use of the proximal iterations may be lost. Adaptive schemes which adjust $p$ based on algorithmic progress may be an interesting possibility for addressing this issue.

Finally, let us note that the analysis of this subsection applies to the algorithms \randomoptpolitz\ and \randomoptpolito\ for other definitions of $T$, not necessarily involving minimization as in Eq.\ \minform. For example, the components $T(i,x)$ of $T(x)$ may be defined by
$$T(i,x)=\max_{\m(i)\in  M(i)}\big\{a\big(i,\m(i)\big)'x+b\big(i,\m(i)\big)\big\},\qquad x\in \rn,\ i=1,\ldots,n,\xdef\maxmap{\lab}\eqnum\show{oneo}$$
or by
$$T(i,x)=\min_{\m(i)\in  M(i)}\max_{\nu(i)\in  N(i)}\big\{a\big(i,\m(i),\nu(i)\big)'x+b\big(i,\m(i),\nu(i)\big)\big\},\qquad x\in \rn,\ i=1,\ldots,n,\xdef\minmaxmap{\lab}\eqnum\show{oneo}$$
where $N(i)$ is a given finite set for each $i$, and $a\big(i,\m(i),\nu(i)\big)$ and $b\big(i,\m(i),\nu(i)\big)$ are vectors and scalars, respectively, that depend on $i,\m(i),\nu(i)$. 

In the case of the mapping \maxmap, the mapping $T_\m$ is defined by $T_\m=A_\m x+b_\m$ [cf.\ Eq.\
\linsys], and we have 
$$T(x)=\max_{\m\in {\cal M}}T_\m x.$$
Moreover, $\m_k$ is updated by maximization, so the algorithm is
$$x_{k+1}=\cases{
T_{\m_k}^{(\l)} x_k& with probability $p$,\cr
T (x_k)& with probability $1-p$,\cr
}\xdef\randomoptpolith{\lab}\eqnum\show{oneo}$$
 where 
 $$\m_k\in\arg\max_{\m\in {\cal M}}T_\m x_k.$$
In the case of the mapping \minmaxmap, the mapping $T_\m$ has (nonlinear) components defined by
$$T_\m(i,x)=\max_{\nu(i)\in  N(i)}\big\{a\big(i,\m(i),\nu(i)\big)'x+b\big(i,\m(i),\nu(i)\big)\big\},\qquad x\in \rn,\ i=1,\ldots,n,$$
 and we have 
$$T(x)=\min_{\m\in {\cal M}}T_\m x.$$
The algorithm is again defined by Eq.\ \randomoptpolith\ with $\m_k$ updated according to 
 $$\m_k\in\arg\min_{\m\in {\cal M}}T_\m x_k.$$
In both cases the convergence result of Prop.\ \assumptionlino\ holds as stated. The mapping \minmaxmap\ arises among others in sequential games, which admit a DP treatment (see, e.g., the books [FiV96], [Ber12], or the survey [RaF91]).

Another interesting and related context arises when we want to solve a fixed point equation of the form
$$x=W T(x),$$
 where $W$ is an $n\times n$ matrix, and $T$ as earlier has the form
$$T(x)=\min_{\m\in{\cal M}}T_\m x,$$
where $T_\m$ is the linear mapping  $T_\m=A_\m x+b_\m$ [cf.\ Eq.\
\linsys]. An example is when $W$ is a projection on a subspace spanned by a set of basis functions [cf.\ Eq.\ \projoblique]. This context occurs in approximate DP, as discussed in Section 3, where the mappings $T_\m$ and $T$ are sup-norm contractions, but the composition mapping $W T$
need not be a contraction of any kind, even when $W$ is nonexpansive with respect to some norm. The reason is a potential ``norm mismatch" problem, whereby  $W$ is nonexpansive with respect to one norm, but $T_\m$ and $T$ are contractions with respect to a different norm (see de Farias and Van Roy [DFV00] for such examples). However, in the favorable case where all the mappings $(W T_\m)^{(\l)}$, $\m\in{\cal M}$, as well as the mapping $W T$ are contractions with respect to a common norm, the line of analysis of the present section applies to the algorithm
$$x_{k+1}=\cases{
(W T_{\m_k})^{(\l)} x_k& with probability $p$,\cr
W T(x_k)& with probability $1-p$,\cr
}\xdef\randomoptpolif{\lab}\eqnum\show{oneo}$$
 where 
 $$\m_k\in\arg\min_{\m\in {\cal M}}T_\m x_k.$$
In particular, assuming that the set ${\cal M}$ is finite, we can simply modify the proof of Prop.\ \propcontactionth\ to show that a sequence $\{x_k\}$ generated by this algorithm converges to the fixed point $x^*$ of $WT$ with probability one. For further discussion of this context, the associated convergence issues, and alternative choices of the matrix $W$, we refer to the author's survey [Ber11b] (Section 3.4) and the textbook [Ber12a]. 

\vskip  -4mm
\section{Concluding Remarks}
\vskip  -2mm

\pn In this paper we have shown that proximal and temporal difference methods for linear fixed point problems are closely related, and their implementations can benefit from each other, in both the exact and the approximate simulation-based setting. In particular, within the context of DP, the TD($\l$) algorithm for exact policy evaluation, can be written as the stochastic proximal algorithm 
$$x_{k+1}=x_k+\g_k\Big({sample}\big({\pb x_k}\big)-x_k\Big),$$
for solving the linear Bellman equation $x=Tx$ corresponding to a policy [in view of Eq.\ \tdlambda, and taking into account the fact that $({T^{(\l)} x_k- x_k})$ is the product of $({\pb x_k} -x_k)$ with the scalar $1/\l$, where $\l={c\over c+1}$; cf.\ Fig.\ \figoneone]. Of course, the probabilistic mechanism used to obtain the sample in the preceding equation is an important algorithmic issue, extensively researched in the approximate DP literature, but our analysis has not dealt with this (or any other issues that relate to stochastic proximal algorithms). Another important issue, which we have not fully explored, is how to connect the TD($\l$) algorithm for approximate policy evaluation and the projected proximal algorithm $x_{k+1}=\P {\pb x_k}$ of Section 3. 

Our Assumption \assumptiona\ is satisfied in broad classes of problems, including problems involving a contraction, and policy evaluation in exact and approximate DP. However, even when the assumption is not satisfied and $I-A$ is just invertible, the proximal equation 
$$x=\pb x=\Bigg({c+1\over c}I-A\Bigg)^{-1}\lf(b+{1\over c} x\ri),\xdef\pbformulat{\lab}\eqnum\show{oneo}$$
makes sense and has the same solution as the original equation $x=Ax+b$, provided the inverse in Eq.\ \pbformulat\ exists. In this case both the proximal equation $x=\pb x$ and its projected version $x=\P\pb x$ can be solved by simulation methods, which have been described in the papers [BeY07], [Yu12], [WaB13], [WaB14]; see also the textbook [Ber12a], Section 7.3.

Aside from the conceptual and analytical value of the connection between proximal and temporal difference methods, we have shown that under our assumptions, a tangible improvement of the proximal algorithm is possible at no cost. This improvement is obtained by a simple extrapolation of the proximal iterate, and provides a guaranteed acceleration of convergence (not just guaranteed convergence, like alternative extrapolation schemes). Moreover, this improvement carries over to nonlinear fixed point problems, as we have shown in Section 4. In addition, our methodology extends naturally to forward-backward splitting and proximal gradient algorithms.

To extend the connection between proximal and temporal difference algorithms, we have also introduced some proximal-like algorithms for nonlinear fixed point problems. These algorithms are based on linearization, bear a resemblance with Newton's method, and admit temporal differences-based implementations. They are inspired by distributed asynchronous policy iteration methods for infinite horizon DP, given in the papers [BeY10], [BeY12], [YuB13], and the monograph [Ber18], Section 2.6. Making a stronger connection of these methods and temporal difference methods [including TD($\l$)] for solving nonlinear Bellman equations exactly or approximately appears to be a fruitful direction of research. 

Some computational experience with the use of simulation to solve large linear systems, beyond those arising in DP, will be helpful in quantifying the potential benefits of the ideas of this paper. Another interesting question is how to generalize the methods of this paper from the linear equation context to the solution of linear variational inequalities, possibly with a large number of constraints, where both the proximal algorithm and multistep DP-type methods have been applied; see the papers [WaB15] and [Ber11a]. 

\mark{References}

\vskip  -4mm
\section{References}
\vskip  -2mm
\mark{References}

\def\ref{\vskip0pt\par\noindent}
\def\ref{\vskip0pt\par\noindent}

\ninepoint

\ref[BBD10] Busoniu, L., Babuska, R., De Schutter, B., and Ernst, D.,  2010.\ Reinforcement Learning and Dynamic Programming Using Function Approximators, CRC Press, NY.

\ref[BBN04] Bertsekas, D.\ P., Borkar, V.\ S., and Nedi\'c, A., 2004.\ ``Improved Temporal Difference Methods
with Linear Function Approximation," in Learning and Approximate Dynamic Programming, by J.\ Si, A.\ Barto,
W.\ Powell, and D.\ Wunsch (Eds.), IEEE Press, NY.

\ref[BDM14] Boutsidis, C., Drineas, P., and Magdon-Ismail, M., 2014.\ ``Near-Optimal Column-Based Matrix Reconstruction," SIAM J.\ on Computing, Vol.\ 43, pp.\ 687-717.

\ref[BaC11] Bauschke, H.\ H., and Combettes, P.\ L., 2011.\ 
Convex Analysis and Monotone Operator Theory in Hilbert Spaces, Springer, NY.

\ref[BeI96] Bertsekas, D.\ P., and Ioffe, S., 1996.\ ``Temporal Differences-Based Policy Iteration and Applications in Neuro-Dynamic Programming," Lab.\ for Info.\ and Decision Systems Report LIDS-P-2349, MIT.

\ref [BeS78]  Bertsekas, D.\ P., and Shreve, S.\ E., 1978.\  Stochastic Optimal
Control:  The Discrete Time Case, Academic Press, N.\ Y.; may be downloaded 
from http://web.mit.edu/dimitrib/www/home.html

\ref [BeT91]  Bertsekas, D.\ P., and Tsitsiklis, J.\ N., 1991.\ ``An Analysis of
Stochastic Shortest Path Problems,"
Math.\ of OR, Vol.\ 16, pp.\ 580-595.

\ref[BeT96] Bertsekas, D.\ P., and Tsitsiklis, J.\ N., 1996.\ Neuro-Dynamic Programming, Athena Scientific, Belmont, MA.

\ref[BeT10] Beck, A., and Teboulle, M., 2010.\ ``Gradient-Based Algorithms with Applications to Signal-Recovery Problems," in Convex
Optimization in Signal Processing and Communications (Y. Eldar and D. Palomar, eds.),
Cambridge Univ.\ Press, pp.\ 42-88.

\ref[BeY07]  Bertsekas, D.\ P., and Yu, H., 2007.\ ``Solution of Large Systems of Equations Using 
Approximate Dynamic Programming Methods," Lab.\ for Information and Decision Systems Report LIDS-P-2754, MIT.

\ref[BeY09] Bertsekas, D.\ P., and Yu, H., 2009.\ ``Projected Equation Methods for Approximate Solution of Large Linear Systems,"  Journal of Computational and Applied Mathematics, Vol.\ 227, pp.\ 27-50.

\ref[BeY10] Bertsekas, D.\ P., and Yu, H., 2010.\ ``Asynchronous Distributed Policy Iteration in Dynamic Programming,"  Proc.\ of Allerton Conf.\ on Communication, Control and Computing,  Allerton Park, Ill, pp.\ 1368-1374.

\ref[BeY12] Bertsekas, D.\ P., and Yu, H., 2012.\ ``Q-Learning and Enhanced Policy Iteration in Discounted 
Dynamic Programming,"  Math.\ of OR, Vol.\ 37, pp.\ 66-94.

\ref [Ber75] Bertsekas, D.\ P., 1975.\  ``On the Method of Multipliers for Convex Programming," IEEE Transactions on Aut.\ Control, Vol.\ 20, pp.\ 385-388.

\ref [Ber82] Bertsekas, D.\ P., 1982.\  Constrained Optimization and Lagrange
Multiplier Methods, Academic Press, NY; republished by Athena Scientific, Belmont,
MA, 1997; may be downloaded from\hfill\break  http://web.mit.edu/dimitrib/www/home.html

\ref[Ber11a] Bertsekas, D.\ P.,  2011.\ ``Temporal Difference Methods for General Projected Equations," IEEE Trans.\ on Automatic Control, Vol.\ 56, pp. 2128 - 2139.

\ref[Ber11b] Bertsekas, D.\ P., 2011.\
``Approximate Policy Iteration: A Survey and Some New Methods," J.\ of Control Theory and Applications, Vol.\ 9, 2011, pp.\ 310-335. 

\ref[Ber12a] Bertsekas, D.\ P., 2012.\ Dynamic Programming and Optimal Control: Approximate Dynamic Programming, 4th Edition, Vol.\  II, Athena
Scientific, Belmont, MA.

\ref[Ber12b] Bertsekas, D.\ P., 2012.\
``$\l$-Policy Iteration: A Review and a New Implementation," in Reinforcement Learning and Approximate Dynamic Programming for Feedback Control, by F.\ Lewis and D.\ Liu (eds.), IEEE Press, 2012.

\ref [Ber15] Bertsekas, D.\ P., 2015.\ Convex Optimization Algorithms, Athena Scientific, Belmont, MA.

\ref[Ber18] Bertsekas, D.\ P., 2018.\ Abstract Dynamic Programming, 2nd Edition, Athena Scientific, Belmont, MA; on-line at http://web.mit.edu/dimitrib/www/home.html.

\ref[Boy02]
Boyan, J.\ A., 2002.\ ``Technical Update: Least-Squares Temporal Difference Learning," Machine Learning, Vol.\ 49, pp.\ 1-15.

\ref[BrB96]
Bradtke, S.\ J., and Barto, A.\ G., 1996.\ ``Linear Least-Squares Algorithms for Temporal Difference Learning,'' Machine Learning, Vol.\ 22, pp.\ 33-57.

\ref[CeH09] Censor, J., Herman, G.\ T., and Jiang, M., 2009. ``A Note on the Behavior of the Randomized Kaczmarz Algorithm of Strohmer and Vershynin," J.\ Fourier Analysis and Applications, Vol.\ 15, pp.\ 431-436.

\ref[Cur54] Curtiss, J.\ H., 1954.\ ``A Theoretical Comparison of the Efficiencies of Two Classical Methods
and a Monte Carlo Method for Computing One Component of the Solution
of a Set of Linear Algebraic Equations," Proc.\ Symposium on Monte Carlo Methods, pp.\ 191-233.

\ref[Cur57] Curtiss, J.\ H., 1957.\ ``A Monte Carlo Methods for the Iteration of Linear Operators," Uspekhi Mat.\ Nauk, Vol.\ 12, pp.\ 149-174.

\ref[DFV00]  de Farias, D.\ P., and Van Roy, B., 2000.\ ``On the Existence of Fixed Points for Approximate 
Value Iteration and Temporal-Difference Learning," J.\ of  Optimization Theory and Applications, Vol.\
105, pp.\ 589-608.

\ref[DKM06a] Drineas, P., Kannan, R., and Mahoney, M.\ W., 2006. ``Fast Monte Carlo Algorithms for Matrices I: Approximating Matrix Multiplication," SIAM J.\ Computing, Vol.\ 35, pp.\ 132-157.

\ref[DKM06b] Drineas, P., Kannan, R., and Mahoney, M.\ W., 2006.\ ``Fast Monte Carlo Algorithms for Matrices II: Computing a Low-Rank Approximation to a Matrix," SIAM J.\ Computing, Vol.\ 36, pp.\ 158-183.

\ref[DMM06] Drineas, P., Mahoney, M.\ W., and Muthukrishnan, S., 2006. ``Sampling Algorithms for L2 Regression and Applications," Proc.\ 17th Annual SODA, pp.\ 1127-1136.

\ref[DMM08] Drineas, P., Mahoney, M.\ W., and Muthukrishnan, S., 2008.\ ``Relative-Error CUR Matrix Decompositions," SIAM J.\ Matrix Anal.\ Appl., Vol.\ 30, pp.\ 844-881.

\ref[DMM11] Drineas, P., Mahoney, M.\ W., Muthukrishnan, S., and Sarlos, T., 2011.\ ``Faster Least Squares Approximation," Numerische Mathematik, Vol.\ 117, pp.\ 219-249.

\ref [Den67] Denardo, E.\ V., 1967.\  ``Contraction Mappings in the Theory Underlying
Dynamic Programming," SIAM Review, Vol.\ 9, pp.\ 165-177.

\ref[DrL16] Drusvyatskiy, D., and Lewis, A.\ S., 2016.\ ``Error Bounds, Quadratic Growth, and Linear Convergence of Proximal Methods," arXiv:1602.06661.

\ref[DuR17] Duchi, J.\ C., and Ryan, F., 2017.\ ``Stochastic Methods for Composite Optimization Problems," arXiv:1703.08570.

\ref[EcB92] Eckstein, J., and Bertsekas, D.\ P., 1992.\ ``On the Douglas-Rachford Splitting
Method and the Proximal Point Algorithm for Maximal Monotone Operators," Math.\ Programming,
Vol.\ 55, pp.\ 293-318.

\ref[FaP03] Facchinei, F., and Pang, J.-S., 2003.\  Finite-Dimensional Variational Inequalities and
Complementarity Problems,  Springer Verlag, NY.

\ref [FiV96] Filar, J., and Vrieze, K., 1996.\ Competitive Markov Decision Processes, Springer,  N.\ Y.

\ref[Fle84] Fletcher, C.\ A.\ J., 1984.\ Computational Galerkin Methods, Springer-Verlag, NY.

\ref[FoL50]
Forsythe, G.\ E., and Leibler, R.\ A., 1950.\ ``Matrix Inversion by a Monte Carlo Method," Mathematical Tables and Other Aids to Computation, Vol.\ 4, pp.\ 127-129.

\ref[GMS13] Gabillon, V., Ghavamzadeh, M., and Scherrer, B., 2013.\ ``Approximate Dynamic Programming Finally Performs Well in the Game of Tetris," In Advances in Neural Information Processing Systems, pp.\ 1754-1762.

\ref[Gab83] Gabay, D., 1983.\ ``Applications of the Method of Multipliers to Variational Inequalities," in M.\ Fortin and R.\ Glowinski, eds., Augmented Lagrangian Methods: Applications to the Solution of Boundary-Value Problems, North-Holland, Amsterdam.

\ref[Hal70] Halton, J.\ H., 1970.\ ``A Retrospective and Prospective Survey of the Monte Carlo Method,"
SIAM Review, Vol.\ 12, pp.\ 1-63.

\ref[Kra72] Krasnoselskii, M.\ A., et. al, 1972.\ Approximate Solution of Operator Equations, Translated by D.\ Louvish, Wolters-Noordhoff Pub., Groningen.

\ref[LaP03] Lagoudakis, M.\ G.,  and Parr, R., 2003.\ ``Least-Squares Policy Iteration," J.\ of Machine Learning Research, Vol.\ 4, pp.\ 1107-1149.

\ref[LeL10]\ Leventhal, D., and Lewis, A.\ S., 2010.	 ``Randomized Methods for Linear Constraints: Convergence Rates and Conditioning,"
Mathematics of Operations Research,  Vol.\ 35, pp.\ 641-654.

\ref[LeL13] Lewis, F.\ L., and Liu, D., (Eds), 2013.\ Reinforcement Learning and Approximate Dynamic Programming for Feedback Control, Wiley, Hoboken, N.\ J.

\ref[LeW08] Lewis, A.\ S., and Wright, S.\ J., 2008.\ ``A Proximal Method for Composite Minimization," arXiv:0812.0423.

\ref[LiM79] Lions, P.\ L., and Mercier, B., 1979.\ ``Splitting Algorithms for the Sum of Two Nonlinear Operators," SIAM J.\ on Numerical Analysis, Vol.\ 16, pp.\ 964-979.

\ref[MKS15] Mnih, V., Kavukcuoglu, K., Silver, D., et al., 2015.\ ``Human-Level Control Through Deep Reinforcement Learning," Nature, Vol.\ 518, pp.\ 529-533.

{
\ref[Mar70]  Martinet, B., 1970.\ ``Regularisation d' Inequations Variationnelles
par Approximations Successives," Rev.\ Francaise Inf.\ Rech.\ Oper., Vol.\ 4,
pp.\ 154-158.
}

\ref[NeB03] Nedi\'c, A., and Bertsekas, D.\ P., 2003.\ ``Least Squares Policy Evaluation Algorithms with Linear Function Approximation," Discrete Event Dynamic Systems: Theory and Applications, Vol.\ 13, pp.\ 79-110.

\ref[PaB13] Parikh, N., and Boyd, S., 2013.\ 
 ``Proximal Algorithms," Foundations and Trends in Optimization, Vol.\ 1, pp.\ 123-231.

\ref [Pow07] Powell, W.\ B., 2007.\  Approximate Dynamic Programming: Solving the Curses of Dimensionality, Wiley, NY.

\ref [RaF91]  Raghavan, T.\ E.\ S., and Filar, J.\ A.,  1991.\ ``Algorithms for
Stochastic Games -- A Survey,'' ZOR -- Methods and Models of Operations
Research, Vol.\ 35, pp.\ 437-472.

\old{
\ref[Roc76]  Rockafellar, R. T., 1976.\ ``Monotone Operators and the Proximal
Point Algorithm," SIAM J.\ on Control and Optimization, Vol.\ 14, pp.\ 877-898.
}

\ref[SBP04] Si, J., Barto, A., Powell, W., and Wunsch, D., (Eds.), 2004.\ Learning and Approximate Dynamic Programming, IEEE Press, NY.

\ref[SHM16] Silver, D., Huang, A., Maddison, C.\ J., Guez, A., Sifre, L., Van Den Driessche, G., Schrittwieser, S., et al.\ 2016.\ ``Mastering the Game of Go with Deep Neural Networks and Tree Search," Nature, Vol.\ 529, pp.\ 484-489.

\ref[SMG15] Scherrer, B., Ghavamzadeh, M., Gabillon, V., Lesner, B., and Geist, M., 2015.\ ``Approximate Modified Policy Iteration and its Application to the Game of Tetris," J.\ of Machine Learning Research, Vol.\ 16, pp.\ 1629-1676.

\ref[Sam59] Samuel, A.\ L., 1959.\ ``Some Studies in Machine Learning
Using the Game of Checkers,''  IBM Journal of Research and Development, Vol.\ 3,
pp.\ 210-229.

\ref[Sam67] Samuel, A.\ L., 1967.\ ``Some Studies in Machine Learning
Using the Game of Checkers.\ II -- Recent Progress,'' 
IBM Journal of Research and Development, Vol.\ 11,
pp.\ 601-617.

\ref[Sch10] Scherrer, B., 2010.\ ``Should One Compute the Temporal Difference Fixed Point or Minimize the Bellman Residual? The Unified Oblique Projection View," Proc.\ of 2010 ICML, Haifa, Israel. 

\ref[Sch13] Scherrer, B., 2013.\ ``Performance Bounds for $\l$-Policy Iteration and Application to the Game of Tetris," J.\ of Machine Learning Research, Vol.\ 14, pp.\ 1181-1227.

\ref[StV09] Strohmer, T., and Vershynin, R., 2009.\ ``A Randomized Kaczmarz Algorithm with Exponential Convergence", J.\ of Fourier Analysis and Applications, Vol.\ 15, pp.\ 262-178.

\ref[SuB98] Sutton, R.\  S., and Barto, A.\ G., 1998.\ Reinforcement Learning, MIT
Press, Cambridge, MA.

\ref[Sut88] Sutton, R.\  S., 1988.\ ``Learning to Predict by the Methods of
Temporal Differences," Machine Learning, Vol.\ 3, pp.\ 9-44.

\ref [Sze10] Szepesvari, C., 2010.\ Algorithms for Reinforcement Learning, Morgan and Claypool Publishers.

\old{
\ref [Sze11] Szepesvari, C., 2011.\ ``Least Squares Temporal Difference Learning and Galerkin?s Method," Report, Univ. of Alberta.
}

\old{\ref[Tes92] Tesauro, G.\ J., 1992.\ ``Practical Issues in Temporal 
Difference Learning,"
Machine Learning, Vol.\ 8, pp.\ 257-277.}

\ref[TeG96] Tesauro, G., and Galperin, G.\ R., 1996.\ ``On-Line Policy Improvement
Using Monte Carlo Search,'' presented at the 1996 Neural Information
Processing Systems Conference, Denver, CO; also in M. Mozer et al.\ (eds.), Advances in
Neural Information Processing Systems 9, MIT Press (1997).

\ref[Tes94] Tesauro, G.\ J., 1994.\ ``TD-Gammon, a Self-Teaching 
Backgammon
Program, Achieves Master-Level Play,'' Neural Computation, Vol.\ 6, pp.\ 
215-219.

\ref [TsV97]
Tsitsiklis, J.\ N., and Van Roy, B., 1997.\ ``An Analysis of Temporal-Difference Learning
with Function Approximation," IEEE Transactions on Automatic Control, Vol.\ 42, pp.\ 674-690.

\ref[Tse91] Tseng, P., 1991.\ ``Applications of a Splitting Algorithm to Decomposition in Convex Programming and Variational Inequalities," SIAM J.\ on Control and Optimization, Vol.\ 29, pp.\ 119-138.

\ref[VVL13] Vrabie, D., Vamvoudakis, K.\ G., and Lewis, F.\ L., 2013.\ Optimal Adaptive Control and Differential Games by Reinforcement Learning Principles,
The Institution of Engineering and Technology, London.

\ref[WaB13] Wang, M., and Bertsekas, D.\ P., 2013.\ ``Stabilization of Stochastic Iterative Methods for Singular and Nearly Singular Linear Systems," Math.\ of Operations Research, Vol.\ 39, pp.\ 1-30.

\ref[WaB14] Wang, M., and Bertsekas, D.\ P., 2014.\ ``Convergence of Iterative Simulation-Based Methods for Singular Linear Systems," Stochastic Systems, Vol.\ 3, pp.\ 39-96.

\ref[WaB15] Wang, M., and Bertsekas, D.\ P., 2015.\ ``Incremental Constraint Projection Methods for Variational Inequalities," Math.\ Programming, Vol.\ 150, pp.\ 321-363.

\ref[WiB93] Williams, R.\ J., and Baird, L.\ C., 1993.\ ``Analysis of Some Incremental Variants of Policy Iteration: First Steps Toward Understanding Actor-Critic Learning Systems,'' Report NU-CCS-93-11, College of Computer Science, Northeastern University, Boston, MA. 

\old{
\ref[Was52] Wasow, W.\ R., 1952.\ ``A Note on Inversion of Matrices by Random Walks," Mathematical Tables and Other Aids to Computation, Vol.\ 6, pp.\ 78-81.
}

\ref[YuB06] Yu, H., and Bertsekas, D.\ P., 2006.\ ``Convergence Results for Some Temporal Difference Methods Based on Least Squares," IEEE Trans.\ on Aut. Control, Vol.\ 54, pp.\ 1515-1531.

\ref[YuB10] Yu, H., and Bertsekas, D.\ P., 2010.\ ``Error Bounds for Approximations from Projected Linear Equations," Math.\ of Operations Research, Vol.\ 35, pp.\ 306-329. 

\ref[YuB12] Yu, H., and Bertsekas, D.\ P., 2012.\ ``Weighted Bellman Equations and their Applications in Dynamic Programming," Lab.\ for Information and Decision Systems Report LIDS-P-2876, MIT. 

\ref[YuB13] Yu, H., and Bertsekas, D.\ P., 2013.\ ``Q-Learning and Policy Iteration Algorithms for Stochastic Shortest Path Problems," Annals of Operations Research, Vol.\ 208, pp.\ 95-132.

\ref[Yu10] Yu, H., 2010.\ ``Convergence of Least Squares Temporal Difference Methods
Under General Conditions," Proc.\ of the 27th ICML, Haifa, Israel.

\ref[Yu12] Yu, H., 2012.\ ``Least Squares Temporal Difference Methods: An Analysis
Under General Conditions," SIAM J.\ on Control and Optimization, Vol.\ 50, pp.\ 3310-3343.

\end